\documentclass[a4paper, british]{amsart}

\usepackage{libertine}
\usepackage[libertine]{newtxmath}

\usepackage{babel}
\usepackage{enumitem}
\usepackage{hyperref}
\usepackage[utf8]{inputenc}
\usepackage{newunicodechar}
\usepackage{mathtools}
\usepackage{varioref}
\usepackage[arrow,curve,matrix]{xy}

\usepackage{colortbl}
\usepackage{graphicx}
\usepackage{tikz}

\usepackage{mathrsfs}

\definecolor{linkred}{rgb}{0.7,0.2,0.2}
\definecolor{linkblue}{rgb}{0,0.2,0.6}

\numberwithin{figure}{section}

\usepackage[hyperpageref]{backref}

\setdescription{labelindent=\parindent, leftmargin=2\parindent}
\setitemize[1]{labelindent=\parindent, leftmargin=2\parindent}
\setenumerate[1]{labelindent=0cm, leftmargin=*, widest=iiii}

\newunicodechar{α}{\ensuremath{\alpha}}
\newunicodechar{β}{\ensuremath{\beta}}
\newunicodechar{χ}{\ensuremath{\chi}}
\newunicodechar{δ}{\ensuremath{\delta}}
\newunicodechar{ε}{\ensuremath{\varepsilon}}
\newunicodechar{Δ}{\ensuremath{\Delta}}
\newunicodechar{η}{\ensuremath{\eta}}
\newunicodechar{γ}{\ensuremath{\gamma}}
\newunicodechar{Γ}{\ensuremath{\Gamma}}
\newunicodechar{ι}{\ensuremath{\iota}}
\newunicodechar{κ}{\ensuremath{\kappa}}
\newunicodechar{λ}{\ensuremath{\lambda}}
\newunicodechar{Λ}{\ensuremath{\Lambda}}
\newunicodechar{ν}{\ensuremath{\nu}}
\newunicodechar{μ}{\ensuremath{\mu}}
\newunicodechar{ω}{\ensuremath{\omega}}
\newunicodechar{Ω}{\ensuremath{\Omega}}
\newunicodechar{π}{\ensuremath{\pi}}
\newunicodechar{Π}{\ensuremath{\Pi}}
\newunicodechar{φ}{\ensuremath{\phi}}
\newunicodechar{Φ}{\ensuremath{\Phi}}
\newunicodechar{ψ}{\ensuremath{\psi}}
\newunicodechar{Ψ}{\ensuremath{\Psi}}
\newunicodechar{ρ}{\ensuremath{\rho}}
\newunicodechar{σ}{\ensuremath{\sigma}}
\newunicodechar{Σ}{\ensuremath{\Sigma}}
\newunicodechar{τ}{\ensuremath{\tau}}
\newunicodechar{θ}{\ensuremath{\theta}}
\newunicodechar{Θ}{\ensuremath{\Theta}}
\newunicodechar{ξ}{\ensuremath{\xi}}
\newunicodechar{Ξ}{\ensuremath{\Xi}}
\newunicodechar{ζ}{\ensuremath{\zeta}}

\newunicodechar{ℓ}{\ensuremath{\ell}}
\newunicodechar{ï}{\"{\i}}

\newunicodechar{⁰}{\ensuremath{^0}}
\newunicodechar{¹}{\ensuremath{^1}}
\newunicodechar{²}{\ensuremath{^2}}
\newunicodechar{³}{\ensuremath{^3}}
\newunicodechar{⁴}{\ensuremath{^4}}
\newunicodechar{⁵}{\ensuremath{^5}}
\newunicodechar{⁶}{\ensuremath{^6}}
\newunicodechar{⁷}{\ensuremath{^7}}
\newunicodechar{⁸}{\ensuremath{^8}}
\newunicodechar{⁹}{\ensuremath{^9}}
\newunicodechar{ⁱ}{\ensuremath{^i}}

\newunicodechar{⌈}{\ensuremath{\lceil}}
\newunicodechar{⌉}{\ensuremath{\rceil}}
\newunicodechar{⌊}{\ensuremath{\lfloor}}
\newunicodechar{⌋}{\ensuremath{\rfloor}}

\newunicodechar{≅}{\ensuremath{\cong}}
\newunicodechar{⇔}{\ensuremath{\Leftrightarrow}}
\newunicodechar{∃}{\ensuremath{\exists}}
\newunicodechar{±}{\ensuremath{\pm}}

\DeclareFontFamily{OMS}{rsfs}{\skewchar\font'60}
\DeclareFontShape{OMS}{rsfs}{m}{n}{<-5>rsfs5 <5-7>rsfs7 <7->rsfs10 }{}
\DeclareSymbolFont{rsfs}{OMS}{rsfs}{m}{n}
\DeclareSymbolFontAlphabet{\scr}{rsfs}
\DeclareSymbolFontAlphabet{\scr}{rsfs}

\DeclareFontFamily{U}{mathx}{\hyphenchar\font45}
\DeclareFontShape{U}{mathx}{m}{n}{
      <5> <6> <7> <8> <9> <10>
      <10.95> <12> <14.4> <17.28> <20.74> <24.88>
      mathx10
      }{}
\DeclareSymbolFont{mathx}{U}{mathx}{m}{n}
\DeclareFontSubstitution{U}{mathx}{m}{n}
\DeclareMathAccent{\wcheck}{0}{mathx}{"71}

\newunicodechar{∂}{\ensuremath{\partial}}
\newunicodechar{∇}{\ensuremath{\nabla}}

\newunicodechar{↺}{\ensuremath{\circlearrowleft}}
\newunicodechar{∞}{\ensuremath{\infty}}
\newunicodechar{⊕}{\ensuremath{\oplus}}
\newunicodechar{⊗}{\ensuremath{\otimes}}
\newunicodechar{•}{\ensuremath{\bullet}}
\newunicodechar{Λ}{\ensuremath{\wedge}}
\newunicodechar{↪}{\ensuremath{\into}}
\newunicodechar{→}{\ensuremath{\to}}
\newunicodechar{↦}{\ensuremath{\mapsto}}
\newunicodechar{⨯}{\ensuremath{\times}}
\newunicodechar{∪}{\ensuremath{\cup}}
\newunicodechar{∩}{\ensuremath{\cap}}
\newunicodechar{⊋}{\ensuremath{\supsetneq}}
\newunicodechar{⊇}{\ensuremath{\supseteq}}
\newunicodechar{⊃}{\ensuremath{\supset}}
\newunicodechar{⊊}{\ensuremath{\subsetneq}}
\newunicodechar{⊆}{\ensuremath{\subseteq}}
\newunicodechar{⊂}{\ensuremath{\subset}}
\newunicodechar{⊄}{\ensuremath{\not \subset}}
\newunicodechar{≥}{\ensuremath{\geq}}
\newunicodechar{≠}{\ensuremath{\neq}}
\newunicodechar{≫}{\ensuremath{\gg}}
\newunicodechar{≪}{\ensuremath{\ll}}

\newunicodechar{≤}{\ensuremath{\leq}}
\newunicodechar{∈}{\ensuremath{\in}}
\newunicodechar{∉}{\ensuremath{\not \in}}
\newunicodechar{∖}{\ensuremath{\setminus}}
\newunicodechar{◦}{\ensuremath{\circ}}
\newunicodechar{°}{\ensuremath{^\circ}}
\newunicodechar{…}{\ifmmode\mathellipsis\else\textellipsis\fi}
\newunicodechar{·}{\ensuremath{\cdot}}
\newunicodechar{⋯}{\ensuremath{\cdots}}
\newunicodechar{∅}{\ensuremath{\emptyset}}
\newunicodechar{⇒}{\ensuremath{\Rightarrow}}

\newunicodechar{⁰}{\ensuremath{^0}}
\newunicodechar{¹}{\ensuremath{^1}}
\newunicodechar{²}{\ensuremath{^2}}
\newunicodechar{³}{\ensuremath{^3}}
\newunicodechar{⁴}{\ensuremath{^4}}
\newunicodechar{⁵}{\ensuremath{^5}}
\newunicodechar{⁶}{\ensuremath{^6}}
\newunicodechar{⁷}{\ensuremath{^7}}
\newunicodechar{⁸}{\ensuremath{^8}}
\newunicodechar{⁹}{\ensuremath{^9}}
\newunicodechar{ⁱ}{\ensuremath{^i}}

\newunicodechar{⌈}{\ensuremath{\lceil}}
\newunicodechar{⌉}{\ensuremath{\rceil}}
\newunicodechar{⌊}{\ensuremath{\lfloor}}
\newunicodechar{⌋}{\ensuremath{\rfloor}}

\newunicodechar{≅}{\ensuremath{\cong}}
\newunicodechar{⇔}{\ensuremath{\Leftrightarrow}}
\newunicodechar{∃}{\ensuremath{\exists}}
\newunicodechar{±}{\ensuremath{\pm}}

\DeclareFontFamily{OMS}{rsfs}{\skewchar\font'60}
\DeclareFontShape{OMS}{rsfs}{m}{n}{<-5>rsfs5 <5-7>rsfs7 <7->rsfs10 }{}
\DeclareSymbolFont{rsfs}{OMS}{rsfs}{m}{n}
\DeclareSymbolFontAlphabet{\scr}{rsfs}
\DeclareSymbolFontAlphabet{\scr}{rsfs}

\DeclareFontFamily{U}{mathx}{\hyphenchar\font45}
\DeclareFontShape{U}{mathx}{m}{n}{
      <5> <6> <7> <8> <9> <10>
      <10.95> <12> <14.4> <17.28> <20.74> <24.88>
      mathx10
      }{}
\DeclareSymbolFont{mathx}{U}{mathx}{m}{n}
\DeclareFontSubstitution{U}{mathx}{m}{n}
\DeclareMathAccent{\wcheck}{0}{mathx}{"71}

\DeclareMathOperator{\sing}{sing}

\theoremstyle{plain}
\newtheorem{Th}{Théorème}[section]

\newtheorem{cor}[Th]{Corollaire}
\newtheorem{defn}[Th]{Définition}

\newtheorem{lem}[Th]{Lemme}

\newtheorem{Prop}[Th]{Proposition}

\theoremstyle{remark}

\newtheorem{c-n-d}[Th]{Claim and Definition}

\newtheorem{notation}[Th]{Notation}
\newtheorem{obs}[Th]{Observation}
\newtheorem{rem}[Th]{\color{blue}{Remarques}}

\newtheorem*{rem-nonumber}{Remark}

\numberwithin{equation}{Th}
\setlist[enumerate]{label=(\thethm.\arabic*), before={\setcounter{enumi}{\value{equation}}}, after={\setcounter{equation}{\value{enumi}}}}
\newcommand{\into}{\hookrightarrow}

\newcommand{\factor}[2]{\left. \raise 2pt\hbox{$#1$} \right/\hskip -2pt\raise -2pt\hbox{$#2$}}

\newcommand{\Publication}[1]{}
\newcommand{\approvals}[2][Approval]{}
\renewcommand{\phi}{\varphi}
\usepackage{tikz}
\usepackage{tikz-cd}
\tikzset{commutative diagrams/arrow style=Latin Modern}
\author{Mohamed Kaddar} %
\email{\href{mailto:mohamed.kaddar@univ-lorraine.fr}{mohamed.kaddar@univ-lorraine.fr}} %
\title[Dualité relative de type Kleiman: Cas propre.]{Dualité relative de type Kleiman I: Cas propre.}%
\date{\today}

\makeatletter
\hypersetup{
  pdfauthor={\authors},
  pdftitle={\@title},
  pdfsubject={\@subjclass},
  pdfkeywords={\@keywords},
  pdfstartview={Fit},
  pdfpagelayout={TwoColumnRight},
  pdfpagemode={UseOutlines},
  bookmarks,
  colorlinks,
  linkcolor=linkblue,
  citecolor=linkred,
  urlcolor=linkred}
\makeatother
\newcommand{\chapref}[1]{\hyperref[#1]{Chapter~\ref*{#1}}}
\newcommand{\lemmaref}[1]{\hyperref[#1]{Lemma~\ref*{#1}}}
\newcommand{\parref}[1]{\hyperref[#1]{Section~\ref*{#1}}}
\newcommand{\theoremref}[1]{\hyperref[#1]{Théorème~\ref*{#1}}}
\newcommand{\definitionref}[1]{\hyperref[#1]{Definition~\ref*{#1}}}
\newcommand{\propositionref}[1]{\hyperref[#1]{Proposition~\ref*{#1}}}
\newcommand{\conjectureref}[1]{\hyperref[#1]{Conjecture~\ref*{#1}}}
\newcommand{\corollaryref}[1]{\hyperref[#1]{Corollary~\ref*{#1}}}
\newcommand{\exampleref}[1]{\hyperref[#1]{Example~\ref*{#1}}}
\newcommand{\exerciseref}[1]{\hyperref[#1]{Exercise~\ref*{#1}}}
\newcommand{\factref}[1]{\hyperref[#1]{Fact~\ref*{#1}}}
\newcommand{\claimref}[1]{\hyperref[#1]{Claim~\ref*{#1}}}
\newcommand{\remarkref}[1]{\hyperref[#1]{Remark~\ref*{#1}}}
\newcommand{\settingref}[1]{\hyperref[#1]{Setting~\ref*{#1}}}
\newcommand{\appendixref}[1]{\hyperref[#1]{Appendix~\ref*{#1}}}

\theoremstyle{plain}

\theoremstyle{remark}

\newcommand{\eq}[1][r]
   {\ar@<-3pt>@{-}[#1]
    \ar@<-1pt>@{}[#1]|<{}="gauche"
    \ar@<+0pt>@{}[#1]|-{}="milieu"
    \ar@<+1pt>@{}[#1]|>{}="droite"
    \ar@/^2pt/@{-}"gauche";"milieu"
    \ar@/_2pt/@{-}"milieu";"droite"}

\begin{document}

\maketitle
\approvals[Approval for Abstract]{Mohamed & yes}
\begin{abstract} Dans cette partie, nous montrons que:\vspace{1mm}

\noindent {\bf(1)}
 si $X$ et $S$ sont des espaces complexes de dimension finie avec $X$ dénombrable à l'infini et $\pi:X\rightarrow S$ un morphisme propre à fibres de dimension $n$,  le foncteur  ${\rm I}\!{\rm R}^{n}\pi_{*}$ admet un adjoint à droite  $\pi^{!}_{\mathcal K}$ induisant une dualité relative au niveau de la catégorie des faisceaux cohérents similaire à celle  décrite par  Kleiman (\cite{Kl}) dans le cadre algébrique et jouissant de toutes les propriétés que satisfait son analogue algébrique,\vspace{1mm}
 
 \noindent
 {\bf(2)} des relations mises en évidence entre la dualité relative générale pour un morphisme propre de Ramis-Ruget-Verdier pour le foncteur ${\rm I}\!{\rm R}\pi_{*}$ \cite{RRV71} et celle introduite par Flenner pour le foncteur ${\rm I}\!{\rm L}\pi_{*}$ \cite{Fle81},  permettent de dégager une condition nécessaire et suffisante garantissant l'identification de ces foncteurs  à un décalage près, 
 \vspace{1mm}
 
 \noindent
 {\bf(3)} Une étude systématique des foncteurs cohomologiques ${\mathcal H}^{-n}(\pi^{!}({\mathcal G})$ et ${\mathcal H}^{0}({\rm I}\!{\rm L}\pi^{*}({\mathcal G})$
 nous conduit à la conclusion de \cite{Kl} disant que la dualité totale n'est réalisable que si et seulement si $\pi$ est de Cohen Macaulay.
\end{abstract}
\maketitle
\tableofcontents

%
%
\section{\color{blue}{Introduction}.}
Le présent  travail s'insère dans une étude générale concernant l'extension au cadre analytique complexe de la dualité de Kleiman \cite{Kl} assurant que, 
 pour tout morphisme propre $\pi:X\rightarrow S$  de présentation
finie de schémas  dont les fibres sont de dimension pure
 au plus égale à un certain entier $n$, le foncteur
  ${\rm I}\!{\rm R}^{n}\pi_{*}:{\rm Q}coh(X) \rightarrow {\rm Q}coh(S)$ défini  de la
catégorie des faisceaux quasi-cohérents sur $X$ vers celle des
quasi-cohérents sur $S$  admet un foncteur  adjoint à droite
covariant
 ( noté $\pi^{!}_{\mathcal K}$)  induisant un isomorphisme de dualité relative
$${\rm I}\!{\rm H}om(X; {\mathcal F},\pi^{!}_{\mathcal K}({\mathcal G}))\simeq
{\rm I}\!{\rm H}om(S; {\rm I}\!{\rm R}^{n}\pi_{*}{\mathcal F}, {\mathcal G})$$
  bifonctoriel en  ${\mathcal F}$ et ${\mathcal G}$ et vérifiant un certain nombre de propriétés fonctorielles. Il montre, de plus,
     que pour $\pi$ plat et projectif, le morphisme canonique appelé de {\emph{dualité totale}}
     $${\rm E}xt^{m}(\pi; {\mathcal F}, \pi^{*}{\mathcal G}\otimes \omega^{n}_{X/S})\rightarrow {\rm I}\!{\rm H}om(S;{\rm I}\!{\rm R}^{n-m}\pi_{*}{\mathcal F}, {\mathcal G})$$
est un isomorphisme si et seulement si $\pi$ est de Cohen-Macaulay (i.e plat et à fibres non vides de Cohen-Macaulay). \vspace{1mm}

\noindent

L'approche de Kleiman consiste à rester au niveau de la catégorie des faisceaux quasi-cohérents  et montrer l'existence d'un  adjoint à droite du foncteur ${\rm I}\!{\rm R}^{n}\pi_{*}$.  Pour ce faire, il s'appuie sur des  considérations générales et les critères donnant
l'existence d'un  adjoint  d'un foncteur donné en se réferrant
à  (\cite{Ml}, Thm 2; p.125).  Sa construction est propre à la géométrie algébrique et ne se transpose pas telle quelle en géométrie
analytique complexe car contrairement au cadre  algébrique où tout schéma noethérien peut être reconstruit à l'aide de sa catégorie
des faisceaux cohérents, le cadre analytique complexe général
ne possède pas cette propriété puisque la catégorie
considérée  n'admet pas de famille compacte de générateurs
rendant de fait la référence (\cite{Ml}, Thm 2, p.125)
inappropriée.  D'ailleurs, il est  bien connu que les variétés complexes générales  possèdent trop peu de fibrés vectoriels
holomorphes pour pouvoir contrôler un faisceau cohérent
arbitraire sur ladite variété (le lecteur peut regarder \cite{Vo} où est donné l'exemple d'un  faisceau d'idéaux d'un point sur un tore générique complexe de dimension au
moins $3$ n'admettant pas de résolution localement libre globale).
\vspace{1mm}

\noindent 
Pour pallier le manque d'objet compact, on pourrait penser,  si le morphisme est propre, à la notion de  catégorie {\emph bien} {\emph{ engendrée}} (i.e admettant un nombre arbitrairement grand de
générateurs) introduite par Neeman dans \cite{Ne2} (cf le  théorème de représentabilité de Brown-Neeman  (\cite{Ne1},  {{\emph theorem (4.1)}})  étendant les critères de représentabilité classiques. Malheureusement,  bien que le foncteur ${\rm I}\!{\rm R}^{n}\pi_{*}$ soit  continue (i.e il commute à la limite inductive), la limite inductive de faisceaux cohérents n'étant généralement jamais cohérente, on se retrouve naturellement devant deux gros problèmes à savoir décrire ${\rm IndCoh}(X)$ et définir une notion de quasi-cohérence dans le cadre analytique. Or, jusqu'à présent, aucune définition raisonnable et générale n'a pu être donnée sachant que ces objets "idéaux" devraient   satisfaire les théorèmes {\bf A} et {\bf B} de Cartan tout en recouvrant le cas cohérent (sans conditions supplémentaires,les limites inductives, même complétées,  de faisceaux cohérents ne sont, généralement,  jamais acycliques sur les ouverts de Stein).\vspace{1mm}

\noindent 
Ce travail se divise en deux parties à savoir le cas d'un morphisme propre et le cas général pour lequel il nous est nécessaire de le développer dans le cadre de la dualité des espaces vectoriels topologiques localement convexes.\vspace{1mm}

\noindent
En ce qui concerne, le cas propre, nous déduisons cette dualité de type Kleiman de la dualité relative pour un morphisme propre (\cite{RRV71}) exprimant que, pour tout complexe ${\mathcal A}^{\bullet}$  de  ${\mathcal O}_{X}$-modules à cohomologie cohérente (${\mathcal A}^{\bullet}\in {\bf  D}_{{\rm coh}}(X)$ ) et tout complexe ${\mathcal B}^{\bullet}$ de ${\mathcal O}_{X}$-modules à cohomologie cohérente et bornée à gauche ( ${\mathcal B}^{\bullet}\in {\bf  D}^{+}_{{\rm coh}}(S)$), la flèche canonique 
$${\rm I}\!{\rm
 R}{\pi}_{*}{\rm I}\!{\rm R}{\rm H}om(X; {\mathcal A}^{\bullet}, \overset{\pi^!}{\overbrace{{\rm I}\!{\rm R}{\rm H}om( {\rm I}\!{\rm
 L}\pi^{*}({\Bbb D}_{S}({\mathcal B}^{\bullet})), {\mathcal D}^{\bullet}_{X})}})\rightarrow
{\rm I}\!{\rm R}{\rm H}om(S; {\rm I}\!{\rm R}\pi_{*}{\mathcal A}^{\bullet}, {\mathcal B}^{\bullet})$$ est un isomorphisme.\vspace{1mm}

\noindent
Dans ce cas, le foncteur $\pi^{!}_{\mathcal K}$ est, par définition ou construction, est lié au foncteur $\pi^{!}$  par la relation $\pi^{!}_{\mathcal K}({\mathcal G})={\mathcal H}^{-n}(\pi^{!}({\mathcal G}))$. Il nous reste, alors, à établir les différentes propriétés fonctorielles et étudier soigneusement celles du faisceau $\pi^{!}_{\mathcal K}({\mathcal G})$ en vue de lui trouver un analogue possédant le plus grand nombre de ses propriétés.
\vspace{1mm}
  
  \noindent
  Ayant pour second objectif l'étude du morphisme de dualité totale, on se place dans la situation d'un morphisme propre  ayant un faisceau canonique relatif $\omega^{n}_{X/S}:={\mathcal H}^{-n}(\pi^{!}({\mathcal O}_{S}))$ plat sur $S$. On fait, alors référence à la dualité de Flenner \cite{Fle81} ({\bf Satz (2.1),p.178}) disant que, pour tout morphisme $\pi:X\rightarrow S$  d'espaces complexes et tout faisceau cohérent ${\mathcal G}$ $X$, $S$-plat et de support propre sur $S$, il existe un foncteur covariant  ${\rm I}\!{\rm L}\pi_{\#}: {\rm D}^{-}_{coh}(X)\rightarrow  {\rm D}^{-}_{coh}(S)$ dépendant de ${\mathcal G}$, commutant aux changements de base arbitraires et induisant l'isomorphisme
  fonctoriel en les arguments,
$${\rm I}\!{\rm R}\pi_{*}{\rm I}\!{\rm R}{\mathcal H}om({\mathcal F}^{\bullet},  \overset{{\rm I}\!{\rm L}\pi^{\#}}{\overbrace{{\mathcal G}\otimes{\rm I}\!{\rm L}\pi^{*}{\mathcal N}^{\bullet}}})\simeq{\rm I}\!{\rm R}{\mathcal H}om({\rm I}\!{\rm L}\pi_{\#}({\mathcal F}^{\bullet}), {\mathcal N}^{\bullet})$$ 
pour tout ${\mathcal F}^{\bullet}\in  {\rm D}^{-}_{coh}(X)$ et tout ${\mathcal N}^{\bullet}\in  {\rm D}^{+}_{coh}(S)$.

\vspace{1mm}
On établi facilement l'existence d'un morphisme de foncteurs
${\rm I}\!{\rm R}\pi_{*}[n]\rightarrow {\rm I}\!{\rm L}\pi_{\#}$ dont on montre que c'est un isomorphisme si et seulement si $\pi$ est de Cohen Macaulay.\vspace{1mm}

\noindent
La paire $({\rm I}\!{\rm R}^{n}\pi_{*}, \pi^{!}_{{\mathcal K}})$ (resp. $({\mathcal H}^{0}({\rm I}\!{\rm L}\pi_{\#}), {\mathcal H}^{0}({\rm I}\!{\rm L}\pi^{\#})) $ ) hérite naturellement sa dualité de la paire $({\rm I}\!{\rm R}\pi_{*}, \pi^{!})$ (resp. $({\rm I}\!{\rm L}\pi_{\#}, {\rm I}\!{\rm L}\pi^{\#})$. On arrive, alors, à la conclusion que pour tout  morphisme propre  d'espaces complexes de dimension finie $\pi:X\rightarrow S$ à fibres de dimension pure $n$ et $X$ dénombrable à l'infini, le foncteur ${\rm I}\!{\rm R}^{n}\pi_{*}$ admet un adjoint à droite $\pi^{!}_{\mathcal K}$ induisant, pour tout  faisceaux  cohérents   ${\mathcal
F}$ et ${\mathcal G}$  sur $X$ et $S$ respectivement, un isomorphisme de
dualité relative
$${\rm I}\!{\rm H}om(X; {\mathcal F},\pi^{!}_{\mathcal K}{\mathcal G})\simeq
{\rm I}\!{\rm H}om(S; {\rm I}\!{\rm R}^{n}\pi_{*}{\mathcal
F}, {\mathcal G})$$
 bifonctoriel en ces faisceaux et jouissant de
propriétés fonctorielles analogues à celle du cas
algébrique.
\vspace{1mm}

\noindent
 On termine, alors, en montrant que, si $\pi^{!}_{{\mathcal K}}$ est exacte à droite. Le morphisme de dualité totale est un isomorphisme si et seulement si $\pi$ est de Cohen Macaulay. 
 \vspace{2mm}

\noindent 
  Dans le cas non propre (sur lequel on ne s'attarde pas trop dans cette introduction car c'est  l'objet de la seconde partie),  on se doit d'étudier le foncteur image directe supérieure à support propre ${\rm I}\!{\rm R}^{n}\pi_{!}$. Mais à de rares exceptions près imposées par des contraintes drastiques sur le faisceau cohérent ou sur le morphisme lui-même,  il est bien connu que ces images directes de faisceaux cohérents ne sont  jamais cohérentes. Comme on ne peut asseoir notre dualité sur une théorie générale déjà existante; on rappelle, au passage, que  la dualité relative de Ramis-Ruget \cite{RR74} est encore, à bien des égards,  sous une forme inachevée et ne permet en aucun d'aborder cette question sous sa forme générale. On obtient un résultat de pseudo-dualité dans le sens suivant: \vspace{1mm}

  \noindent si ${\rm Qcoh}_{R}(S)$ est la catégorie quasi-abélienne des faisceaux quasi-cohérents au sens de Ramis-Ruget\footnote{Insistons sur le fait que, dans ce contexte ``simpliciel'', il se dégage, tout naturellement, une certaine  notion de  quasi-cohérence  mais dont la
définition reste encore plus ou moins floue dans le sens où elle
dépend  de la situation et de la problématique considérée.
A ce titre, rappelons qu'elle fût  introduite, pour la première,
dans la dualité relative de Ramis et Ruget \cite{RR74}, étudiée sous
l'aspect topologique par J.L Stehlé \cite{St},  reprise, dans un cadre
"catégoriel" par Orlov \cite{O}, défini, pour les besoins de la
géométrie analytique rigide, par Conrad [C] ou dans le cadre de
la transformation de Mukai par Ben-Bassat, Block et Pantev
\cite{BBBP}  en géométrie analytique complexe proche de celle de la géométrie algébrique et tenant compte de la topologie}, on montre qu'à tout morphisme $\pi:X\rightarrow S$ universellement $n$-équidimensionnel (i.e ouvert à fibres de dimension constante $n$) est associé un couple de foncteurs  ${\rm I}\!{\rm R}^{n}\pi_{!}: {\rm Coh}(X)\rightarrow {\rm Qcoh_R}(S)$ et $\overline{\pi}:{\rm Coh}(S)\rightarrow {\rm Coh}(X)$ tel que le morphisme canonique
$${\rm I}\!{\rm H}om(X; {\mathcal F}, \overline{\pi}{\mathcal G})\simeq
{\rm I}\!{\rm H}om(S; {\rm I}\!{\rm R}^{n}\pi_{!}{\mathcal
F}, {\mathcal G})$$
 soit  un isomorphisme bifonctoriel en les faisceaux cohérents ${\mathcal F}$ et ${\mathcal G}$ et jouissant de
propriétés fonctorielles comparables à celles de son analogue propre. Il est bien évident que l'on doit agrémenter de topologies adéquates ces objets car ils sont loin d'être de type fini.
\vspace{1mm}

\noindent
L'idéal serait de donner à cet isomorphisme de faisceaux quasi-cohérents au sens de Ramis-Ruget une consistance plus homogène en prenant pour argument des faisceaux quasi-cohérents de cette nature ce qui revient à savoir si  ces foncteurs se prolongent aux catégories ${\rm Qcoh_R}(X)$ et ${\rm Qcoh_R}(S)$ de sorte à compléter le diagramme
$$\xymatrix{{\rm Coh}(S)\ar[d]\ar[r]^{\overline\pi}&{\rm Coh}(X)\ar[ld]_{{\rm I}\!{\rm R}^{n}\pi_{!}}\ar[d]\\
{\rm Qcoh_R}(S)\ar@/_1pc/@{.>}[r]_{T'}&{\rm Qcoh_R}(X)\ar@/_1pc/@{.>}[l]^{T}}$$
$T$ et $T'$ devant être les prolongements de ${\rm I}\!{\rm R}^{n}\pi_{!}$ et $\overline\pi$ respectivement.
\vspace{1mm}

\noindent
Ce cadre nous impose  de prendre en considération  l'aspect topologique et une lecture homologique dans le cadre des espaces vectoriels topologiques localement convexes ou, dans certains, des modules topologiques et une "atomisation" du morphisme  dans le champs des systèmes de Forster-Knorr [F.K] et {\emph résolvant de morphismes d'espaces complexes} ([Pa]. \cite{BuF1}, \cite{BuF2})  et des  techniques de descente maintenant standards,
initiées par Grothendieck, mises en formes par Deligne et jouant
un rôle important dans beaucoup d'articles et spécialement ceux de Forster-Knorr, Grauert, Kiehl-Verdier, Ramis- Ruget et  Buchweitz-Flenner. Les articles \cite{BuF1}, \cite{BuF2}  montrent la puissance de ces techniques et donnent  les ingrédients essentiels à la mise en place d'une telle stratégie dans le cadre analytique complexe. \vspace{2mm}

\noindent 

\noindent Notre objectif est de montrer les résultats suivants:
 \Th{}{}\label{thm1} Soient $n$ un entier naturel et $S$ un espace analytique complexe de dimension finie. On note\vspace{1mm}
 
 \noindent
$\bullet$  ${\Bbb E}(S; n)$ l'ensemble des morphismes propres à
fibres de dimension pure $n$, de source un espace analytique
complexe de dimension finie et paracompact  $X$, et de but $S$ supposés réduits ou à partie
régulières denses.\vspace{1mm}

\noindent
$\bullet$  ${\rm C}oh(X)$ (resp. ${\rm C}oh(S)$) la catégorie des
faisceaux cohérents sur $X$ (resp. $S$).\vspace{2mm}

\noindent Alors,
à tout élément   $\pi:X\rightarrow S$ de  ${\Bbb E}(S; n)$ est associé un unique (à isomorphisme près) foncteur $\pi^{!}_{\mathcal K}:{\rm C}oh(S)\rightarrow {\rm C}oh(X)$ vérifiant les propriétés suivantes:\vspace{1mm}

\noindent
\noindent {\bf(i)} il est covariant, exact à gauche et adjoint à droite du foncteur 
${\rm I}\!{\rm R}^{n}\pi_{*}:{\rm C}oh(X)\rightarrow {\rm
C}oh(S)$  induisant un isomorphisme bifonctoriel naturel 
$$\Theta^{\pi}_{{\mathcal F}, {\mathcal G}}:{\rm I}\!{\rm H}om(X; {\mathcal F},\pi^{!}_{\mathcal K}{\mathcal G})
\rightarrow{\rm I}\!{\rm H}om(S; {\rm I}\!{\rm R}^{n}\pi_{*}{\mathcal
F}, {\mathcal G})$$
en les arguments ${\mathcal F}$ et ${\mathcal G}$ et ``topologisé'' de façon naturelle si l'on munit ces groupes de cohomologies d'une structure topologique de type Fréchet- Schwartz.
\vspace{1mm}

\noindent {\bf(ii)} {\bf Propriétés fonctorielles relativement à  ${\Bbb E}(S;
n)$:} \vspace{1mm}

{\bf(a)}  il est de nature locale relativement aux éléments de ${\Bbb  E}(S;
n)$ dans le sens suivant:\vspace{1mm}\noindent
 si   $\pi_{i}:X_{i}\rightarrow S,\,\,i=1,2$ sont deux objets de  ${\Bbb E}(S; n)$  et $U$ un ouvert de $X_{1}$ muni
  de deux inclusions ouvertes $j_{i}:U\rightarrow X_{i}$
  tels  que le diagramme
  $$\xymatrix{&U\ar[ld]_{j_{1}}\ar[rd]^{j_{2}}&\\
  X_{1}\ar[rd]_{\pi_{1}}&&X_{2}\ar[ld]^{\pi_{2}}\\
  &S&}$$
  soit commutatif, on a
  $j_{1}^{*}{\pi_{1}}^{!}_{\mathcal K}=
  j_{2}^{*}{\pi_{2}}^{!}_{\mathcal K}$\vspace{1mm}

{\bf(b)} il est
 de formation compatible aux images directes propres entre éléments de ${\Bbb E}(S; n)$ ce qui signifie que, pour  tout diagramme
commutatif d'espaces analytiques complexes
$$\xymatrix{X_{2}\ar[rr]^{\Psi}\ar[rd]_{\pi_{2}}&&X_{1}\ar[ld]^{\pi_{1}}\\
&S&}$$
 avec $\pi_{1}$ (resp. $\pi_{2}$) appartenant à ${\Bbb  E}(S;
n)$  et $\Psi$ propre,   il existe un morphisme canonique de faisceaux cohérents
 $\Psi_{*}({{\pi_{2}}^{!}_{\mathcal K}}({\mathcal G}))\rightarrow {{\pi_1}^{!}_{\mathcal K}}({\mathcal G})$  
 injectif fibre par fibre,\vspace{2mm}
 
{\bf(c)} il est compatible avec la composition des morphismes dans le sens où, pour tout élément
  $f:X\rightarrow Y$ (resp. $g:Y\rightarrow S$) de ${\Bbb E}(Y; m)$ (resp. ${\Bbb E}(S; n)$), on a
   $$(g\circ f)^{!}_{{\mathcal K}}={f^{!}_{{\mathcal K}}}\circ {g^{!}_{{\mathcal K}}}$$
   \vspace{1mm}
   
   \noindent {\bf(iii)} {\bf Propriétés fonctorielles relativement à la base:}\vspace{1mm}
   
{\bf(a)} il est de formation compatible aux changements de base  (et donc local sur $S$), dans le sens où  pour tout morphisme  d'espaces complexes 
$\nu:S'\rightarrow S$ de diagramme de changement de base
$$\xymatrix{{X'}\ar[r]^{\nu'}\ar[d]_{\tilde{\pi}}&X\ar[d]^{\pi}\\
{S'}\ar[r]_{\nu}&S}$$
on a un morphisme naturel 
$$\nu'^{*}(\pi^{!}_{\mathcal K}({\mathcal G}))\rightarrow {\pi'^{!}_{\mathcal K}}(\nu^{*}({\mathcal G})),\,\,\,\forall\,\, {\mathcal
G}\in {\rm
 C}oh(S)$$
 qui est un isomorphisme si $\nu$ est plat.\vspace{1mm}

{\bf(b)} il est compatible avec l'image directe propre dans le sens où pour tout diagramme du type précédent mais non nécessairement de changement de base avec $\nu$ propre, on a
$$\nu'_{*}\tilde{\pi}^{!}_{\mathcal
K}(\tilde{\mathcal G})\simeq {\pi}^{!}_{{\mathcal K}}\nu_{*}(\tilde{\mathcal
G}),\,\,\,\forall\,\, \tilde{\mathcal G}\in {\rm
 C}oh(\tilde{S})$$\rm
\Th{}{}\label{thm2} Soit $\pi:X\rightarrow S$ un élément de ${\Bbb E}(S;n)$ tel que $\omega^{n}_{X/S}:=\pi^{!}_{\mathcal K}({\mathcal O}_{S})$ soit plat sur $S$. Soit $\pi^{\#}:{\rm Coh}(S)\rightarrow{\rm Coh}(X) $ le foncteur défini par $\pi^{\#}({\mathcal G}):=\pi^{*}{\mathcal G}\otimes_{{\mathcal O}_{X}}
 \omega^{n}_{X/S}$. Alors, il existe un unique (à isomorphisme près) foncteur $\pi_{\#}:{\rm Coh}(X)\rightarrow {\rm Coh}(S)$ adjoint à gauche du foncteur $\pi^{\#}$ induisant un isomorphisme de dualité
 $$\pi_{*}{\mathcal H}om({\mathcal F}, \pi^{\#}({\mathcal G}))\simeq{\mathcal H}om( \pi_{\#}{\mathcal F}, {\mathcal G})$$
 ayant toute les propriétes fonctorielles de la dualité donnée par la paire $({\rm I}\!{\rm R}^{n}\pi_{*}, \pi^{!}_{\mathcal K})$ du \theoremref{thm1}. De plus, on a  le diagramme commutatif (reliant les deux dualités) 
$$\xymatrix{\pi_{*}{\mathcal H}om({\mathcal F}, {\pi^{\#}}({\mathcal G}))\ar[d]\eq[r]&{\mathcal H}om( {\pi_{\#}}{\mathcal F}, {\mathcal G})\ar[d]\\
{\pi_{*}}{\mathcal H}om({\mathcal F}, {\pi^{!}_{\mathcal K}}{\mathcal G})\eq[r]&{\mathcal H}om({\rm I}\!{\rm R}^{n}\pi_{*}{\mathcal F}, {\mathcal G})}$$ 

 \Th{}{}\label{thm3} Avec les mêmes  notations et hypothèses, il existe un unique morphisme de foncteurs $\Xi_{\pi}:
\pi^{\#}\rightarrow \pi^{!}_{{\mathcal K}}$, génériquement bijectif en les arguments $\pi$ et ${\mathcal G}$,  de formation compatible aux restrictions ouvertes sur
$X$ et aux changements  de base plats. De plus :\vspace{2mm} 

\noindent
{\bf(i)} $\Xi_{\pi}$ est injectif si et seulement si $\pi^{\#}$ est exact à gauche,\vspace{1mm}

\noindent
{\bf(ii)} $\Xi_{\pi}$ est surjectif si et seulement si $\pi^{!}_{{\mathcal K}}$ est exact à droite,\vspace{1mm}

\noindent
{\bf(iii)} De plus, si $\pi$ est de type $S_2$, $\Xi_{\pi}$ est bijectif si et seulement si $\pi$ est un morphisme de Cohen-Macaulay.\rm\vspace{2mm}

\noindent Dans le dernier cas, sa formation
est compatible  aux changements de base arbitraires. \rm

\Th{}{}\label{thm4} Soient $\pi\in {\Bbb E}(S,n)$ et ${\rm I}\!{\rm L}\pi_{\#}$ le foncteur de 
Flenner. Alors, les assertions suivantes sont équivalentes:\vspace{2mm}

\noindent
{\bf(i)}  Le morphisme canonique $${\rm I}\!{\rm R}\pi_{*}{\mathcal F}^{\bullet}[n]\rightarrow{\rm I}\!{\rm L}\pi_{\#}({\mathcal F}^{\bullet})$$
est un isomorphisme,\vspace{2mm}

\noindent
{\bf(ii)} le morphisme $\pi$ est de Cohen Macaulay.\rm

\vfill\eject \section{\color{blue}{ Préliminaires}.}\vspace{2mm}

\noindent
\subsection{ Conventions, notations et rappels.}\rm
\vspace{2mm}

\noindent
 \subsubsection{} Dans ce qui suit, les notions citées ont été introduites  et abordées dans les détails dans \cite{K1}.\vspace{1mm}
 
 \noindent
$\bullet$ Les espaces analytiques complexes considérés ne sont
pas nécessairement réduits, sauf mention expresse du contraire,
mais on les suppose toujours  de dimension finie, dénombrables à
l'infini avec une partie régulière dense. Si $S$ est un espace
analytique complexe et $n$ un entier naturel, rappelons les
notations de \cite{K1} et \cite{K2}:\vspace{1mm}

 $\star$ ${\mathcal E}(S; n)$ (resp. ${\Bbb E}(S;n)$) est l'ensemble
des morphismes  {\emph {universellement}} $n$-{\emph {équidimensionnels}}
c'est-à-dire des morphismes d'espaces complexes de but $S$ qui
sont ouverts, surjectifs et de fibres de dimension pure $n$. ${\Bbb E}(S;n)$ désignera l'ensemble des morphismes  propres de ${\mathcal E}(S; n)$.\vspace{1mm}

$\star$ ${\mathcal E}_{eq}(S;n)$
 désignera l'ensemble des morphismes  d'espaces complexes équidimensionnels à fibres de dimension pure $n$ et ${\Bbb E}_{eq}(S;n)$ les élements propres de ${\mathcal E}_{eq}(S;n)$. \vspace{1mm}

$\star$ ${\mathcal G}_{a}(S; n)$ (resp. ${\mathcal G}_{c}(S; n)$) désigne
les éléments de ${\mathcal E}(S; n)$ qui sont {\emph {analytiquement}}
$n$-{\emph {géométriquement plats}} (resp. {\emph {continûment
$n$-géométriquement plat}}). \vspace{2mm}

 $\star$ Si $Z$ est un espace complexe de dimension bornée,
on désigne par ${\mathcal D}^{\bullet}_{Z}$ son complexe dualisant de
Ramis-Ruget (\cite{RR70}) et ${\Bbb  D}_{Z}({\mathcal B}^{\bullet}):={\rm
I}\!{\rm R}{\mathcal H}om({\mathcal B}^{\bullet}, {\mathcal D}^{\bullet}_{Z})$
 le ${\mathcal D}^{\bullet}_{Z}$-dual d'un complexe à cohomologie
cohérente donné ${\mathcal B}^{\bullet}$.
\vspace{2mm}

$\star$ A tout morphisme universellement $n$-équidimensionnel (pas
nécessairement propre!) est associé  un  faisceau particulièrement
intéressant noté $\omega^{n}_{\pi}$ et entièrement caractérisé par les propriétés  données par le {\color{red}
{théorème 1}} de \cite{K2}. Si $X$ est réduit, il coincide avec
$\tilde{\omega}_{X/S}$,  le faisceau des formes méromorphes
régulières relatives ou faisceau de Kunz-Waldi-Kersken tel qu'il est exhibé dans le {\color{red}
{théorème 3}} de \cite{K2}. Si $\pi$
est propre, ce faisceau est le faisceau dualisant relatif canonique donné par la $n$-ème homologie du complexe à cohomologie cohérente $\pi^{!}({\mathcal O}_{S})$ noté $\omega^{n}_{X/S}:={\mathcal H}^{-n}(
\pi^{!}({\mathcal O}_{S}))$ et dont les propriétés peuvent se déduire de la dualité relative comme on peut le voir dans \cite{K1},  {\emph {corollaire 3}}, p73. Si
${\mathcal F}$
 est  l'un de ces trois faisceaux, rappelons que l'on dispose d'un
morphisme d'intégration  ${\mathcal O}_{S}$-linéaire et continue
${\rm I}\!{\rm
 R}^{n}{\pi}_{!}{\mathcal F}\rightarrow {\mathcal O}_{S}$, ${\mathcal F}$
  de formation
   compatible aux restrictions ouvertes sur $X$ et aux changements
   de base plats.\vspace{2mm}
  
   \noindent
\vspace{2mm}
  
   \noindent
\subsubsection{\color{blue}{A propos de la torsion}.}\vspace{2mm}

\noindent 

\defn{}{}\vspace{1mm}

\noindent
\paragraph {\bf (a)} Si ${\mathcal A}$ est un anneau et  ${\mathcal M}$ un  ${\mathcal A}$-module, un élément $m$ de ${\mathcal M}$ est dit de {\emph {torsion}} s'il existe un élément $a$ de ${\mathcal A}$ non diviseur de zéro et tel que $a.m=0$.\vspace{2mm}
\noindent
\paragraph{\bf(b)} Soit ${\mathcal F}$  un  faisceau cohérent sur un espace analytique complexe $X$. Une section $\sigma$ de ${\mathcal F}$ est dite de {\emph {torsion}} si, pour tout $x\in X$, 
 il existe $f_{x}$ un germe de fonction holomorphe dans ${\mathcal O}_{X,x}$ tel que $f_{x}$ soit non diviseur de zéro  dans la réduction ${\mathcal O}_{{\rm Red}(X),x}$ et $f_x .\sigma_x=0$
\paragraph {\bf(c)} On appelle {\emph {sous faisceau de torsion}} d'un faisceau cohérent ${\mathcal F}$, le faisceau engendré par la réunion $\displaystyle{\bigcup_{\sigma_{x}, x\in X}\sigma_{x}}$ ou le faisceau 
associé au préfaisceau 
$$U\rightarrow \{\sigma\in \Gamma(U, {\mathcal F}): \exists\,  f\in \Gamma(U, {\mathcal O}_{X})\,\,{\rm non\,diviseur\,de\,0}\,\,{\rm tel\,que}\,f.\sigma=0\}$$\rm
\Prop{}{}\label{P1} Soit $X$ un espace analytique
complexe réduit\footnote{Rappelons qu'en géométrie analytique complexe,
un espace analytique $X$ est réduit si et seulement si la
restriction naturelle $\Gamma(X, {\mathcal O}_{X})\rightarrow
\Gamma(X-{\rm Sing}(X), {\mathcal O}_{X})$ est injective ou, ce qui
revient au même, que le lieu singulier ${\rm Sing}(X)$ est analytiquement rare (cf \cite{Fi}). Cela entraine que,  pour tout sous ensemble fermé d'intérieur
$\Sigma$ de $X$, ${\mathcal H}^{0}_{\Sigma}({\mathcal O}_{X})=
0$ .} de dimension finie.  Soit ${\mathcal
F}$  un faisceau cohérent sur $X$ et ${\mathcal T}({\mathcal F})$ son sous faisceau de torsion. On désigne par ${\mathfrak M}_{X}$  le faisceau des fonctions méromorphes au sens de \cite{Fuj}.  
Soient $\phi:{\mathcal
F}\rightarrow ({\mathcal F}^{*})^{*}:={\mathcal
H}om({\mathcal H}om({\mathcal F}, {\mathcal O}_{X}), {\mathcal O}_{X})$ et $\psi:{\mathcal F}\rightarrow {\mathcal F}\otimes_{{\mathcal O}_{X}}{\mathfrak M}_{X}$ les morphismes canoniques. Alors, on a:
$${\mathcal T}({\mathcal F})={\rm Ker}(\phi)={\rm Ker}(\psi)$$\rm
\begin{proof}
Pour vérifier l'assertion, il suffit
d'utiliser l'isomorphisme canonique
$${\mathcal H}om({\mathcal F}, {\mathcal O}_{X})\otimes_{{\mathcal O}_{X}}{\mathfrak
M}_{X}\simeq {\mathcal H}om({\mathcal F}\otimes_{{\mathcal O}_{X}}{\mathfrak M}_{X},
{\mathfrak M}_{X})$$ et de voir que, $X$ étant réduit, le
morphisme canonique ${\mathcal F}\rightarrow ({\mathcal F}^{*})^{*}$ devient
un isomorphisme une fois tensorisé par ${\mathfrak M}_{X}$ sur ${\mathcal
O}_{X}$ (cf \cite{Hou}, \cite{A60} ou \cite{Fuj} dans
lequel il trouvera la construction générale du morphisme
canonique ${\mathcal F}\rightarrow {\mathcal F}\otimes_{{\mathcal O}_{X}}{\mathfrak
M}_{X}$ au moyen des normalisations des composantes irréductibles
de $X$)$\,\blacksquare$
\end{proof}
\vspace{2mm}
\noindent
 Une autre manière de formuler cette propriété est 
\Prop{}{}\label{P2} Soit $X$ un espace analytique complexe réduit. Alors,  un  faisceau de ${\mathcal O}_{X}$-modules cohérent ${\mathcal F}$ est sans ${\mathcal O}_{X}$-torsion si et seulement si ${\mathcal H}^{0}_{\Sigma}({\mathcal
F}))=0$ pour tout fermé d'intérieur vide $\Sigma$ dans
$X$.\rm
\begin{proof} En effet, on peut, sans enfreindre la
généralité, supposer que   ${\rm Supp}({\mathcal F})=X$. Comme $X$ est
réduit, le bidual $({\mathcal F}^{*})^{*}$ est toujours sans torsion.
\vspace{1mm}

\noindent Soit ${\mathcal N}$ le sous faisceau de torsion de ${\mathcal
F}$ donné par la suite exacte
$$0\rightarrow{\mathcal N}\rightarrow {\mathcal F}\rightarrow ({\mathcal
F}^{*})^{*}$$ Comme un faisceau cohérent est toujours
génériquement localement libre sur $X$, le support du faisceau
cohérent ${\mathcal N}$ est un fermé $T$ d'intérieur vide dans
$X$. Alors cette propriété d'annulation cohomologique  donne, en
particulier, pour $\Sigma=T$, ${\mathcal H}^{0}_{T}({\mathcal N})={\mathcal
N}=0$ et, par suite, assure l'absence de torsion pour ${\mathcal F}$.ou l'injectivité de la flèche naturelle ${\mathcal F}\rightarrow ({\mathcal F}^{*})^{*}$.
\vspace{1mm}

\noindent Réciproquement, supposons cette dernière injective. Alors, il s'en suit immédiatement, au vu de l'isomorphisme
$${\mathcal H}^{0}_{\Sigma}({\mathcal H}om({\mathcal F}^{*}, {\mathcal
O}_{X}))\simeq {\mathcal H}om({\mathcal F}^{*}, {\mathcal H}^{0}_{\Sigma}({\mathcal
O}_{X})), $$ que,   pour tout ferméd'intérieur vide $\Sigma$ dans 
$X$, on a  ${\mathcal H}^{0}_{\Sigma}({\mathcal F})=0\,\blacksquare$
\end{proof}

\noindent
Dans la suite, on fera souvent usage du résultat classique suivant:
\lem{}{}\label{L1} Soient $X$, un espace analytique complexe réduit, ${\mathcal F}$ et ${\mathcal G}$ deux faisceaux  ${\mathcal O}_{X}$-cohérents. On suppose qu'il existe un morphisme $\phi:{\mathcal F}\rightarrow {\mathcal G}$ et $U$  un ouvert  dense de $X$ sur lequel $\phi$ est bijectif. Alors,\vspace{2mm}

\noindent
{\bf(i)} si ${\mathcal F}$ est  sans ${\mathcal O}_{X}$-torsion, $\phi$ est injectif.\vspace{1mm}

\noindent
{\bf(ii)} si ${\mathcal G}$ est  sans sans ${\mathcal O}_{X}$-torsion et  $\phi$  injectif, ${\mathcal F}$  est  sans ${\mathcal O}_{X}$-torsion.\rm\vspace{1mm}

\noindent
{\bf(iii)} Si ${\mathcal G}$ est sans ${\mathcal O}_{X}$-torsion, ${\mathcal H}om_{{\mathcal O}_{X}}({\mathcal F}, {\mathcal G})$  aussi.\rm
\begin{proof} Les deux premières assertions sont évidentes, la dernière découle de l'isomorphisme
$${\mathcal H}^{0}_{\Sigma}({\mathcal H}om_{{\mathcal O}_{X}}({\mathcal F}, ({\mathcal G}))\simeq ({\mathcal H}om_{{\mathcal O}_{X}}({\mathcal F}, {\mathcal H}^{0}_{\Sigma}({\mathcal G}))$$ et, si ${\mathcal G}$ est sans torsion, de l'injection
$${\mathcal H}^{1}_{\Sigma}({\mathcal H}om_{{\mathcal O}_{X}}({\mathcal F}, ({\mathcal G}))\simeq ({\mathcal H}om_{{\mathcal O}_{X}}({\mathcal F}, {\mathcal H}^{1}_{\Sigma}({\mathcal G}))$$
déduits des suites exacte courtes de bas degré  des suites spectrales de seconds termes
$${\rm E}^{i,j}_{2}:={\mathcal H}^{i}_{\Sigma}({\mathcal E}xt^{j}({\mathcal F}, {\mathcal G}),\,\,{\rm E'}^{i,j}_{2}:={\mathcal E}xt^{j}({\mathcal F}, {\mathcal H}^{i}_{\Sigma}({\mathcal G}))$$
et d'aboutissement ${\mathcal E}xt^{i+j}_{\Sigma}({\mathcal F}, {\mathcal G})$.
En effet, on a
$$0\rightarrow{\rm E}^{1,0}_{2}\rightarrow{\mathcal E}xt^{1}_{\Sigma}({\mathcal F}, {\mathcal G})\rightarrow {\rm E}^{0,1}_{2}$$ 
$$0\rightarrow{\rm E'}^{1,0}_{2}\rightarrow{\mathcal E}xt^{1}_{\Sigma}({\mathcal F}, {\mathcal G})\rightarrow {\rm E'}^{0,1}_{2}$$ 
L'injection se trouve justifiée par l'annulation  ${\rm E'}^{1,0}_{2}=$ si ${\mathcal G}$ est sans torsion$\,\blacksquare$ \end{proof}
\vspace{2mm}

\noindent
   
\subsection{Torsion et morphismes.}\vspace{2mm}

\defn{}{} Soit $\pi:X\rightarrow S$ un morphisme d'espaces complexes et ${\mathcal F}$ un faisceau cohérent sur $X$, on dira qu'il est sans $\pi$- {\emph {torsion}} ou sans {\emph {torsion fibre par
fibre}} si pour tout sous ensemble analytique $T$ de $X$ rencontrant
les fibres en des fermés d'intérieur vide, on a ${\mathcal
H}^{0}_{T}({\mathcal F})=0$ ou, ce qui revient au même, pour tout
ouvert dense fibre par fibre $U$ de $X$, le morphisme naturel de
restriction $\Gamma(X, {\mathcal F})\rightarrow \Gamma(U,{\mathcal F})$ est 
injectif.\rm\vspace{2mm}

\noindent On peut rapprocher cette définition de la
définition plus générale
 de faisceaux cohérents de type ${\rm S}_{d}$ relativementà un certain morphisme
 (\cite{EGA4}, \S 6; \cite{Ha1}, \S 3.3, ex.4). 
 \defn{}{} Soit $\pi:X\rightarrow S$
 un morphisme plat  de schémas.  Un faisceau
cohérent  ${\mathcal F}$ sur $X$ et plat sur $S$, est dit de {\emph type
$S_{r}$ relativement à $\pi$} si pour chaque $x$ de $X$ et
$\pi(x)=s$, on a
$${\rm Prof}_{{\mathcal O}_{X,x}}({\mathcal F})={\rm Min}\{k, {\mathcal H}^{k}_{x}({\mathcal F})\not=0\}\geq {\rm Prof}_{{\mathcal O}_{S,s}}({\mathcal O}_{S,s})+( min\{r, {\rm dim}{\mathcal O}_{X,x}-
 {\rm dim}{\mathcal O}_{S,s} \}$$
 \rm\vspace{1mm}
\lem{}{}\label{L2} Soit $f:X\rightarrow Y$ un morphisme fini, ouvert et surjectif d'espaces complexes réduits. Soit ${\mathcal F}$ un faisceau cohérent sur $X$. Alors, on a l'équivalence\vspace{2mm}

\noindent
\centerline{${\mathcal F}$  sans torsion si et seulement si $f_{*}{\mathcal F}$ est sans torsion.}\rm
\begin{proof}Dans tout ce qui suit et sans enfreindre la généralité, on peut supposer $Y$ irréductible. Soit, donc,  $f:X\rightarrow Y$ un morphisme ouvert, fini et  surjectif.  Soit ${\mathcal F}$ un faisceau
 cohérent sur $X$ et montrons que $f_{*}{\mathcal F}$
est sans torsion si et seulement si ${\mathcal F}$ est sans torsion.
\vspace{1mm}

\noindent $\bullet$ Supposons ${\mathcal F}$ sans torsion sur $X$. Soit
$\Sigma$ un fermé d'intérieur vide dans $S$. Alors, $f$ étant
ouvert, $f^{-1}(\Sigma)$ est aussi d'intérieur vide dans $X$.
Comme $f_{*}{\mathcal H}^{0}_{f^{-1}(\Sigma)}({\mathcal F})\simeq {\mathcal
H}^{0}_{\Sigma}(f_{*}{\mathcal F})$, on en déduit aussitôt que
$f_{*}{\mathcal F}$ est sans torsion sur $S$ puisqu'il est bien connu que pour tout faisceau analytique ${\mathcal G}$ 
et tout morphisme fini $f$, la surjection canonique $f^{*}f_{*}{\mathcal G}\rightarrow{\mathcal G}$ montre que l'annulation de $f_{*}{\mathcal G}$ entraine celle de ${\mathcal G}$
\vspace{1mm}

\noindent Réciproquement, supposons $f_{*}{\mathcal F}$ sans torsion
sur $Y$ et soit $\Sigma$ un fermé d'intérieur vide dans $X$. Par surjectivité et pour des raisons de dimension  $f(\Sigma)$ est d'intérieur vide dans $Y$. On déduit, de l'isomorphisme $f_{*}{\mathcal H}^{0}_{f^{-1}(f(\Sigma))}({\mathcal
F})\simeq {\mathcal H}^{0}_{f(\Sigma)}(f_{*}{\mathcal F})$ et de l'ouverture de
$f$ entrainant que $f^{-1}(f(\Sigma))$ est d'intérieur vide dans $X$,
 ${\mathcal H}^{0}_{f^{-1}(f(\Sigma))}({\mathcal F})=0$.  Mais l'inclusion $\Sigma\subset f^{-1}(f(\Sigma))$ et la
suite exacte d'extension des supports nous donnent finalement ${\mathcal H}^{0}_{\Sigma}({\mathcal
F})=0\,\blacksquare$
\end{proof}
\lem{}{}\label{L2'} Soient ${\mathcal F}$ et ${\mathcal G}$  deux faisceaux cohérents sur un espace complexe de sous faisceaux de torsion ${\mathcal S}$ et ${\mathcal T}$ respectivement. Alors tout morphisme de faisceaux cohérents $\phi:{\mathcal F}\rightarrow{\mathcal G}$ induit un morphisme naturel $\bar\phi:{\mathcal F}/{\mathcal S}\rightarrow{\mathcal G}/{\mathcal T}$. \rm
\begin{proof} Evidente.
\end{proof}
\vspace{2mm}

\noindent L'assertion suivante met en évidence la relation entre platitude et torsion.
\Prop{}{}\label{P3} Si
$\pi:X\rightarrow S$ est un morphisme ouvert d'espaces complexes
réduits de dimension finie à  fibres de dimension constante $n$.
Alors, on a les équivalences :\vspace{1mm}

\noindent {\bf(i)} $\pi$ est un
morphisme plat,\vspace{1mm}

\noindent {\bf(ii)} pour tout faisceau cohérent
sans torsion ${\mathcal G}$ sur $S$, $\pi^{*}{\mathcal G}$ est sans torsion
sur $X$.\rm\vspace{2mm} 

\begin{proof} Rappelons qu'un faisceau cohérent
${\mathcal F}$ sur un espace réduit  $X$ est sans  torsion si, en
tout point $x$ de $X$, la multiplication par un élément $\alpha$
 non diviseur de zéro dans l'anneau local ${\mathcal O}_{X,x}$,
 est injective sur la fibre ${\mathcal F}_{x}$ ou de façon
équivalente, si les éléments premiers de ${\mathcal O}_{X, x}$
  associés dans ${\mathcal F}_{x}$ sont exactement les éléments
  minimaux premiers de ${\mathcal O}_{X,x}$.\vspace{2mm}
  
  \noindent
{\bf (i)}$\Longrightarrow${\bf(ii)}:\vspace{1mm}

\noindent
Cette implication est évidente et peut se voir de deux façons au moins. En effet, si  ${\mathcal G}$ est un faisceau cohérent sur $S$, dire qu'il est sans torsion  revient à dire que le morphisme canonique
${\mathcal G}\rightarrow {\mathcal H}om( {\mathcal H}om({\mathcal G}, {\mathcal O}_{S}),  {\mathcal O}_{S})$ est injectif. Alors, l'assertion découle immédiatement de la platitude de $\pi$ puisque $\pi^{*}$ préserve l' injectivité et, de plus, il est bien connu que $\pi^{*} {\mathcal H}om(A, B)\simeq  {\mathcal H}om(\pi^{*}A, \pi^{*}B)$.\vspace{1mm}

\noindent On peut aussi remarquer, tout comme il a été fait dans le {\emph lemme (2.0.7)}, que ${\mathcal G}$ étant sans torsion s'injecte dans un faisceau localement libre de rang fini et, par conséquent, $\pi$étant plat, $\pi^{*}{\mathcal G}$ s'injecte aussi dans un faisceau localement libre de rang fini.\vspace{2mm}

\noindent
{\bf (ii)}$\Longrightarrow${\bf(i)}: \vspace{1mm}

\noindent
Comme le problème est de nature locale sur $X$, on est essentiellement ramené à un problème d'algèbre commutative dont le décor est fixé par
 un point  $x$ de $X$, $s:=\pi(x)$ , $A:={\mathcal O}_{S,s}$,
$B:={\mathcal O}_{X,x}$,  ${\mathfrak M}$ l'idéal maximal de $A$ (qui est la
fibre du faisceau d'idéaux sans torsion ${\mathcal I}_{S,s}$) et $k:=A/{\mathfrak M}$ le corps résiduel.\vspace{1mm}

\noindent
 Notons que,étant ouvert,  $\pi$  applique chaque idéal premier minimal de $B$ sur un idéal
minimal de $A$ ( i.e $\pi$ applique les composantes irréductibles
de $X$ sur celles de $S$). \vspace{1mm}

\noindent En vertu du critère de
platitude fibre par fibre de Grothendieck, il suffit de voir que
${\rm Tor}^{A}_{1}(k, B)=0$ ou, ce qui revient au même, de voir
que le morphisme naturel ${\mathfrak M}\otimes_{A}B\rightarrow B$ est injectif.
Mais il est déjà injectif en dehors de l'image réciproque des
fibres de l'idéal maximal de $A$, qui ne contient aucun
élément premier minimal $B$ et qui est, par conséquent,
d'intérieur vide. Comme ${\mathfrak M}\otimes_{A} B$ est un $B$-module sans
torsion, il ne peut contenir de sous modules non triviaux dont le
support est un sous ensemble d'intérieur vide; cela montre que le
noyau de ce morphisme naturel est trivial, d'où le
résultat.$\blacksquare$\end{proof}\rm\vspace{2mm}

\noindent 
On peut aussi signaler un résultat intéressant dans la pratique:
\Prop{}{} \label{P4}(\cite{Fu},{\emph lemma (5.6)}, p.43). Soit
$\pi:X\rightarrow S$ un morphisme plat et surjectif
 d'espaces complexes réduits  avec $X$ de dimension pure et $S$  irréductible.
 Soit ${\mathcal F}$ un faisceau ${\mathcal O}_{X}$-cohérent et $S$-plat. Alors, on a les équivalences suivantes:\vspace{2mm}

 {\bf (1)} Il existe un point $s$
de $S$ tel que ${\mathcal F}_{s}$ (la restriction analytique de ${\mathcal
F}$à la fibre schématique $X_{s}$)  soit sans ${\mathcal
O}_{X_{s}}$-torsion.\vspace{1mm}  

{\bf (2)} L'ensemble $U:=\{s\in S/
{\mathcal F}_{s}\,{\rm est\,sans}\, {\mathcal O}_{X_{s}}-{\rm torsion}\}$
est un ouvert dense de Zariski.\vspace{1mm}  

\indent{\bf (3)} ${\mathcal F}$
est sans ${\mathcal O}_{X}$-torsion.\rm\vspace{2mm}

\noindent
\begin{rem} 
{\bf(i)} Si les espaces sont seulement supposés sans composantes immergées et si l'on remplaçait l'expression {\emph sans ${\mathcal O}_{X}$-torsion} par {\emph sans cycles premiers immergés}, les analogues des assertions   {\bf(1)} et {\bf(2)} ne sont pas équivalentes ( pour qu'elles le soient, il faut imposer des conditions supplémentaires  par exemples  de  type ${\rm S}_{n}$ sur les fibres).\vspace{1mm}

\noindent
{\bf(ii)} Tous ces résultats sont de nature locale et reposent, en général, sur des résultats plus ou moins connus d'algèbre commutative. Dans le cas non réduit, on peut les prouver en utilisant la notion d'éléments actifs en invoquant (\cite{Fi}, {\emph{Lemma 1, Lemma 2}}, p. 144).
\end{rem}

Nous ferons souvent usage du résultat  suivant sans doute bien connu en géométrie
algébrique.
\lem{}{}\label{L3} Soient $X$, $Y$ et $S$ des espaces analytiques
complexes de Stein de dimension finie installés dans le
 diagramme commutatif d'espaces complexes sur $S$
$$\xymatrix{X\ar[rr]^{f}\ar[rd]_{\pi_{1}}&&Y\ar[ld]^{\pi_{2}}\\
&S&}$$ dans lequel $f$  est fini, ouvert et surjectif, $\pi_{1}$ et
$\pi_{2}$ sont universellement $n$-équidimensionnel.  Soit ${\mathcal
F}$ un faisceau cohérent sur $X$. Alors,\vspace{2mm}

\noindent
{\bf(i)}   ${\mathcal F}$ est sans $\pi_{2}$-torsion si et seulement si
$f_{*}{\mathcal F}$ est sans $\pi_{1}$-torsion\footnote{$^{(6)}$}{Dans
le cas absolu, le résultat connu dit que  si $g:Z\rightarrow T$
est un morphisme fini et surjectif d'espaces complexes et ${\mathcal F}$
un faisceau cohérent sur $Z$ alors ${\rm Prof}({\mathcal F})={\rm
Prof}(g_{*}{\mathcal F})$}.\vspace{1mm}

\noindent {\bf(ii)} ${\mathcal F}$ admet une
représentation locale du type
$$f^{*}({\mathcal F}_{1})\rightarrow f^{*}({\mathcal F}_{2})\rightarrow {\mathcal F}\rightarrow 0$$
dans laquelle ${\mathcal F}_{1}$ et ${\mathcal F}_{2}$ sont des faisceaux
cohérents sur $Y$.\vspace{1mm}

\noindent {\bf(iii)} Si $f$ n'est pas
nécessairement fini mais de Stein et ouvert surjectif, {\bf(ii)}
est encore vrai si l'on suppose les faisceaux quasi-cohérents au
sens de la géométrie analytique complexe.\rm\vspace{2mm}

\begin{proof} Le cas où $S$ est réduit à un point est traité dans le \lemmaref{L1}.\vspace{1mm}

\noindent Le passage au cas $S$-relatif nécessite d'adapter les arguments du cas absolu  en exploitant
l'isomorphisme $\Gamma(X, {\mathcal F})\simeq \Gamma(Y, f_{*}{\mathcal F})$
et les principales propriétés d'un morphisme ouvert  de Stein.
Il faut vérifier que les ouverts $\pi_{1}$ denses correspondent
biunivoquement aux ouverts $\pi_{2}$-dense. Pour cela,
considérons un ouvert $V$ $\pi_{2}$-dense dans $Y$ ce qui signifie
que $V\cap Y_{s}=V\cap{\pi_{2}}^{-1}(s)$ est un ouvert dense de la
fibre $Y_{s}$. L'ouverture étant  une propriété préservée
dans n'importe quel changement de base, les applications de
restrictions $f:f^{-1}(V)\rightarrow V$ et $f_{s}:X_{s}\rightarrow
Y_{s}$ sont ouvertes. Alors, les égalités  $f^{-1}(V\cap Y_{s})=
f^{-1}(V)\cap f^{-1}(Y_{s})=f^{-1}(V)\cap X_{s}$ montrent que
$f^{-1}(V)\cap X_{s}$ est un ouvert dense dans  la fibre $X_{s}$.
Alors, du diagramme commutatif
$$\xymatrix{\Gamma(Y, f_{*}{\mathcal F})\ar[r]^{\simeq}\ar[d]&\Gamma(X,
{\mathcal
F}) \ar[d]\\
\Gamma(V, f_{*}{\mathcal F})\ar[r]^{\simeq}&\Gamma(f^{-1}(V), {\mathcal F})
}$$ et de  l'absence de  $\pi_{2}$-torsion du faisceau  $f_{*}{\mathcal
F}$ entrainant l'injectivité de la première flèche verticale,
on en déduit que  ${\mathcal F}$ est sans  $\pi_{1}$-torsion.\vspace{1mm}

\noindent
Réciproquement, soit $U$ un ouvert $\pi_{1}$-dense de $X$. Alors,
$f$étant ouverte et finie, on a
$$f(U\cap X_{s})=f(U\cap f^{-1}(Y_{s}))=f(U)\cap Y_{s}$$
avec $f(U)$ dense dans $Y_{s}$. Comme  $f^{-1}(f(U))$ est $\pi_{1}$-
dense, un diagramme analogue au précédent montre que la flèche
de restriction $\Gamma(Y, f_{*}{\mathcal F})\rightarrow \Gamma(f(U),
f_{*}{\mathcal F})$ est nécessairement injective. Par surjectivité
de $f$, il en résulte que  $f_{*}{\mathcal F}$ est $\pi_{2}$-dense.
\vspace{1mm}
          
 Pour le point
{\bf(ii)}, il suffit , grâce à la surjection canonique
$f^{*}f_{*}{\mathcal F}\rightarrow {\mathcal F}$ dont le noyau est un
certain faisceau cohérent ${\mathcal K}$, de prendre ${\mathcal
F}_{2}:=f_{*}{\mathcal F}$ et ${\mathcal F}_{1}:=f_{*}{\mathcal
K}$.\vspace{1mm}\noindent Cette notion de quasi-cohérence en géométrie
analytique sera largement sollicitée dans le cas non propre. Les
arguments du
 point précédent peuvent être appliqués
  puisque la quasi-cohérence est préservée par
  image directe de Stein et le noyau d'un morphisme
   de faisceaux quasi-cohérents est encore quasi-cohérent $\blacksquare$\end{proof}\vspace{2mm}
   
   \noindent
\lem{}{}\label{L4} Soit
$\pi:X\rightarrow S$ un morphisme propre d'espaces analytiques
complexes dont les fibres sont de dimension $n$. Alors, pour tout
faisceau cohérent ${\mathcal F}$ sur $X$, \vspace{2mm}

\noindent   {\bf(i)} ${\rm
Supp}({\rm I}\!{\rm R}^{k}\pi_{*}{\mathcal F})\subset\{s\in S/ {\rm
dim}(\pi^{-1}(s)\cap{\rm Supp}({\mathcal F}) )\geq k\}$\vspace{1mm}

\noindent
{\bf(ii)} ${\rm I}\!{\rm R}^{k}{\pi_{*}}{\mathcal F}=0$ pour tout
entier $k>n$.\vspace{1mm}

\noindent
{\bf(iii)}
le foncteur ${\mathcal F}\rightarrow  {\rm I}\!{\rm R}^{n}{\pi_{*}}{\mathcal
F}$ commute aux changements de  base,\vspace{1mm}

\noindent
{\bf(iv)} Pour tout couple de faisceaux cohérents ${\mathcal F}$ et ${\mathcal G}$ sur $X$ et $S$ respectivement, le morphisme naturel
$${\mathcal G}\otimes_{{\mathcal O}_{S}}{\rm I}\!{\rm R}^{n}{\pi_{*}}{\mathcal
F}\rightarrow {\rm I}\!{\rm R}^{n}{\pi_{*}}(\pi^{*}({\mathcal G})\otimes_{{\mathcal O}_{X}}{\mathcal F})$$ 
est un isomorphisme fonctoriel en les arguments.
\rm 
\begin{proof}\vspace{1mm}

  \noindent Rappelons que le théorème de
cohérence des images directes Grauert \cite{GR60} (dont on fera un usage constant) garantit, pour tout morphisme propre $\pi:X\rightarrow S$ d'espaces complexes et
     tout faisceau cohérent ${\mathcal F}$ sur $X$,  la cohérence des faisceaux
     images directes supérieures ${\rm I}\!{\rm
R}^{k}\pi_{*}{\mathcal F}$ pour tout entier $k$. Rappelons, d'autre part, que le produit tensoriel au dessus d'une base $S$ de deux faisceaux cohérents ${\mathcal F}$ et ${\mathcal G}$ donnés sur  $X$ et $S$ respectivement, est le faisceau noté ${\mathcal F}\otimes_{{\mathcal O}_{S}}{\mathcal G}$ dont le germe en un point $\{s\}$ de $S$ est défini par la section, sur un voisinage ouvert suffisamment petit  $V$  de $s$, pour lequel on note 
$\pi^{V}:\pi^{-1}(V)\rightarrow V$ la restriction de $\pi$,   ${\mathcal F}|_{\pi^{-1}(V)}\otimes_{{\mathcal
O}_{\pi^{-1}(V)}}(\pi^{V})^{*}{\mathcal  G}$ (en fait, c'est le faisceau associé au préfaisceau  décrit par la donnée précédente). Il est clair que, pour tout entier naturel  $k$,
${\rm I}\!{\rm R}^{k}\pi_{*}({\mathcal F}\otimes_{{\mathcal O}_{S}}{\mathcal G})$
est naturellement muni d'une structure naturelle de ${\mathcal O}_{S}$-module
indépendante du choix de ${\mathcal G}$ et de type fini puisque  $\pi$ est
propre.\vspace{1mm}

\noindent On a un morphisme canonique (défini fibre par fibre) 
$${\mathcal  G}\otimes_{{\mathcal O}_{S}}({\rm I}\!{\rm R}^{k}f_{*}{\mathcal
F})_{s}\rightarrow {\rm I}\!{\rm R}^{k}f_{*}({\mathcal F}\otimes_{{\mathcal
O}_{S}}{\mathcal G})$$
En utilisant le morphisme naturel ${\mathcal  G}\rightarrow \pi_{*}\pi^{*}{\mathcal G}$, cette flèche  peut se voir comme  un morphisme de type cup-produit 
$${\rm I}\!{\rm R}^{0}\pi_{*}\pi^{*}{\mathcal  G}\otimes_{{\mathcal O}_{S}}({\rm I}\!{\rm R}^{k}\pi_{*}{\mathcal
F})_{s}\rightarrow {\rm I}\!{\rm R}^{k}\pi_{*}({\mathcal F}\otimes_{{\mathcal
O}_{S}}{\mathcal G})$$
On en déduit que les foncteurs ${\rm I}\!{\rm R}^{k}\pi_{*}({\mathcal F}$ sont stables par changement de base arbitraire mais ne commutent généralement pas à ces opérations. 
En d'autres termes, si $\eta:T\rightarrow
S$ est un morphisme d'espaces complexes muni du diagramme de changement de base
$$\xymatrix{X_{T}\ar[d]_{\pi_{T}}\ar[r]^{\Theta}&X\ar[d]_{\pi}\\
T\ar[r]_{\eta}&S}$$
La stabilité par changement de base signifie l'existence d'une flèche canonique
$$\eta^{*}{\rm I}\!{\rm R}^{k}{\pi_{*}}{\mathcal F}\rightarrow{\rm I}\!{\rm R}^{k}{{\pi_{T}}_{*}}\Theta^{*}{\mathcal F}$$
Il y'a commutation quand cette flèche est un isomorphisme. Comme on peut facilement s'en convaincre, il en est ainsi  si et seulement si ${\rm I}\!{\rm R}^{k}f_{*}$ est exacte à droite. \vspace{1mm}

\noindent Or, pour $k=n$, cette propriété est vérifiée en vertu du théorème de Reiffen \cite{R} ou du théorème des voisinages $n$-complets de \cite{B3}.  On pourra trouver un énoncé "faisceautique" dans \cite{Ta}, et sous
  une forme générale,
 dans \cite{B.V} donnant explicitement l'annulation des images directes supérieures ${\rm
I}\!{\rm
 R}^{j}{\pi}_{*}{\mathcal F}$ pour tout $j>n$ et tout faisceau
 cohérent
 ${\mathcal F}$. On peut remarquer qu'en vertu de  \cite{GR60} les supports des faisceaux ${\rm Supp}({\rm I}\!{\rm
R}^{k}\pi_{*}{\mathcal F})$ sont des sous ensemble analytique de $S$ contenu dans les sous ensembles analytiques ${\rm S}_{k}:=\{s\in S/ {\rm dim}(\pi^{-1}(s))\geq k\}$ ( on ne regarde que les fibres $\pi^{-1}(s))$ dont les supports rencontre le support de ${\mathcal F}$ \vspace{1mm}
 
 \noindent
L'inclusion  ${\rm Supp}({\rm I}\!{\rm R}^{n}\pi_{*}{\mathcal F})\subset
{\rm S}_{n}$ est une conséquence du théorème des voisinages
$n$-complets de \cite{B3} (ou sa version générale donnée dans \cite{Dem}). En effet, soit $s\notin  {\rm S}_{n}$ et ${\rm
I}\!{\rm R}^{n}\pi_{*}{\mathcal F}(s)$ la fibre du faisceau cohérent
associé au préfaisceau $U\rightarrow {\rm H}^{n}(\pi^{-1}(U),
{\mathcal F})$. Pour montrer que cette fibre en $s$ est nulle, il suffit
de voir qu'il existe un voisinage ouvert $U$ de $s$ et un voisinage
ouvert $n$-complet voisinage $V$ de $\pi^{-1}(s)\cap {\rm
Supp}({\mathcal F})$ tel que le diagramme commutatif
$$\xymatrix{ {\rm H}^{n}(\pi^{-1}(U), {\mathcal F})\ar[rr]\ar[rd]&&
{\rm I}\!{\rm R}^{n}\pi_{*}{\mathcal F}(s)\\
&{\rm H}^{n}(V, {\mathcal F})\ar[ru]&}$$ soit commutatif. Or ceci est
vrai puisque $\pi^{-1}(s)\cap {\rm Supp}({\mathcal F})$ étant compact
de dimension au plus égale à $n-1$  admet une base  de
voisinages $n$-complets. En particulier,  ${\rm H}^{n}(V, {\mathcal
F})=0$ et par conséquent ${\rm I}\!{\rm R}^{n}\pi_{*}{\mathcal
F}(s)=0$.\vspace{1mm}

\noindent
 Il en résulte que le foncteur ${\mathcal
G}\rightarrow {\rm I}\!{\rm R}^{n}\pi_{*}({\mathcal F}\otimes_{{\mathcal
O}_{S}}{\mathcal G})$ est aussi exacte à
droite. D'où la commutation aux changements de base arbitraires du faisceau  ${\rm I}\!{\rm
R}^{n}f_{*}{\mathcal F}\blacksquare$.\end{proof}
\vspace{1mm}

\noindent
\begin{rem}\vspace{1mm}

\noindent
{\bf(i)} En caractéristisque nulle, l'analogue algébrique du
théorème de changement de base n'est pas toujours vrai  sans
hypothèses de quasi-compacité et de quasi-séparabilité sur
le  morphisme donné (cf [A.H.K]). \vspace{1mm}

\noindent {\bf(ii)}
L'inclusion ${\rm Supp}({\rm I}\!{\rm R}^{n}\pi_{*}{\mathcal F})\subset
{\rm S}_{n}$ est généralement stricte. En effet, il suffit de
considérer le cas absolu donné par  $X={\Bbb P}_{1}$, $S=\{s\}$
et ${\mathcal O}_{X}$
 puisqu'alors ${\rm H}^{1}({\Bbb P}_{1},{\mathcal O}_{{\Bbb P}_{1}})= 0$ et
  ${\rm Supp}({\mathcal O}_{X})= X$. En utilisant une fibration en ${\Bbb P}_{1}$ cet exemple montre encore que l'inclusion est stricte. Dans cette situation, on a, avec  $X_{s}={\Bbb P}_{1}$, $X$ l'espace total de la famille triviale $(X_{s},s)_{s\in S}$ et ${\mathcal F}:=\sigma_{*}{\mathcal O}_{X_{s}}$, $\sigma$étant le plongement naturel de chaque fibre dans son espace total,  ${\rm H}^{1}({\Bbb P}_{1},{\mathcal O}_{{\Bbb P}_{1}})={\rm H}^{1}(X, \sigma_{*}({\mathcal O}_{X_{s}}))= 0 $  et ${\rm Supp}( \sigma_{*}({\mathcal O}_{X_{s}}))=X_{s}$. Cela montre que dans le cas absolu ou pour  des familles triviales, l'inclusion est stricte.
 \end{rem}\vspace{2mm}
 
 \noindent
\vfill\eject

\section{\color{blue}{La dualité relative pour le foncteur dérivé droit du foncteur  image directe propre et 
le foncteur $\pi^{!}_{{\mathcal K}}$.}}
\subsection{Dualité relative propre.}
\rm\vspace{4mm}

\noindent
 \paragraph { Rappelons l'important résultat de dualité relative en géométrie analytique complexe:}
 
\Th{}{} \cite{RRV71}\label{T'1}.\vspace{1mm}

\noindent
  Soit $\pi:X\rightarrow S$ un morphisme propre d'espaces analytiques complexes avec $X$ paracompact et de dimension de Zariski bornée de complexe dualisant respectifs ${\mathcal D}^{\bullet}_{X}$ et  ${\mathcal D}^{\bullet}_{S}$ .  Alors, il existe un unique foncteur $\pi^{!}:{\bf  D}^{+}_{{\rm coh}}(S)\rightarrow {\bf  D}^{+}_{{\rm coh}}(X)$ dont la formation commute aux changements de base plats et aux restrictions ouvertes sur $X$, compatible à la composition des morphismes propres, adjoint à droite du  foncteur ${\rm I}\!{\rm
 R}{\pi}_{*}:{\bf  D}_{{\rm coh}}(X)\rightarrow {\bf  D}_{{\rm coh}}(S)$. La paire $( {\rm I}\!{\rm
 R}{\pi}_{*}, \pi^{!})$ est munie d'un morphisme canonique de foncteurs
(à isomorphisme près)  ${\rm I}\!{\rm
 R}{\pi}_{*}\pi^{!}\rightarrow {\rm I}d $\footnote{ Ce morphisme est naturellement déduit du  morphisme trace   ${\rm I}\!{\rm R}\pi_{*}{\mathcal D}^{\bullet}_{X}\rightarrow{\mathcal D}^{\bullet}_{S} $ construit aussi bien dans \cite {RRV71} que
 dans \cite{RR74} ou \cite{RR70}.}
 et induit, pour tout complexe ${\mathcal A}^{\bullet}$  de  ${\mathcal O}_{X}$-modules à cohomologie cohérente (${\mathcal A}^{\bullet}\in {\bf  D}_{{\rm coh}}(X)$ ) et tout complexe ${\mathcal B}^{\bullet}$ de ${\mathcal O}_{X}$-modules à cohomologie cohérente et bornée à gauche ( ${\mathcal B}^{\bullet}\in {\bf  D}^{+}_{{\rm coh}}(S)$), un isomorphisme canonique\footnote{Signalons que
l'on dispose, pour un morphisme quelconque, de deux flèches
naturelles
$${\rm I}\!{\rm
 R}{\pi}_{*}{\rm I}\!{\rm R}{\mathcal H}om({\mathcal A}^{\bullet},{\rm I}\!{\rm R}{\rm H}om( {\rm I}\!{\rm
 L}\pi^{*}({\Bbb D}_{S}({\mathcal B}^{\bullet})), {\mathcal D}^{\bullet}_{X}))\rightarrow
{\rm I}\!{\rm R}{\mathcal H}om( {\rm I}\!{\rm R}\pi_{!}{\mathcal A}^{\bullet},
{\mathcal B}^{\bullet}),$$
l'autre étant obtenue  en échangeant ${\rm I}\!{\rm
 R}{\pi}_{!}$ et ${\rm I}\!{\rm
 R}{\pi}_{*}$.} 
$${\rm I}\!{\rm
 R}{\pi}_{*}{\rm I}\!{\rm R}{\mathcal H}om({\mathcal A}^{\bullet}, \pi^{!}({\mathcal B}^{\bullet}))\rightarrow
{\rm I}\!{\rm R}{\mathcal H}om({\rm I}\!{\rm R}\pi_{*}{\mathcal A}^{\bullet},
{\mathcal B}^{\bullet})$$ \rm
\vspace{1mm}

 \noindent
Par construction,  $\pi^{!}({\mathcal G}^{\bullet}):={\rm
I}\!{\rm R}{\mathcal H}om({\rm I}\!{\rm L}\pi^{*}{\Bbb D}_{S}({\mathcal G}^{\bullet}),
 {\mathcal D}^{\bullet}_{X}):={\Bbb D}_{X}({\rm I}\!{\rm L}\pi^{*}{\Bbb D}_{S}({\mathcal G}^{\bullet}))$  
\vspace{2mm}

\noindent
\begin{rem}
Ce théorème dit, entre autres, que pour tout morphisme propre $\pi:X\rightarrow S$ d'espaces complexes, la paire de foncteurs 
$({\rm I}\!{\rm
 R}{\pi}_{*},{\Bbb D}_{X}({\rm I}\!{\rm L}\pi^{*}{\Bbb D}_{S})(-))$ est adjointe dans les catégories précisées ci-dessus.
\end{rem}
\begin{notation} Pour un morphisme $\pi:X\rightarrow S$ de ${\Bbb
E}(S;n)$, on note ${\pi}^{!}_{{\mathcal K}}:{\rm Coh}(X)\rightarrow{\rm Coh}(S)$ le foncteur défini de la   catégorie des faisceaux cohérents sur $S$ sur celle de $X$ par  $ {\pi}^{!}_{{\mathcal K}}({\mathcal G}):={\mathcal H}^{-n}(\pi^{!}({\mathcal G})$ et on pose $\omega^{n}_{X/S}:={\mathcal H}^{-n}(\pi^{!}{\mathcal O}_{S})$.\vspace{1mm}
\end{notation}
 \centerline{\color{blue}{3.3. La preuve du \theoremref{thm1}}}\vspace{3mm}

 \noindent
On va la décomposer en deux propositions dont l'une concernera l'existence de la dualité relative de type Kleiman et l'autre les propriétés fonctorielles du foncteur $\pi^{!}_{{\mathcal K}}$ en insistant sur son comportement via-à-vis de  cette classe de morphismes (universellement) équidimensionnels.

\Prop{}{} \label{P5}\vspace{1mm}
Pour tout morphisme $\pi\in {\Bbb E}(S,n)$, il existe un unique foncteur(à isomorphisme près) $\pi^{!}_{{\mathcal K}}:{\rm Coh}(S)\rightarrow{\rm Coh}(X)$  covariant, exact à gauche  tel que la paire $({\rm I}\!{\rm R}^{n}{\pi}_{*}, \pi^{!}_{{\mathcal K}})$   muni d'un morphisme de foncteurs
$${\mathcal T}^{\pi}_{\mathcal K}:{\rm I}\!{\rm R}^{n}{\pi}_{*}\pi^{!}_{{\mathcal K}}\rightarrow {\rm Id}$$ 
 induisant, pour tout faisceaux cohérents ${\mathcal F}$ et  ${\mathcal G}$ sur $X$ et $S$ respectivement, un isomorphisme canonique
$${\rm I}\!{\rm H}om_{X}({\mathcal F}, \pi^{!}_{{\mathcal K}}({\mathcal G}))\rightarrow{\rm I}\!{\rm H}om_{S}({\rm I}\!{\rm R}^{n}{\pi}_{*}{\mathcal F}, {\mathcal G}) $$
 fonctoriel en ces arguments dont la formation commute aux changements de base plats et aux restriction ouvertes sur $X$.
 De plus, il se localise sur $S$ en l'isomorphisme de faisceaux cohérents
$$\pi_{*}{\mathcal H}om_{X}({\mathcal F}, \pi^{!}_{{\mathcal K}}({\mathcal G}))\simeq{\mathcal  H}om_{S}({\rm I}\!{\rm R}^{n}{\pi}_{*}{\mathcal F}, {\mathcal G}) $$

\begin{proof}\vspace{1mm}

\noindent

\indent
{\bf(i) Sa construction.}\vspace{1mm}

\noindent
Rappelons que pour tout morphisme propre $\pi:X\rightarrow S$, le \theoremref{T'1}
de \cite{RRV71} montre  que le foncteur $$\pi^{!}({\mathcal G}^{\bullet}):={\rm I}\!{\rm R}{\mathcal H}om({\rm I}\!{\rm L}\pi^{*}{\rm I}\!{\rm R}{\mathcal H}om({\mathcal G}, {\mathcal D}^{\bullet}_{S}), {\mathcal D}^{\bullet}_{X})$$
est  un adjoint à droite du foncteur ${\rm I}\!{\rm R}\pi_{*}$ dans la catégorie dérivée des complexes de faisceaux à cohomologie cohérente bornée inférieurement.\vspace{1mm}

\noindent
Pour un morphisme propre $\pi$ à fibres de dimension au plus égale un certain entier $n$, on note $\pi^{!}_{{\mathcal K}}$ le foncteur associé à la $n$-ème homologie donné par $\pi^{!}_{{\mathcal K}}({\mathcal G}):={\mathcal H}^{-n}(\pi^{!}{\mathcal G})$). On va voir que la paire adjointe $({\rm I}\!{\rm R}\pi_{*}, \pi^{!}) $ fournit naturellement une paire adjointe $({\mathcal H}^{0}({\rm I}\!{\rm R}\pi_{*}(-)[n]), {\mathcal H}^{0}(\pi^{!}(-)[-n]))$.
 \vspace{1mm}

\indent 
{\bf(a) Adjonction et dualité.}\vspace{1mm}

\noindent

Le morphisme trace ${\rm I}\!{\rm R}\pi_{*}{\mathcal D}^{\bullet}_{X}\rightarrow{\mathcal D}^{\bullet}_{S} $, tel qu'il est décrit dans \cite{RRV71} ou \cite{RR70}, permet de construire, pour tout complexe ${\mathcal G}^{\bullet}\in {\rm D}^{+}_{coh}(S)$, un  morphisme canonique
$${\rm T}_{\pi}: {\rm I}\!{\rm R}\pi_{*}\pi^{!}({\mathcal
  G}^{\bullet})\rightarrow {\mathcal G}^{\bullet}$$
  comme on s'en convainc aisément en écrivant 
 $${\rm I}\!{\rm R}\pi_{*}\pi^{!}{\mathcal G}^{\bullet}={\rm I}\!{\rm R}{\mathcal H}om({\rm I}\!{\rm R}{\mathcal H}om({\mathcal G}^{\bullet}, {\mathcal D}^{\bullet}_{S}), {\rm I}\!{\rm R}\pi_{*}{\mathcal D}^{\bullet}_{X})\rightarrow {\rm I}\!{\rm R}{\mathcal H}om({\rm I}\!{\rm R}{\mathcal H}om({\mathcal G}^{\bullet}, {\mathcal D}^{\bullet}_{S}), {\mathcal D}^{\bullet}_{S})={\mathcal G}^{\bullet}$$
Alors,  étant donné que, pour tout complexe ${\rm K}^{\bullet}$, on a une  flèche naturelle ${\mathcal H}^{0}({\rm K}^{\bullet})\rightarrow{\rm K}^{\bullet}$, on en déduit un  diagramme commutatif
 $$\xymatrix{{\rm I}\!{\rm R}\pi_{*}{\mathcal H}^{0}(\pi^{!}{\mathcal G}[-n])\ar[rd]\ar[r]&{\rm I}\!{\rm R}\pi_{*}\pi^{!}{\mathcal G}[-n]\ar[d]\\
 &{\mathcal G}[-n]}$$
 et, donc, le morphisme ${\rm I}\!{\rm R}\pi_{*}{\mathcal H}^{-n}(\pi^{!}{\mathcal G})[n]\rightarrow{\mathcal G}$ dont on prend la cohomologie de degré $0$, pour avoir le morphisme   (qui est déterminé par ce dernier et le détermine aussi!)\footnote{Signalons qu'aux vues des annulations des faisceaux de
cohomologie ${\mathcal H}^{j}(\pi^{!}{\mathcal G})$ pour tout $j<-n$ et ${\rm I}\!{\rm R}^{j}\pi_{*}{\mathcal F}=0$ pour tout $j>n$, on peut utiliser la suite spectrale 
$${\rm I}\!{\rm R}^{i}\pi_{*}{\mathcal H}^{j}(\pi^{!}{\mathcal G})\Longrightarrow\,{\mathcal H}^{i+j}({\rm I}\!{\rm R}\pi_{*}\pi^{!}{\mathcal G})$$
et exhiber le morphisme voulue grâce à un morphisme latéral de cette suite spectrale que l'on compose avec le morphisme trace de la dualité relative globale.}
$${\rm I}\!{\rm R}^{n}\pi_{*}{\mathcal H}^{-n}(\pi^{!}{\mathcal G})\rightarrow{\mathcal G}$$
et, donc, le morphisme de foncteurs 
 $${\mathcal T}^{\pi}_{\mathcal K}:{\rm I}\!{\rm R}^{n}{\pi}_{*}\pi^{!}_{{\mathcal K}}\rightarrow {\rm Id} $$ 
On va, à présent, montrer que le
morphisme canonique
$$\Theta^{\pi}_{{\mathcal F},{\mathcal G}}:{\rm I}\!{\rm H}om(X; {\mathcal F},{\mathcal H}^{-n}(\pi^{!}({\mathcal G})))\rightarrow  {\rm I}\!{\rm H}om(S; {\rm I}\!{\rm R}^{n}\pi_{*}{\mathcal F},{\mathcal G})$$
est un isomorphisme fonctoriel en les arguments ${\mathcal F}$ et ${\mathcal
G}$. \vspace{1mm}

\noindent Tout d'abord, remarquons que cette flèche est naturelle puisque composée
des morphismes canoniques suivants:
$${\rm I}\!{\rm H}om(X; {\mathcal F},{\mathcal H}^{-n}(\pi^{!}({\mathcal G})))\rightarrow
 {\rm I}\!{\rm H}om(S; {\rm I}\!{\rm R}^{n}\pi_{*}{\mathcal F},{\rm I}\!{\rm R}^{n}\pi_{*}{\mathcal H}^{-n}(\pi^{!}({\mathcal G}))$$
et
 $${\rm I}\!{\rm H}om(S; {\rm I}\!{\rm R}^{n}\pi_{*}{\mathcal F},
 {\rm I}\!{\rm R}^{n}\pi_{*}{\mathcal H}^{-n}(\pi^{!}({\mathcal G}))
 \rightarrow {\rm I}\!{\rm H}om(S; {\rm I}\!{\rm R}^{n}\pi_{*}{\mathcal F},{\mathcal G}),$$
ce dernier étant déduit du morphisme {\emph trace}
ou d'{\emph intégration} $${\mathcal T}^{\pi}_{\mathcal K}({\mathcal G}):{\rm I}\!{\rm
R}^{n}\pi_{*}{\mathcal H}^{-n}(\pi^{!}({\mathcal
  G}))\rightarrow {\mathcal G}$$ 
  construit précédemment.\vspace{1mm}
  
  \noindent
 {\bf{Bijectivité de $\Theta^{\pi}_{{\mathcal F},{\mathcal G}}$.} }  Elle découle du  théorème de dualité relative de \cite{RRV71} donnant, en particulier, les isomorphismes
 $${\rm I}\!{\rm R}\pi_{*}{\rm I}\!{\rm R}{\mathcal H}om({\mathcal F}, \pi^{!}({\mathcal G})[-n])\simeq {\rm I}\!{\rm R}{\mathcal H}om({\rm I}\!{\rm R}\pi_{*}{\mathcal F}, {\mathcal G}[-n])$$
 $${\rm I}\!{\rm H}om_{{\mathcal D}(X)}({\mathcal F}, \pi^{!}({\mathcal G})[-n])\simeq{\rm I}\!{\rm H}om_{{\mathcal D}(S)}({\rm I}\!{\rm R}\pi_{*}{\mathcal F}[n], {\mathcal G})$$
 qui, en vertu des annulations de cohomologie  ${\mathcal H}^{j}({\rm I}\!{\rm R}\pi_{*}{\mathcal F}[n])=0$ pour tout $j>0$, établie dans le  \lemmaref{L4} conduisent au diagramme commutatif
 $$\xymatrix{{\rm I}\!{\rm H}om_{{\mathcal D}(X)}({\mathcal F}, \pi^{!}({\mathcal G})[-n])\eq[d]\ar[r]&{\rm I}\!{\rm H}om_{X}({\mathcal H}^{0}({\mathcal F}), {\mathcal H}^{0}(\pi^{!}({\mathcal G}))[-n]))\ar[d]_{\Theta^{\pi}_{{\mathcal F},{\mathcal G}}}\\
{\rm I}\!{\rm H}om_{{\mathcal D}(S)}({\rm I}\!{\rm R}\pi_{*}{\mathcal F}[n], {\mathcal G})\eq[r]&{\rm I}\!{\rm H}om_{S}({\rm I}\!{\rm R}^{n}{\pi}_{*}{\mathcal F}, {\mathcal G})} $$
 pour lequel on montre facilement que toutes ses flèches sont des isomorphismes le réduisant au diagramme
$$\xymatrix{{\rm I}\!{\rm H}om_{{\mathcal D}(X)}({\mathcal F}, \pi^{!}({\mathcal G})[-n])\eq[d]\eq[r]&{\rm I}\!{\rm H}om_{X}({\mathcal F}, \pi^{!}_{{\mathcal K}}({\mathcal G}))\eq[d]_{\Theta^{\pi}_{{\mathcal F},{\mathcal G}}}\\
{\rm I}\!{\rm H}om_{{\mathcal D}(S)}({\rm I}\!{\rm R}\pi_{*}{\mathcal F}[n], {\mathcal G})\eq[r]&{\rm I}\!{\rm H}om_{S}({\rm I}\!{\rm R}^{n}{\pi}_{*}{\mathcal F}, {\mathcal G})} $$
On peut remarquer que cela impose les nnulations ${\mathcal H}^{j}(\pi^{!}({\mathcal G})[-n])=0$ pour tout $j<0$ et, donc, l'exactitude à gauche du foncteur $\pi^{!}_{\mathcal K}$.\vspace{1mm}

\noindent

Pour le lecteur peu familier avec ces notions, on peut dire que cela revient à utiliser la définition 
$${\rm E}xt^{p}({\mathcal F}^{\bullet}, {\mathcal G}^{\bullet})\simeq
 {\rm I}\!{\rm H}om_{{\mathcal D}(X)}({\mathcal F}^{\bullet}, {\mathcal G}^{\bullet})$$
 et les suites spectrales (pour un complexe ${\mathcal F}^{\bullet}$ borné supérieurement et un faisceau ${\mathcal G}$)
 $${\rm E}^{p;q}_{2}= {\rm E}xt^{p}({\mathcal H}^{-q}({\mathcal F}^{\bullet}), {\mathcal G})\Longrightarrow {\rm E}xt^{p+q}({\mathcal F}^{\bullet}, {\mathcal G})$$
 $${\rm E'}^{p;q}_{2}= {\rm E}xt^{p}({\mathcal H}^{i}({\mathcal F}^{\bullet}), {\mathcal H}^{q}({\mathcal G}^{\bullet})\Longrightarrow {\rm E}xt^{p+q}({\mathcal F}^{\bullet}, {\mathcal G}^{\bullet})$$
 que l'on applique à ${\mathcal F}^{\bullet}:={\rm I}\!{\rm R}\pi_{*}{\mathcal F}[n]$, dans  la première et à  ${\mathcal F}^{\bullet}:={\mathcal F}$, $i=0$ et ${\mathcal G}^{\bullet}:=\pi^{!}{\mathcal G}[-n]$ dans la seconde pour avoir les isomorphismes
 $${\rm I}\!{\rm H}om_{{\mathcal D}(S)}({\rm I}\!{\rm R}\pi_{*}{\mathcal F}[n], {\mathcal G})\simeq{\rm I}\!{\rm H}om_{S}({\rm I}\!{\rm R}^{n}{\pi}_{*}{\mathcal F}, {\mathcal G})$$
 $${\rm I}\!{\rm H}om_{{\mathcal D}(X)}({\mathcal F}, \pi^{!}({\mathcal G})[-n])\simeq{\rm I}\!{\rm H}om_{X}({\mathcal F}, \pi^{!}_{{\mathcal K}}({\mathcal G}))$$
 L'isomorphisme de  dualité relative met, alors, en bijection ces objets.
 On déduit facilement, par localisation sur $S$, l'isomorphisme
 $$\pi_{*}{\mathcal H}om_{X}({\mathcal F}, \pi^{!}_{{\mathcal K}}({\mathcal G}))\simeq{\mathcal  H}om_{S}({\rm I}\!{\rm R}^{n}{\pi}_{*}{\mathcal F}, {\mathcal G}) $$
 La fonctorialité de $$\Theta^{\pi}_{{\mathcal F},{\mathcal G}}:{\rm I}\!{\rm H}om(X; {\mathcal F},{\mathcal H}^{-n}(\pi^{!}({\mathcal G})))\rightarrow  {\rm I}\!{\rm H}om(S; {\rm I}\!{\rm R}^{n}\pi_{*}{\mathcal F},{\mathcal G})$$
en les arguments ${\mathcal F}$ et ${\mathcal
G}$ résulte de celle  de l'isomorphisme de dualité de \cite{RRV71}.\vspace{1mm}

\noindent
{\bf(ii) Propriétés intrinsèques: covariance et exactitude à gauche} \vspace{1mm}

\noindent 
Son exactitude à gauche est une conséquence de son adjonction avec le foncteur ${\rm I}\!{\rm R}\pi_{*}$ qui est exact à droite ( à noter qu'elle est aussi une conséquence des annulations  ${\mathcal H}^{-n-j}(\pi^{!}{\mathcal G})=0$ pour $j>0$ que l'on montrera plus loin dans la \propositionref{P6}) qui garantit la préservation de l'injectivité  par application du foncteur  $\pi^{!}_{{\mathcal K}}$).  \vspace{1mm}
 \noindent 
 En effet, supposons donnée  une  suite exacte de faisceaux cohérents
  sur $S$
$$0\rightarrow
 {\mathcal G}_{1}\rightarrow {\mathcal G}_{2}$$
 dont on déduit  le diagramme commutatif
 $$\xymatrix{&{\rm I}\!{\rm H}om(X;{\mathcal
 F}, \pi^{!}_{\mathcal K}({\mathcal G}_{1}))\ar[r]\eq[d]_{\Theta^{\pi}_{{\mathcal F},{\mathcal G}_{1}}}&{\rm I}\!{\rm H}om(X;{\mathcal F},
 \pi^{!}_{\mathcal K}({\mathcal G}_{2}))\eq[d]_{\Theta^{\pi}_{{\mathcal F},{\mathcal G}_{2}}}\\
0\ar[r]&{\rm I}\!{\rm H}om(S; {\rm I}\!{\rm R}^{n}\pi_{*}{\mathcal F},
  {\mathcal G}_{1})\ar[r]&{\rm I}\!{\rm H}om(S; {\rm I}\!{\rm R}^{n}\pi_{*}{\mathcal F},
 {\mathcal G}_{2})}$$
  D'où l'injectivité de la flèche horizontale supérieure et, par suite, pour ${\mathcal F}={\mathcal O}_{X}$, la suite exacte
 $$0\rightarrow\pi^{!}_{\mathcal K}({\mathcal G}_{1})\rightarrow\pi^{!}_{\mathcal K}( {\mathcal G}_{2})$$
 
 \noindent
  Comme il a été dit ce foncteur hérite des principales propriétés fonctorielles de $\pi^{!}$
 et plus précisemment de ([V],{\emph {theorem }2}, p.394; {\emph corollary 1},
 p.395 ) et de l'unicité
(à isomorphisme canonique près) de ${\pi}^{!}_{\mathcal K}\,\blacksquare$
\end{proof}\vspace{2mm}

\noindent
\rm A présent, nous allons étudier les propriétes fonctorielles comme le comportement de cette dualité vis à vis des changements de bases, des images directes ou réciproques pour la classe de morphismes universellement équidimensionnels.
\Prop{}{}\label{P6} {\it Propriétés fonctorielles.} \vspace{1mm}

\noindent La formation du foncteur $\pi^{!}_{{\mathcal K}}$ est compatible: \vspace{2mm}

\noindent
 {\bf(i) aux changements de base}:
  Si  $\nu:S'\rightarrow S$ est un morphisme d'espaces complexes induisant  le diagramme de changement de base
  $$\xymatrix{X'\ar[d]_{\pi'}\ar[r]^{\nu'}&X\ar[d]^{\pi}\\
 S'\ar[r]_{\nu}&S}$$
  on a deux morphismes canoniques de foncteurs:
 $$\nu'^{*}\pi^{!}_{{\mathcal K}} \rightarrow\pi'^{!}_{{\mathcal K}}\nu^{*}$$
 $$\nu'_{*}\pi'^{!}_{{\mathcal K}} \rightarrow\pi^{!}_{{\mathcal K}}\nu_{*}$$
 le second étant défini pour $\nu$ propre. De plus, ceux sont des isomorphismes si $\nu$ est plat.\vspace{1mm}
 
 \noindent
 Pour tout ouvert $U$ de $X$ muni d'un morphisme d'inclusion $j_U$ et de diagramme associé
$$\xymatrix{U\ar[d]_{\pi_U}\ar[r]^{j_U}&X\ar[d]^{\pi}\\
 V:=\pi(U)\ar[r]_{j_V}&S}$$
on a, $$j^{*}_U\circ (\pi^{!}_{{\mathcal K}})=({\pi_{U}})^{!}_{{\mathcal K}}\circ j^{*}_V$$
{\bf(ii) avec la
composition des morphismes universellement équidimensionnels propres} c'est-à-dire que pour tout diagramme commutatif $$\xymatrix{X\ar[rr]^{f}\ar[rd]_{\pi}&&Y\ar[ld]^{g}\\
&S&}$$ de $S$-espaces complexes avec $f\in
{\Bbb E}(Y; m)$ et $g\in{\Bbb E}(S; n)$, on a
$$\pi^{!}_{{\mathcal K}}= f^{!}_{{\mathcal
K}}\circ g^{!}_{{\mathcal K}}$$
  {\bf(iii) aux images directes
propres entre éléments de ${\Bbb E}(S; n)$}  c'est-à-dire que, pour
tout diagramme commutatif d'espaces complexes sur $S$,
 $$\xymatrix{\tilde{X}\ar[rr]^{\theta}\ar[rd]_{\tilde{\pi}}&&X\ar[ld]^{\pi}\\
&S&}$$ dans lequel  ${\pi}$ (resp. $\tilde\pi$) est propre à
fibres de
 dimension pure $n$, ${\theta}$ est propre ( donc nécessairement génériquement
 fini), on a un morphisme canonique
 $$\theta_{*}{\tilde{\pi}}^{!}_{{\mathcal K}}\rightarrow \pi^{!}_{\mathcal
 K}$$ \rm\vspace{3mm}
 
 \noindent                              
\begin{proof}

\vspace{1mm}

\noindent
{\bf (i) Propriétés de compatibilité relativement aux morphismes de changement de bases.} \vspace{1mm}

Dans la preuve qui suit, on  exploite la dualité relative de la \propositionref{P5} pour éviter d'utiliser les propriétés du foncteur $\pi^{!}$. Nous développons, néammoins cette approche, dans les preuves proposées dans la remarque.\vspace{1mm}

\noindent
 Soit $\nu:S'\rightarrow S$ un morphisme d'espaces complexes   et
$$\xymatrix{X'\ar[d]_{\pi'}\ar[r]^{\nu'}&X\ar[d]^{\pi}\\
 S'\ar[r]_{\nu}&S}$$
le diagramme de changement de base associé. Commençons par établir l'existence du morphisme canonique
$$\nu'_{*}\circ\pi'^{!}_{\mathcal K}\rightarrow\pi^{!}_{{\mathcal K}}\circ\nu_{*}$$
pour un changement de base propre. Dans ce cas, $\nu'$ est aussi propre car il est déduit de $\nu$ dans le changement de base donné par $\pi$.
La construction de ce morphisme découle du  diagramme commutatif
$$\xymatrix{{\rm I}\!{\rm H}om({X}; {{\mathcal F}}, \nu'_{*}\pi'^{!}_{\mathcal K}({\mathcal G}'))\eq[r]_{\alpha}\ar[d]_{\psi}&{\rm I}\!{\rm H}om({X'}; \nu'^{*}{{\mathcal F}}, \pi'^{!}_{\mathcal K}({\mathcal G}'))\eq[r]^{\beta}&{\rm I}\!{\rm H}om(S'; {\rm I}\!{\rm R}^{n}{\pi}'_{*}{\nu'^{*}}{{\mathcal F}}, {\mathcal G}')\ar[d]_{\gamma }\\
{\rm I}\!{\rm H}om({X}; {\mathcal F}, 
\pi^{!}_{\mathcal K}(\nu_{*}({\mathcal G}'))&\!{\rm H}om(S; {\rm I}\!{\rm R}^{n}{\pi}_{*}{\mathcal F}, \nu_{*}({\mathcal G}')\eq[l]^{\epsilon}&{\rm I}\!{\rm H}om(S'; \nu^{*}{\rm I}\!{\rm R}^{n}{\pi}_{*}{\mathcal F}, {\mathcal G}')\eq[l]^{\delta}}$$
dans  lequel $\alpha$ et $\delta$ sont les isomorphismes  donné par l'adjonction  entre $\nu^{*}$ (resp. $\nu'^{*}$) et $\nu_{*}$ (resp. $\nu'_{*}$), $\beta$ et $\epsilon$ sont des isomorphismes donnés par  la dualité de la \propositionref{P5}, $\gamma$ est un morphisme naturel déduit du changement de base 
$$\nu^{*}{\rm I}\!{\rm R}^{n}{{\pi}_{*}}{\mathcal F}\rightarrow {\rm I}\!{\rm R}^{n}{{\pi'}_{*}}\nu'^{*}{\mathcal F}$$ 
que l'on construit  en considérant les suites spectrales de Leray associées  au morphisme $\psi:=\pi\circ\nu'=\nu\circ\pi'$, 
${\rm I}\!{\rm R}^{i}{{\pi}_{*}}{\rm I}\!{\rm R}^{j}{\theta_{*}},\,\,\,{\rm I}\!{\rm R}^{i}{{\nu}_{*}}{\rm I}\!{\rm R}^{j}{\pi}'_{*}$
d'aboutissement ${\rm I}\!{\rm R}^{i+j}{{\psi}_{*}}$.  Les deux morphismes latéraux
${\rm I}\!{\rm R}^{n}{{\pi}_{*}}{\nu'_{*}}{\mathcal H}\rightarrow
{\rm I}\!{\rm R}^{n}{{\psi}_{*}}{\mathcal H}\rightarrow{{\nu}_{*}}{\rm I}\!{\rm R}^{n}{{\pi'}_{*}}{\mathcal H}$ nous donnent pour ${\mathcal H}:=\nu'^{*}{\mathcal F}$, un morphisme naturel ${\rm I}\!{\rm R}^{n}{{\pi}_{*}}{\mathcal F}\rightarrow {{\nu}_{*}}{\rm I}\!{\rm R}^{n}{{\pi'}_{*}}\nu'^{*}{\mathcal F}$ et, donc, le morphisme de changement de base.\vspace{1mm}

\noindent
En prenant, en particulier, ${\mathcal F}:={\mathcal O}_{X}$, on en déduit le morphisme de foncteurs désiré
$$\nu'_{*}\circ\pi'^{!}_{\mathcal K}\rightarrow\pi^{!}_{{\mathcal K}}\circ\nu_{*}$$
Si, de plus,  $\nu$ est plat ce qui sera le cas de $\nu'$ aussi (puisque déduit de $\nu$ dans le changement de base donné par $\pi$), le morphisme de changement de base ci-dessus est un isomorphisme et, par conséquent, $\beta$ aussi. Le morphisme de foncteurs précédemment construit est alors un  isomorphisme. \vspace{1mm}

\noindent Le second morphisme $$\nu'^{*}\circ\pi^{!}_{\mathcal K}\rightarrow\pi'^{!}_{{\mathcal K}}\circ \nu^{*}$$
ne nécessite par la propreté du changement de base  $\nu$ et se construit localement sur $X'$.\vspace{1mm}

\noindent
 Soit $x'$ un point de $X'$ et $x:=\nu'(x')$ son image dans $X$. Comme les morphisme sont équidimensionnels, on peut choisir des paramétrisations locales adaptées à notre situation. Quitte à rétrécir les données, on peut construire une installation locale factorisant $\pi'$ et $\pi$ et donnée par le diagramme commutatif
$$\xymatrix{X'\ar@/_2pc/[dd]_{\pi'}\ar[r]^{\nu'}\ar[d]_{f'}&X\ar[d]_{f}\ar@/^2pc/[dd]^{\pi}\\
S'\times U\ar[r]^{\nu"}\ar[d]_{q'}&S\times U\ar[d]_{q}\\
S'\ar[r]_{\nu}&S}$$
Avec les abus de notations que l'on connait, $U$ est un polydisque relativement compact d'un certain ${\Bbb C}^{n}$, $f$ (resp. $f'$) est un revêtement ramifié (morphisme fini, surjectif et ouvert), $q$ (resp. $q'$) la projection canonique et  $\nu":=\nu\times {\rm id}_{U}$.\vspace{1mm}

\noindent
Dans cette situation locale, on peut signaler l'isomorphisme (valable pour toute factorisation de ce type)
$$f_{*}\pi^{!}_{\mathcal K}({\mathcal F})={\mathcal H}^{-n}(\pi^{!}{\mathcal F})\simeq {\mathcal H}om(f_{*}{\mathcal O}_{X}, q^{*}{\mathcal F}\otimes \Omega^{n}_{S\times U/S})$$
que l'on peut justifier très simplement comme suit:  d'abord, par finitude de $f$, on a  
$$f_{*}{\mathcal H}^{-n}(\pi^{!}{\mathcal F})\simeq {\mathcal H}^{-n}(f_{*}\pi^{!}{\mathcal F})\simeq{\mathcal H}^{-n}({\rm I}\!{\rm R}{\mathcal H}om(f_{*}{\mathcal O}_{X}, q^{!}({\mathcal F}))\simeq {\mathcal H}om(f_{*}{\mathcal O}_{X}, q^{*}({\mathcal F})\otimes \Omega^{n}_{S\times U/S})$$
puisque 
$$\Omega^{n}_{S\times U/S}[n]\simeq{\rm I}\!{\rm R}{\mathcal H}om(q^{*}{\mathcal D}^{\bullet}_{S}, {\mathcal D}^{\bullet}_{S\times U})\simeq q^{!}({\mathcal O}_{S})$$
et que le morphisme canonique $q^{*}({\mathcal F})\otimes^{{\rm I}\!{\rm L}}q^{!}({\mathcal O}_{S})\rightarrow q^{!}({\mathcal F}) $ est un isomorphisme par platitude de  $q$. Alors, comme les morphismes finis commutent à tout changement de base, on a le diagramme
$$\xymatrix{f'_{*}\nu'^{*}(\pi^{!}_{\mathcal K}({\mathcal F}))\ar[dd]\simeq\nu"^{*}f_{*}\pi^{!}_{\mathcal K}({\mathcal F})\eq[r]&\nu"^{*}{\mathcal H}om(f_{*}{\mathcal O}_{X}, q^{*}{\mathcal F}\otimes \Omega^{n}_{S\times U/S})\ar[d]\\
&{\mathcal H}om(\nu"^{*}f_{*}{\mathcal O}_{X}, \nu"^{*}q^{*}{\mathcal F}\otimes \nu"^{*}\Omega^{n}_{S\times U/S})\eq[d]\\
f'_{*}\pi'^{!}_{{\mathcal K}}(\nu^{*}{\mathcal F})\eq[r]&{\mathcal H}om(f'_{*}(\nu'^{*}{\mathcal O}_{X}), q^{*}(\nu^{*}{\mathcal F})\otimes \Omega^{n}_{S'\times U/S'})}$$
D'où le morphisme 
$$f'_{*}\nu'^{*}(\pi^{!}_{\mathcal K}({\mathcal F}))\rightarrow f'_{*}\pi'^{!}_{{\mathcal K}}(\nu^{*}{\mathcal F})$$
qui, au vu du choix arbitraire de la factorisation locale, donne le morphisme désiré
$$\nu'^{*}\circ\pi^{!}_{\mathcal K}\rightarrow\pi'^{!}_{\mathcal K}\circ \nu^{*}$$
Si  $\nu$ est plat auquel cas  $\nu'$ aussi, la seule flèche de ce diagramme qui ne soit pas un isomorphisme le devient et, donc: 
$$\nu'^{*}\circ\pi^{!}_{\mathcal K}\simeq\pi'^{!}_{\mathcal K}\circ \nu^{*}$$
D'ailleurs, on peut retrouver cet isomorphisme directement en utilisant les isomorphismes  dualité de la  \propositionref{P5} et la commutation des images réciproques plates  aux images directes supérieures propres  et au foncteur ${\rm I}\!{\rm R}{\mathcal H}om(-,-)$. On obtient facilement l'isomorphisme 
$${\rm I}\!{\rm H}om({X'}; \nu'^{*}{{\mathcal F}}, \nu'^{*}\pi^{!}_{\mathcal K}({\mathcal G}))\rightarrow
{\rm I}\!{\rm H}om(X'; {\nu'}^{*}{{\mathcal F}}, \pi'^{!}_{\mathcal K}(\nu^{*}{\mathcal G}))$$
déduit du diagramme 
$$\xymatrix{\nu^{*}{\rm I}\!{\rm H}om(S; {\rm I}\!{\rm R}^{n}{\pi}_{*}{{\mathcal F}}, {\mathcal G})\eq[r]\eq[d]&{\rm I}\!{\rm H}om(S'; \nu^{*}{\rm I}\!{\rm R}^{n}{\pi}_{*}{{\mathcal F}}, \nu^{*}{\mathcal G})\eq[r]
&{\rm I}\!{\rm H}om(S'; {\rm I}\!{\rm R}^{n}{\pi}'_{*}{\nu'^{*}}{{\mathcal F}}, \nu^{*}{\mathcal G})\eq[d]\\
\nu'^{*}{\rm I}\!{\rm H}om({X}; {{\mathcal F}},
\pi^{!}_{\mathcal K}({\mathcal G}))\eq[r]&{\rm I}\!{\rm H}om({X'}; \nu'^{*}{{\mathcal F}}, \nu'^{*}\pi^{!}_{\mathcal K}({\mathcal G}))\ar[r]&
{\rm I}\!{\rm H}om(X'; {\nu'}^{*}{{\mathcal F}}, \pi'^{!}_{\mathcal K}(\nu^{*}{\mathcal G}))}$$
Il nous suffira de prendre encore une fois  ${\mathcal F}:={\mathcal O}_{X}$ en laissant ${\mathcal G}$ arbitraire,  pour obtenir l'isomorphisme de foncteurs souhaité.\vspace{1mm}
 
 \noindent
La compatibilité aux restrictions
ouvertes relativement aux éléments de ${\Bbb E}(S; n)$  dans le sens
  précisé dans le théorème se ramène, en con sidérant les produits fibrés sur $S$, à vérifier que pour tout
  ouvert $U$  de $X$ muni de l'immersion ouverte  $j_{U}:U\rightarrow X$, on
a un isomorphisme
 ${j_{U}}^{*}\circ\pi^{!}_{\mathcal K}() ={\pi_{U}}^{!}_{\mathcal
 K}\circ{j_{U}}^{*}$ avec $\pi_{U}=\pi|_{U}$.
 Or ceci n'est qu'un cas particulier du changement de base plat. En effet, comme $\pi$ est ouvert, $\pi(U)$ est un ouvert $V$ de $S$ et l'inclusion de cet ouvert est un changement de base plat.
\vspace{1mm}

\noindent
{\bf (ii) la  compatibilité avec la
composition des morphismes universellement équidimensionnels}.\vspace{2mm}

\noindent
Soit  $$\xymatrix{X\ar[rr]^{f}\ar[rd]_{\pi}&&Y\ar[ld]^{g}\\
&S&}$$ un diagramme commutatif de $S$-espaces complexes avec $f\in
{\Bbb E}(S; m)$ et $g\in{\Bbb E}(S; n)$.\vspace{1mm}

\noindent Comme, de toute
évidence,  $\pi\in{\Bbb E}(S; m+n)$, on va appliquer plusieurs fois la dualité de la  \propositionref{P5}. On a, alors, 
 $${\rm I}\!{\rm H}om(X; {\mathcal F},\pi^{!}_{{\mathcal K}}({\mathcal G}))\simeq
{\rm I}\!{\rm H}om(X; {\rm I}\!{\rm R}^{n+m}{\pi}_{*}{\mathcal F}, {\mathcal
G})$$ Mais l'isomorphisme de foncteurs (dû à l'annulation des images directes supérieures au delà de la dimension des fibres) ${\rm I}\!{\rm
R}^{n+m}{\pi}_{*}\simeq {\rm I}\!{\rm R}^{n}g_{*} {\rm I}\!{\rm
R}^{m}f_{*}$ et la dualité pour $f$ et $g$, nous donnent
$${\rm I}\!{\rm H}om(X; {\rm I}\!{\rm R}^{n+m}{\pi}_{*}{\mathcal F}, {\mathcal G})\simeq {\rm I}\!{\rm H}om(Y;{\rm I}\!{\rm R}^{m}f_{*}{\mathcal F},g^{!}_{{\mathcal K}}({\mathcal G}))\simeq  {\rm I}\!{\rm H}om(X; {\mathcal F}, f^{!}_{{\mathcal K}}\circ g^{!}_{{\mathcal K}}({\mathcal G}))$$
et donc  $\pi^{!}_{{\mathcal K}}= f^{!}_{{\mathcal
K}}\circ g^{!}_{{\mathcal K}}$\vspace{1mm}

\noindent
{\bf(iii) la compatibilité aux images directes
propres entre éléments de}  ${\Bbb E}(S; n)$.\vspace{1mm}

\noindent
$\bullet${\bf{Sur la source:} }On se donne un diagramme commutatif d'espaces complexes sur $S$,
 $$\xymatrix{\tilde{X}\ar[rr]^{\phi}\ar[rd]_{\tilde{\pi}}&&X\ar[ld]^{\pi}\\
&S&}$$ dans lequel  ${\pi}$ (resp. $\tilde\pi$) est propre à
fibres de
 dimension pure $n$, ${\phi}$ est propre ( donc nécessairement génériquement
 fini), on a un morphisme canonique
 $$\phi_{*}{\tilde{\pi}}^{!}_{{\mathcal K}}\rightarrow \pi^{!}_{\mathcal
 K}$$ \vspace{1mm}
 
 \noindent
Ce morphisme est construit par composition de morphismes naturels donnés dans le  diagramme
 
 $$\xymatrix{{\rm I}\!{\rm H}om(X; {\mathcal F}, {\phi}_{*} \tilde{\pi}^{!}_{{\mathcal
K}}({\mathcal G}))\ar[dd]_{\alpha_{0}}\ar[r]^{\alpha_{1}}&{\rm I}\!{\rm H}om(\tilde{X}; \phi^{*}{\mathcal F},
 \tilde{\pi}^{!}_{{\mathcal K}}({\mathcal G}))\ar[d]^{\alpha_{2}}\\
 &{\rm I}\!{\rm H}om(S; {\rm I}\!{\rm R}^{n}{\tilde\pi}_{*}{\phi}^{*}
 {\mathcal F}, {\mathcal G})\ar[d]^{\alpha_{3}}\\
{\rm I}\!{\rm H}om(X; {\mathcal F},\pi^{!}_{{\mathcal K}}({\mathcal G}))&{\rm I}\!{\rm H}om(S; {\rm I}\!{\rm R}^{n}{\pi}_{*}{\mathcal F}, {\mathcal G})\ar[l]^{\alpha_{4}}}$$ 
dans lequel: \vspace{2mm}
\noindent 

$\bullet$ $\alpha_{1}$  est l' isomorphisme d'adjonction usuel, \vspace{1mm}

\noindent

$\bullet$ $\alpha_{2}$ et $\alpha_{4}$ sont les  isomorphismes de dualité $\Theta^{\tilde\pi}_{\phi^{*}{\mathcal F},{\mathcal G}}$ et $\Theta^{\pi}_{{\mathcal F},{\mathcal G}}$ de la \propositionref{P5},\vspace{1mm}

\noindent
$\bullet$ $\alpha_{3}$ est déduit de l'identité de foncteurs ${\rm I}\!{\rm
R}{\tilde{\pi}}_{*}= {\rm I}\!{\rm R}{\pi}_{*}{\rm I}\!{\rm
 R}{\phi_{*}}$ et du morphisme de foncteurs  ${\rm Id}\rightarrow
 {\rm I}\!{\rm R}{\phi}_{*}{\phi}^{*}$,\vspace{1mm}
 
 \noindent
 
La composée de ces morphismes notée $\alpha_{0}$ donne, en particulier, pour ${\mathcal
F}={\mathcal O}_{X}$,  la flèche
$${\phi}_{*}{\tilde{\pi}}^{!}_{{\mathcal
K}}({\mathcal G}))\rightarrow {\pi}^{!}_{{\mathcal
K}}({\mathcal G})$$ 

$\bullet$ {\bf{Sur la base}}:\vspace{2mm}

\noindent
On se donne un diagramme commutatif d'espaces complexes 
$$\xymatrix{X'\ar[d]_{\pi'}\ar[r]^{\nu'}&X\ar[d]^{\pi}\\
 S'\ar[r]_{\nu}&S}$$
 dans lequel $\pi$ (resp. $\pi'$) est un  morphisme de ${\Bbb E}(S; n)$ (resp. ${\Bbb E}(S'; n')$
Alors, on a un isomorphisme de foncteurs
$$\nu'_{*}\circ{\pi'}^{!}_{\mathcal
K}\simeq {\pi}^{!}_{{\mathcal K}}\circ\nu_{*}$$
dont la construction résulte du diagramme commutatif
$$\xymatrix{{\rm I}\!{\rm H}om(X; {\mathcal F}, \nu'_{*}{\pi'}^{!}_{\mathcal
K}(\tilde{\mathcal G}) )\ar[r]\ar[d]&{\rm I}\!{\rm H}om({X'}; \nu'^{*}{\mathcal F},
{\pi'}^{!}_{\mathcal
K}(\tilde{\mathcal G}))\ar[r]&{\rm I}\!{\rm H}om({S'}; {\rm I}\!{\rm R}^{n}{\pi'}_{*}\nu'^{*}{\mathcal F}, \tilde{\mathcal G})\ar[d]\\
{\rm I}\!{\rm H}om(X; {\mathcal F}, {\pi}^{!}_{\mathcal
K}(\nu_{*}\tilde{\mathcal G}) )&
{\rm I}\!{\rm H}om(S; {\rm I}\!{\rm R}^{n}{{\pi}}_{*}{\mathcal F}, \nu_{*}\tilde{\mathcal G})\ar[l]&
{\rm I}\!{\rm H}om(S'; \nu^{*}{\rm I}\!{\rm R}^{n}{{\pi}}_{*}{\mathcal F}, \tilde{\mathcal G})\ar[l]}$$
dans lequel les flèches sont données par les isomorphismes d'adjonction, de dualité et de changement de base plat pour la cohomologie 
${\rm I}\!{\rm R}^{n}{\pi_{*}}\nu'^{*}{\mathcal F}\simeq \nu^{*}{\rm I}\!{\rm R}^{n}{\pi}_{*}{\mathcal F}$\,$\blacksquare$\end{proof}
\vspace{1mm}

\noindent Ainsi s'achève la preuve du \theoremref{thm1}.

\cor{}{}\label{C1} Avec les notations et hypothèses de la \propositionref{P6}, on a:\vspace{2mm}

\noindent
{\bf(i)} pour tout changement de base plat,  $$\theta^{*}(\Theta^{\pi}_{{\mathcal F},{\mathcal G}})=\Theta^{\tilde\pi}_{\theta^{*}{\mathcal F}, \nu^{*}{\mathcal G}}$$
{\bf(ii)} pour tout $S$-morphisme propre $\phi$,  $$\phi_{*}(\Theta^{\tilde\pi}_{\phi^{*}{\mathcal F},{\mathcal G}})=\Theta^{\pi}_{{\mathcal F},{\mathcal G}}$$
\rm

\begin{rem} Au vu du diagramme commutatif naturel
$$\xymatrix{ {\rm I}\!{\rm R}\pi_{*}{\rm I}\!{\rm R}{\mathcal H}om({\mathcal F}, \pi^{!}_{\mathcal K}({\mathcal G}))\ar[r]\ar[d]_{\chi^{\pi}_{{\mathcal F}, {\mathcal G}}}&
{\rm I}\!{\rm R}\pi_{*}{\rm I}\!{\rm R}{\mathcal H}om({\mathcal F}, \pi^{!}({\mathcal G})[-n])\eq[d]\\
{\rm I}\!{\rm R}{\mathcal H}om({\rm I}\!{\rm R}^{n}\pi_{*}{\mathcal F}, {\mathcal G})&{\rm I}\!{\rm R}{\mathcal H}om({\rm I}\!{\rm R}\pi_{*}{\mathcal F}[n], {\mathcal G})\ar[l]}$$
on peut légitimement se poser la question de savoir sous  quelles conditions peut on espérer que  $\chi^{\pi}_{{\mathcal F}, {\mathcal G}}$ soit un quasi isomorphisme. Mais cela impose de fortes  contraintes sur le morphisme car même pour $\pi$  de Cohen Macaulay sur une base lisse, on obtient un certain nombre de conditions parmi lesquelles 
$${\rm I}\!{\rm R}^{j}\pi_{*}{\mathcal F}=0$$
pour tout entier $j\not=n$ et tout faisceau cohérent ${\mathcal F}$ sur $X$. Cela peut conduire rapidement vers des impossibilités comme,  par exemple, si $X:=S\times Z$ avec $Z$ projective lisse de dimension $N$ et $S$ de Stein ou tout simplement un point, alors la dualité de Serre pour ${\mathcal F}={\mathcal O}_{Z}$ et pour ${\mathcal F}=\Omega^{N}_{Z}$ donne des conditions incompatibles; ce qui enlève tout espoir d'avoir le quasi-isomorphisme souhaité.
\end{rem}
\section{\color{blue}{Propriétés du faisceau ${\mathcal
H}^{-n}(\pi^{!}{\mathcal G})$}}
\vspace{1mm}

\noindent
L'objet de ce qui suit est de dégager, pour un faisceau cohérent ${\mathcal G}$,  un certain nombre de  propriétés du faisceau $\pi^{!}_{{\mathcal K}}({\mathcal G}):={\mathcal H}^{-n}(\pi^{!}{\mathcal G})$. Les espaces analytiques considérés seront supposés à partie régulière dense sauf mention expresse du contraire.
 \Prop{}{}\label{P6'}Soit $\pi\in {\Bbb
E}(S;n)$.  On a, alors :
\vspace{2mm}

\noindent {\bf(i)} ${\mathcal H}^{j}(\pi^{!}({\mathcal G}))= 0$ pour tout faisceau cohérent ${\mathcal G}$ sur $S$ et tout entier  $j<-n$.
\vspace{2mm}

\noindent
{\bf(ii)}\indent$\bullet$ ${\mathcal H}^{-n}(\pi^{!}{\mathcal G})$ est sans torsion si et seulement si  ${\mathcal G}$  est sans torsion pour $S$ réduit,\vspace{1mm}

\indent\indent $\bullet$ ${\mathcal H}^{-n}(\pi^{!}{\mathcal G})$ est de profondeur au moins deux si et seulement si ${\mathcal G}$ est de profondeur au moins deux pour $S$ normal, \vspace{2mm}

\noindent
{\bf(iii)} ${\mathcal H}^{-n}(\pi^{!}{\mathcal G})$ est de profondeur $S$-relative au moins deux pour tout faisceau cohérent ${\mathcal G}$ c'est-à-dire, pour tout fermé $\Sigma$ de $X$ tel que, pour tout $s\in S$,  $\Sigma\cap \pi^{-1}(s)$ soit de codimension supérieure ou égale à deux dans la fibre $\pi^{-1}(s)$, on a
$${\mathcal H}^{0}_{\Sigma}({\mathcal H}^{-n}(\pi^{!}{\mathcal G}))={\mathcal H}^{1}_{\Sigma}({\mathcal H}^{-n}(\pi^{!}{\mathcal G}))=0$$
\rm\vspace{1mm}
\noindent
 \lem{}{}\label{L5} Soient $S$ un espaces analytique complexe réduit et $V$ une variété complexe lisse. Soit $\Sigma$ un fermé de $S\times V$ tel que $\Sigma\cap (\{s\}\times V )$
soit de codimension au moins deux dans $V$ pour tout point $s$ de $S$. Alors, pour tout faisceau ${\mathcal F}$ localement libre sur $S\times V$, on a
$${\mathcal H}^{0}_{\Sigma}({\mathcal F})={\mathcal H}^{1}_{\Sigma}({\mathcal F})=0$$
\begin{proof}
 Comme le problème est de nature locale, on peut supposer ${\mathcal F}:={\mathcal O}_{S\times V}$ et on se ramène à montrer les annulations  
$${\mathcal H}^{0}_{\Sigma}({\mathcal O}_{S\times V})={\mathcal H}^{1}_{\Sigma}({\mathcal O}_{S\times V})=0$$
 L'annulation de ${\mathcal H}^{0}_{\Sigma}({\mathcal O}_{S\times V})$ est claire puisque  $S$ est réduit\footnote {Les annulations annoncées sont satisfaites si $S$ est normal puisque $\Sigma$ est, à fortiori, de codimension au moins deux dans $S\times V$}.\vspace{1mm}
 
 \noindent
On peut raisonner par récurrence descendante sur la dimension de $S$ (ou ascendante sur la codimension relative de $\Sigma$) sachant que dans le cas absolu, c'est évident par normalité de $V$. \vspace{1mm}
 
 \noindent  
 Il existe un ouvert dense $S_1$ (resp. un fermé d'intérieur vide $S'_1$) au dessus duquel $\Sigma$ est un produit $S_{1}\times \Sigma_1$ (le problème est de nature locale!)
avec $\Sigma_1$ de codimension plus grande que $2$ dans $V$. Dans ce cas, une formule de type Künneth donne facilement
$${\mathcal H}^{1}_{S_{1}\times \Sigma_{1}}({\mathcal O}_{S\times V})\simeq{\mathcal H}^{0}({\mathcal O}_{S})\widehat{\otimes}_{\Bbb C}{\mathcal H}^{1}_{ \Sigma_{1}}({\mathcal O}_{V})$$
 qui est nul en raison de la codimension de $\Sigma_1$ et de la normalité de $V$.\vspace{1mm}
 
 \noindent
Si $\Sigma'_1$ est la partie au dessus de $S'_1$, la suite exacte longue de cohomologie à support nous donne l'isomorphisme
$${\mathcal H}^{1}_{\Sigma'_1}({\mathcal O}_{S'_{1}\times Z})\simeq{\mathcal H}^{1}_{\Sigma}({\mathcal O}_{S\times Z})$$
qui permet de  conclure par hypothèses de récurrence puisque $S'_{1}$ est de dimension strictement plus petite que celle de $S$ et $\Sigma'_1$ vérifie la condition d'incidence sur ${S'_{1}}\,\blacksquare$
\end{proof}
\vspace{2mm}

\noindent

\vspace{1mm}

\noindent
 \lem{}{}\label{L6} Pour tout diagramme de factorisation locale 
 $$\xymatrix{X\ar[r]^{f}\ar[rd]_{\pi}&Y\ar[d]^{q}\\
 &S}$$
 dans lequel, comme de coutume, $f$ désigne un morphisme fini, ouvert et surjectif sur  $Y:=S\times U$ muni de la projection canonique $q$  sur $S$ et tout faisceau cohérent ${\mathcal G}$ sur $S$, le morphisme canonique:\vspace{1mm}
$$f_{*}{\mathcal H}^{j}(\pi^{!}({\mathcal G}))\simeq {\mathcal
H}^{j}(f_{*}\pi^{!}({\mathcal G}))\rightarrow{\mathcal E}xt^{n+j}_{{\mathcal O}_{Y}}
 (f_{*}{\mathcal O}_{X}, q^{*}{\mathcal G}\otimes \Omega^{n}_{Y/S})\,\,\,(\spadesuit)$$
 est un isomorphisme pour tout $j\in {\Bbb Z}$.
 \begin{proof}
Comme $\pi$ est propre, le foncteur $\pi^{!}$ est bien défini et on a
$$\pi^{!}({\mathcal G})={\rm I}\!{\rm R}{\mathcal H}om({\rm I}\!{\rm L}\pi^{*}({\mathcal D}_{S}({\mathcal G})), {\mathcal D}^{\bullet}_{X} )$$
la dualité relative pour $f$ et la formule de projection nous donnent
$$f_{*}(\pi^{!}({\mathcal G}))\simeq{\rm I}\!{\rm R}{\mathcal H}om(f_{*}{\mathcal O}_{X}\otimes^{{\rm I}\!{\rm L}}{\rm I}\!{\rm L}q^{*}({\mathcal D}_{S}({\mathcal G}), {\mathcal D}^{\bullet}_{Y} )\simeq{\rm I}\!{\rm R}{\mathcal H}om(f_{*}{\mathcal O}_{X}, {\rm I}\!{\rm R}{\mathcal H}om(q^{*}({\mathcal D}_{S}({\mathcal G}), {\mathcal D}^{\bullet}_{Y}))$$
et, par platitude de $q$, on a 
$${\rm I}\!{\rm R}{\mathcal H}om(q^{*}({\mathcal D}_{S}({\mathcal G})), {\mathcal D}^{\bullet}_{Y})=q^{!}({\mathcal G})\simeq q^{*}({\mathcal G})\otimes^{{\rm I}\!{\rm L}}q^{!}({\mathcal O}_{S})\simeq  q^{*}({\mathcal G})\otimes \Omega^{n}_{Y/S}[n] $$ 
puisque $$\Omega^{n}_{S\times U/S}[n]\simeq{\rm I}\!{\rm R}{\mathcal H}om(q^{*}{\mathcal D}^{\bullet}_{S}, {\mathcal D}^{\bullet}_{S\times U})\simeq q^{!}({\mathcal O}_{S})$$
D'où 
$$f_{*}\pi^{!}({\mathcal G})\simeq{\rm I}\!{\rm R}{\mathcal H}om(f_{*}{\mathcal O}_{X}, q^{*}({\mathcal G})\otimes \Omega^{n}_{Y/S}[n])$$
et, donc, en particulier,
$$f_{*}{\mathcal H}^{j}(\pi^{!}({\mathcal G}))\simeq {\mathcal
H}^{j}(f_{*}\pi^{!}({\mathcal G}))\simeq{\mathcal E}xt^{n+j}_{{\mathcal O}_{Y}}
 (f_{*}{\mathcal O}_{X}, q^{*}{\mathcal G}\otimes \Omega^{n}_{Y/S})$$
$$f_{*}\pi^{!}_{\mathcal K}({\mathcal G})=f_{*}{\mathcal H}^{-n}(\pi^{!}{\mathcal G})\simeq {\mathcal H}om(f_{*}{\mathcal O}_{X}, q^{*}{\mathcal G}\otimes \Omega^{n}_{S\times U/S})\,\,\blacksquare$$

\end{proof}\vspace{2mm}

\noindent
\centerline{{\bf{3.3. {Preuve de la \propositionref{P6'}}}}}\rm \vspace{3mm}

\noindent 
Comme ces propriétés sont de nature locale sur $X$, on les vérifie au voisinage de chaque point de $X$ ce qui permet de supposer $\pi$ ouvert en vertu de (\cite{Fi}, {\emph corollary}, p.138). Fixons nous un point  $x$ de $X$ et   une paramétrisation
locale (abusivement notée)  $\xymatrix{X\ar@/_/[rr]_{\pi}\ar[r]^{f}&Y:=S\times
U\ar[r]^{q}&S}$ dans laquelle $f$ est fini, ouvert, surjectif et $q$ la projection usuelle. Soit $\Sigma$ un fermé d'intérieur vide dans
$X$ et  $\tilde\Sigma:=f(\Sigma)$ son image par $f$.\vspace{1mm}

\noindent
En vertu du \lemmaref{L6}, on a, pour tout entier relatif $j$ 
$$f_{*}{\mathcal H}^{j}(\pi^{!}({\mathcal G}))\simeq {\mathcal
H}^{j}(f_{*}\pi^{!}({\mathcal G}))\simeq {\mathcal E}xt^{n+j}_{{\mathcal O}_{Y}}
 (f_{*}{\mathcal O}_{X}, q^{*}{\mathcal G}\otimes \Omega^{n}_{Y/S})$$
 et, par suite, les annulations évidentes pour $j+n<0$; ce qui prouve le premier point. 
 \vspace{1mm}

\noindent Dans la suite, on utilisera intensivement  l'isomorphisme 
$$f_{*}{\mathcal H}^{-n}(\pi^{!}{\mathcal G})\simeq {\mathcal H}om(f_{*}{\mathcal O}_{X},
q^{*}({\mathcal G})\otimes \Omega^{n}_{S\times U/S})$$
qui, grâce au \lemmaref{L2},  permet de mettre en évidence  les équivalences
$$\xymatrix{\!\!\!\!\!\!\!\!\!\!\!\!\!\!\!\!\!\!{\mathcal
H}^{-n}(\pi^{!}{\mathcal G})\,{\rm sans\, torsion}\,{\rm{(resp. de \,prof\geq2)}}{\ar@{<=>}[d]}\ar@{<=>}[r]&f_{*}{\mathcal
H}^{-n}(\pi^{!}{\mathcal G})\,{\rm sans\, torsion}\,{\rm{(resp. de \,prof\geq 2)}}\ar@{<=>}[ld]\\
{\mathcal H}om(f_{*}{\mathcal O}_{X}, q^{*}{\mathcal G}\otimes\Omega^{n}_{Y/S})\,{\rm sans\,torsion}\,{\rm{(resp. de \,prof\geq2)}}}$$
ramenant l'assertion à l'équivalence
$${\mathcal G}\,{\rm sans\, torsion (resp. de \,prof\geq2)}\,\Longleftrightarrow\,{\mathcal H}om(f_{*}{\mathcal O}_{X}, q^{*}{\mathcal G}\otimes\Omega^{n}_{Y/S})\,{\rm sans\, torsion (resp. de \,prof\geq2)}$$ 

{\bf{(i) ${\mathcal H}^{j}(\pi^{!}({\mathcal G}))=0$ pour $j<-n$.}} \vspace{1mm}

\noindent
L'isomorphisme  $$f_{*}{\mathcal H}^{j}(\pi^{!}({\mathcal G}))\simeq {\mathcal
H}^{j}(f_{*}\pi^{!}({\mathcal G}))\simeq {\mathcal E}xt^{n+j}_{{\mathcal O}_{Y}}
 (f_{*}{\mathcal O}_{X}, q^{*}{\mathcal G}\otimes \Omega^{n}_{Y/S})$$
montre  clairement  l'annulation de cette homologie pour  $j+n<0$. Comme $f$ est fini et arbitraire, la surjectivité de la flèche naturelle $f^{*}f_{*}{\mathcal F}\rightarrow {\mathcal F}$ pour tout faisceau cohérent garantit l'annulation de ${\mathcal H}^{j}(\pi^{!}({\mathcal G}))$ pour $j<-n$ sur $X$. \vspace{1mm}

\noindent
{\bf(ii) {\bf(a)} ${\mathcal H}^{-n}(\pi^{!}{\mathcal G})$\,{sans torsion}\,$\Longleftrightarrow\,{\mathcal G}$\,{sans\, torsion.}}
\vspace{1mm}

\noindent 
On a les équivalences 
$${{\mathcal G}\,{\rm sans\,torsion}\,\Longleftrightarrow\,q^{*}({\mathcal G})\,{\rm sans\,torsion}\,\Longleftrightarrow\,q^{*}{\mathcal G}\otimes \Omega^{n}_{Y/S}\,{\rm sans\,torsion}}$$
en vertu du résultat général de la  \propositionref{P3} et de la locale liberté de $\Omega^{n}_{Y/S}$.\vspace{1mm}

\noindent On en déduit, alors, en vertu du \lemmaref{L1}, que si ${\mathcal G}$ est sans torsion,  le faisceau ${\mathcal H}om(f_{*}{\mathcal O}_{X}, q^{*}{\mathcal G}\otimes\Omega^{n}_{Y/S})$ non trivial (car $f_{*}{\mathcal O}_{X}$ a pour support  $S\times U$) l'est aussi puisque  
$${\mathcal H}^{0}_{\tilde\Sigma}({\mathcal H}om(f_{*}{\mathcal O}_{X},
q^{*}({\mathcal G})\otimes \Omega^{n}_{S\times U/S})\simeq{\mathcal H}om(f_{*}{\mathcal O}_{X}, {\mathcal H}^{0}_{\tilde\Sigma}(q^{*}({\mathcal G})\otimes \Omega^{n}_{S\times U/S})) $$
et donc 
$${\mathcal H}^{0}_{\tilde\Sigma}(f_{*}{\mathcal H}^{-n}(\pi^{!}{\mathcal G}))=0$$
\vspace{1mm}

\noindent
Or l'inclusion naturelle $\Sigma\subset f^{-1}f(\Sigma)$ et l'exactitude à gauche du foncteur $f_{*}$ nous donne l'injection
$$0\rightarrow f_{*}{\mathcal H}^{0}_{\Sigma}({\mathcal H}^{-n}(\pi^{!}{\mathcal G}))\rightarrow f_{*}{\mathcal H}^{0}_{f^{-1}f(\Sigma)}({\mathcal H}^{-n}(\pi^{!}{\mathcal G}))\simeq{\mathcal H}^{0}_{\tilde\Sigma}(f_{*}{\mathcal H}^{-n}(\pi^{!}{\mathcal G})) $$
de laquelle résulte 
$${\mathcal H}^{0}_{\Sigma}({\mathcal H}^{-n}(\pi^{!}{\mathcal G}))=0$$
\vspace{1mm}

\noindent
Reciproquement,  si ${\mathcal H}om(f_{*}{\mathcal O}_{X}, q^{*}{\mathcal G}\otimes\Omega^{n}_{Y/S})$ est sans torsion sur $Y$ ou, de façon  équivalente  ${\mathcal H}om(f_{*}{\mathcal O}_{X}, q^{*}{\mathcal G})$ sans torsion.\vspace{1mm}

\noindent Comme  $f$ est un revêtement ramifié, on dispose d'un morphisme trace $\displaystyle{f_{*}{\mathcal O}_{X}\rightarrow {\mathcal O}_{S\times U}}$ qui, au vu de l'image réciproque $\displaystyle{ {\mathcal O}_{S\times U}\rightarrow  f_{*}{\mathcal O}_{X}}$ réalise ${\mathcal H}om({\mathcal O}_{S\times U}, q^{*}{\mathcal G})\simeq q^{*}{\mathcal G} $ comme un facteur direct de   ${\mathcal H}om(f_{*}{\mathcal O}_{X}, q^{*}{\mathcal G})$ qui est, par hypothèse, sans torsion. Ainsi,  $q^{*}{\mathcal G}$ est sans torsion c'est-à-dire ${\mathcal G}$ sans torsion.\vspace{2mm}

\noindent 
{\bf(ii) {\bf(b)} Si $S$ est normal, ${\mathcal H}^{-n}(\pi^{!}{\mathcal G})$\,{de profondeur au moins deux}\,$\Longleftrightarrow\,{\mathcal G}$\,{\bf de profondeur au moins deux.}} \vspace{1mm}

\noindent
Supposons ${\mathcal G}$ de profondeur au moins deux (en particulier  sans torsion). 
On peut l'installer localement dans une suite exacte courte
$$0\rightarrow {\mathcal G}\rightarrow {\mathcal L}\rightarrow {\mathcal K}\rightarrow 0$$
avec ${\mathcal L}$ localement libre et ${\mathcal K}$ sans torsion en raison de la profondeur de ${\mathcal G}$.
Comme $q$ est plat et $\Omega^{n}_{Y/S}$ localement libre, en résulte la suite exacte courte 
$$0\rightarrow q^{*}({\mathcal G})\otimes\Omega^{n}_{Y/S} \rightarrow q^{*}({\mathcal L})\otimes\Omega^{n}_{Y/S}\rightarrow q^{*}({\mathcal K})\otimes\Omega^{n}_{Y/S}\rightarrow 0$$
qui, pour tout fermé $\tilde\Sigma$ de codimension au moins deux dans $Y$, donne la suite exacte longue de cohomologie
$$0\rightarrow {\mathcal H}^{0}_{\tilde\Sigma}(q^{*}({\mathcal K})\otimes\Omega^{n}_{Y/S})\rightarrow {\mathcal H}^{1}_{\tilde\Sigma}(q^{*}({\mathcal G})\otimes\Omega^{n}_{Y/S})\rightarrow {\mathcal H}^{1}_{\tilde\Sigma}(q^{*}({\mathcal L})\otimes\Omega^{n}_{Y/S})\rightarrow\cdots$$
Mais le premier faisceau est nul puisque ${\mathcal K}$ est sans torsion, le troisième l'est aussi puisque $S$ est normal. Ainsi,
$${\mathcal H}^{1}_{\tilde\Sigma}(q^{*}({\mathcal G})\otimes\Omega^{n}_{Y/S})=0$$
et, par suite, celle de 
$${\mathcal H}^{1}_{\Sigma}({\mathcal H}^{-n}(\pi^{!}{\mathcal G}))=0$$
pour tout fermé de codimension au moins deux dans $X$. \vspace{1mm}

\noindent
Réciproquement, si ${\mathcal H}om(f_{*}{\mathcal O}_{X}, q^{*}{\mathcal G}\otimes\Omega^{n}_{Y/S})$ est de profondeur au moins deux sur $Y$ ou, de façon  équivalente  ${\mathcal H}om(f_{*}{\mathcal O}_{X}, q^{*}{\mathcal G})$ de profondeur au moins deux, le raisonnement ci-dessus utilisant le morphisme trace montre que $q^{*}{\mathcal G}$ est nécessairement de profondeur au moins deux sur $Y$ puisque facteur direct d'un faisceau ayant cette propriété. En effet, si, pour une suite exacte scindée
de faisceaux cohérents $$0\rightarrow {\mathcal A'}\rightarrow {\mathcal A}\rightarrow {\mathcal A"}\rightarrow 0$$
${\mathcal A}$ est de profondeur au moins deux et sans torsion, les facteurs directs sont  automatiquement sans torsion. De plus, les conditions ${\mathcal A"}$ sans torsion et  ${\mathcal A}$  de profondeur au moins deux imposent à  ${\mathcal A'}$ d'être de profondeur au moins deux.\vspace{2mm}

\noindent
{\bf(iii) ${\mathcal H}^{-n}(\pi^{!}{\mathcal G})$ est de profondeur $S$-relative au moins deux pour tout faisceau cohérent 
${\mathcal G}$.}\vspace{1mm}

\noindent
D'après ce que l'on a dit dans les preuves qui précèdent, il nous suffit de voir que la propriété est satisfaite pour le faisceau $q^{*}{\mathcal G}\otimes\Omega^{n}_{Y/S}$ et, ce, pour toute factorisation locale. \vspace{1mm}

\indent {\bf(a)} ${\mathcal H}^{0}_{\Sigma}(q^{*}{\mathcal G}\otimes\Omega^{n}_{Y/S})=0.$
\vspace{1mm}

\indent $\bullet$ Si ${\mathcal G}$ est  localement libre, cela découle du \lemmaref{L5}. \vspace{1mm}

\indent
$\bullet$ Si ${\mathcal G}$ est sans torsion, ce n'est qu'un cas particulier de {\bf(ii),a)}. D'ailleurs, étant sans torsion, il est localement sous faisceau d'un faisceau localement libre de rang fini ${\mathcal L}$. On a, alors, $q$ étant plat et $\Omega^{n}_{Y/S}$ $S$-plat, l'injection
$$\xymatrix{{\mathcal H}^{0}_{\Sigma}(q^{*}{\mathcal G}\otimes\Omega^{n}_{Y/S})\ar@{^{(}->}[r]&{\mathcal H}^{0}_{\Sigma}(q^{*}{\mathcal L}\otimes\Omega^{n}_{Y/S})}$$
et l'annulation désirée.\vspace{1mm}

\indent 
$\bullet$ Si ${\mathcal G}$ est cohérent quelconque, on considère la suite exacte courte
$$ 0\rightarrow {\mathcal T}\rightarrow {\mathcal G}\rightarrow {\mathcal K}:={\mathcal G}/{\mathcal T}\rightarrow 0$$
où ${\mathcal T}$ est le sous faisceau de torsion de ${\mathcal G}$ et qui induit 
 la suite exacte courte
$$0\rightarrow {\mathcal H}^{0}_{\Sigma}(q^{*}{\mathcal T}\otimes\Omega^{n}_{Y/S})\rightarrow {\mathcal H}^{0}_{\Sigma}(q^{*}{\mathcal G}\otimes\Omega^{n}_{Y/S})\rightarrow{\mathcal H}^{0}_{\Sigma}(q^{*}{\mathcal K}\otimes\Omega^{n}_{Y/S})\rightarrow 0$$
En notant $S_0:={\rm Supp}({\mathcal T})$, ${{\mathcal T}_{0}}:={\mathcal T}$, 
 $\pi_{0}:X_{0}\rightarrow S_{0}$ (resp. $q_{0}:Y_{0}:=S_{0}\times U\rightarrow S_{0}$) le morphisme déduit de $\pi$ (resp. $q$) par changement de base et $\Sigma_{0}:=\Sigma\cap X_{0}$, on a
$${\mathcal H}^{0}_{\Sigma_{0}}({q_{0}}^{*}{\mathcal T}_{0}\otimes\Omega^{n}_{Y_{0}/S_{0}})\simeq {\mathcal H}^{0}_{\Sigma}(q^{*}{\mathcal T}\otimes\Omega^{n}_{Y/S})\simeq{\mathcal H}^{0}_{\Sigma}(q^{*}{\mathcal G}\otimes\Omega^{n}_{Y/S})$$
Si ${\mathcal T}_{0}$  est sans torsion, l'assertion  est prouvée sinon on procède au dévissage de ${\mathcal T}_{0}$ jusqu'à obtenir un indice $k$ pour lequel  ${\mathcal T}_{k}$ est sans torsion sur $S_k$.\vspace{1mm}

\noindent 
\indent {\bf(b)} ${\mathcal H}^{1}_{\Sigma}(q^{*}{\mathcal G}\otimes\Omega^{n}_{Y/S})=0$
\vspace{1mm}

\indent
$\bullet$ Si ${\mathcal G}$ est  localement libre, on invoque le \lemmaref{L5}, \vspace{1mm}

\indent si ${\mathcal G}$ est sans torsion, on reprend la suite exacte courte locale 
$$ 0\rightarrow {\mathcal G}\rightarrow {\mathcal L}\rightarrow {\mathcal G'}\rightarrow 0$$
dont on déduit la suite exacte
$$\xymatrix{0\ar[r]&{\mathcal H}^{0}_{\Sigma}(q^{*}{\mathcal G'}\otimes\Omega^{n}_{Y/S})\ar[r]& {\mathcal H}^{1}_{\Sigma}(q^{*}{\mathcal G}\otimes\Omega^{n}_{Y/S})\ar[r]&{\mathcal H}^{1}_{\Sigma}(q^{*}{\mathcal L}\otimes\Omega^{n}_{Y/S})\ar[r]&\cdots}$$
et, par suite, l'annulation ${\mathcal H}^{1}_{\Sigma}(q^{*}{\mathcal G}\otimes\Omega^{n}_{Y/S})=0$ en vertu du {\bf(a)} et du lemme cité ci-dessus.\vspace{1mm}

\indent Si ${\mathcal G}$ est cohérent quelconque ayant un sous faisceau de torsion ${\mathcal T}$, on a, avec les mêmes notations que ci-dessus
$${\mathcal H}^{1}_{\Sigma}(q^{*}{\mathcal T}\otimes\Omega^{n}_{Y/S})\simeq{\mathcal H}^{1}_{\Sigma}(q^{*}{\mathcal G}\otimes\Omega^{n}_{Y/S})$$
en raison des annulations
$${\mathcal H}^{0}_{\Sigma}(q^{*}{\mathcal K}\otimes\Omega^{n}_{Y/S})={\mathcal H}^{1}_{\Sigma}(q^{*}{\mathcal K}\otimes\Omega^{n}_{Y/S})=0$$
c'est-à-dire
$${\mathcal H}^{1}_{\Sigma_{0}}({q_{0}}^{*}{\mathcal T}_{0}\otimes\Omega^{n}_{Y_{0}/S_{0}})\simeq {\mathcal H}^{1}_{\Sigma}(q^{*}{\mathcal G}\otimes\Omega^{n}_{Y/S})$$
et on termine par récurrence sur la dimension du support du sous de torsion ou en réitérant le dévissage le long de la torsion de ${\mathcal T}\,\blacksquare$
\vspace{2mm}

\noindent
\begin{rem}
    On peut remarquer que le  \lemmaref{L2} montre que l'on peut toujours supposer $X$ réduit si $S$ l'est.  En effet, si  ${\mathcal N}$ est le sous faisceau des nilpotents de ${\mathcal
O}_{X}$ dont on sait qu'il a un support  d'intérieur vide dans $X$ ou rare puisque  la partie régulière de $X$ (qui est
réduite) est dense dans $X$, on considère   la suite exacte courte
$$ 0\rightarrow {\mathcal N}\rightarrow {\mathcal O}_{X}\rightarrow {{\mathcal O}_{X}}/{\mathcal N}\rightarrow 0$$
dont on déduit, par cohérence des faisceaux considérés et finitude de
$f$,  la suite exacte
$$0\rightarrow f_{*}{\mathcal N}\rightarrow f_{*}{\mathcal O}_{X}\rightarrow f_{*}({{\mathcal O}_{X}}/{\mathcal N})\rightarrow 0$$
et de laquelle résulte la suite exacte
$$0\rightarrow {\mathcal H}om(f_{*}({\mathcal O}_{X}/{\mathcal N}),
 q^{*}{\mathcal G})\rightarrow {\mathcal H}om(f_{*}{\mathcal O}_{X}, q^{*}{\mathcal F})\rightarrow {\mathcal H}om(f_{*}({\mathcal N}), q^{*}{\mathcal G})$$
 Comme  $q^{*}{\mathcal G}$ est sans torsion pour  ${\mathcal G}$ sans torsion et que, par hypothèse de densité sur  la partie régulière de $S$, le support de $f_{*}({\mathcal N})$ (qui est  un sous
ensemble analytique strict de $Y$)  est d'intérieur vide dans $Y$, on a ${\mathcal H}om(f_{*}({\mathcal N}), q^{*}{\mathcal G})=0$. On en déduit immédiatement  par tensorisation par le faisceau localement libre $\Omega^{n}_{Y/S}$, l'isomorphisme
 $${\mathcal H}om(f_{*}({\mathcal O}_{X}/{\mathcal N}),
 q^{*}{\mathcal G}\otimes\Omega^{n}_{Y/S} )\simeq{\mathcal H}om(f_{*}{\mathcal O}_{X}, q^{*}{\mathcal G}\otimes\Omega^{n}_{Y/S})$$
\end{rem}

\cor{}{}\label{c2}
Le foncteur ${\mathcal G}\rightarrow \pi^{!}_{{\mathcal K}}({\mathcal G}):={\mathcal H}^{-n}(\pi^{!}{{\mathcal G}})$ est exact à gauche.\rm
\begin{proof}
C'est une conséquence immédiate de ce qui précède puisque l'annulation ${\mathcal H}^{-n-1}(\pi^{!}{{\mathcal G}})=0$ pour tout faisceau cohérent ${\mathcal G}$ équivaut au fait que le foncteur $\pi^{!}_{{\mathcal K}}$ préserve l'injectivité$\,\blacksquare$
\end{proof}
\cor{}{}\label{c3} Soit $\pi:X\rightarrow S$ un élément de ${\Bbb
E}(S;n)$.  Alors, le faisceau ${\mathcal O}_{X}$-cohérent
$\omega^{n}_{X/S}:={\mathcal H}^{-n}(\pi^{!}{\mathcal O}_{S})$ satisfait les propriétés suivantes:\vspace{2mm}

\noindent
{\bf(i)} il est sans ${\mathcal
O}_{X}$-torsion pour $S$ réduit et  de
profondeur au moins deux sur $X$ si $S$ est normal,\vspace{1mm}

\noindent
{\bf(ii)} il est de profondeur au moins deux fibre
par fibre, en particulier, pour tout fermé $\Sigma$ de $X$
rencontrant les fibres en des fermés de codimension au moins $k$
($k=1,2$),  ${\mathcal H}^{j}_{\Sigma}({\omega}^{n}_{X/S})= 0$
pour $j<k$.\rm
\begin{proof} Immédiate d'après la \propositionref{P6} pour ${\mathcal G}:={\mathcal O}_{S}\,\blacksquare$
\end{proof}
\cor{}{}\label{c4} Soit $\pi:X\rightarrow S$ propre, plat à fibres de dimension $n$ et de type ${\rm S}_{1}$. Soit $Y$ le fermé de $X$ constitué des fibres non de Cohen Macaulay. Alors
$${\mathcal H}^{i}_{Y}(\omega^{n}_{X/S})=0,\, {\rm pour}\,i=0,1$$
\vspace{1mm} 

\noindent
\begin{proof} Comme le problème est de nature locale, on peut supposer donné une factorisation du type $\xymatrix{X\ar@/_/[rr]_{\pi}\ar[r]^{\sigma}&S\times Z\ar[r]^{q}&S}$
 avec  $Z$ variété de Stein de dimension $n+p$, $\sigma$ un plongement et $q$ la projection canonique.\vspace{1mm}

 \noindent
Le résultat reposent, en grande partie, sur les lemmes {\bf(5.2.2)} et {\bf(5.2.5)} de \cite{A.L} se laissant facilement  transposés à  notre cadre. Sous ces conditions, il y'est montré que 
 l'ensemble $V:=\{x\in X: X_{\pi(x)}\,{\rm est\, de\, Cohen\, Macaulay\,en\,x}\}$ est un ouvert de $X$ dont le complémentaire $Y$ intersecte chaque fibre en codimension au moins deux dans cette fibre;  $(Y\cap X_{s})$ est de codimension au moins deux dans  $X_{s}$. Dans ce cas, c'est immédiat d'après le \corollaryref{c3} $\,\blacksquare$
 \end{proof}\vspace{2mm}

\noindent
dont on déduit le
\cor{}{}\label{c5}{}
Avec les notations et hypothèses précédentes, on a, pour tout faisceau cohérent ${\mathcal G}$ de profondeur au moins deux sur $S$:
$${\mathcal H}^{i}_{Y}(\pi_{{\mathcal K}}^{!}({\mathcal G}))=0,\,\, {\rm pour}\, i=0,1$$
\rm
 \begin{proof}
 \vspace{1mm}
 
 \noindent
 $\bullet$ Supposons, dans un premier temps, ${\mathcal G}$ réfléxif  et considérons  un début de résolution localement libre du ${\mathcal O}_{S}$-dual ${\mathcal G}^{*}$
$${\mathcal L}^{1}\rightarrow{\mathcal L}^{0}\rightarrow {\mathcal G}^{*}\rightarrow 0$$ 
En appliquant le foncteur ${\mathcal H}om(-, {\mathcal O}_{S})$ et utilisant la réfléxivité, il est facile de mettre en évidence  une suite exacte courte du type
 $$ 0\rightarrow{\mathcal G}\rightarrow{\mathcal L}\rightarrow{\mathcal H}\rightarrow 0$$
 avec ${\mathcal L}$ localement libre et ${\mathcal H}$ sans torsion. Alors, en vertu des propriétés du foncteur $\pi_{{\mathcal K}}^{!}$ données par la \propositionref{P5} et le \corollaryref{C1} (en particulier celle disant que ce foncteur est stable sur les faisceaux sans torsion ),  on déduit  une suite exacte courte
 $$ 0\rightarrow\pi_{{\mathcal K}}^{!}({\mathcal G})\rightarrow\pi_{{\mathcal K}}^{!}({\mathcal L})\rightarrow{\mathcal K}\rightarrow 0$$
 dans laquelle ${\mathcal K}$ est le faisceau image du morphisme $\pi_{{\mathcal K}}^{!}({\mathcal L})\rightarrow \pi_{{\mathcal K}}^{!}({\mathcal H})$ qui est  un faisceau cohérent sans torsion puisque $\pi_{{\mathcal K}}^{!}({\mathcal H})$ l'est. \vspace{1mm}
 
 \noindent
 Il suffit, alors, d'écrire la suite  exacte longue de cohomologie à support dans $Y$,  d'utiliser l'annulation ${\mathcal H}^{0}_{Y}({\mathcal K})=0$ et le point {\bf(ii)} de la proposition dont découle l'annulation ${\mathcal H}^{1}_{Y}(\pi_{{\mathcal K}}^{!}({\mathcal L}))=0$ eu égard à l'isomorphisme naturel
 $${\mathcal H}^{i}_{Y}(\pi_{{\mathcal K}}^{!}({\mathcal L}))\simeq{\mathcal L}\otimes {\mathcal H}^{i}_{Y}(\pi_{{\mathcal K}}^{!}({\mathcal O}_{S}))={\mathcal L}\otimes{\mathcal H}^{i}_{Y}(\omega^{n}_{X/S})$$ 
 $\bullet$ Si ${\mathcal G}$ est un faisceau de profondeur au moins deux, on applique ce qui précède à la suite exacte courte naturelle
$$0\rightarrow{\mathcal G}\rightarrow{\mathcal G}^{**}\rightarrow {\mathcal H}\rightarrow 0$$
où ${\mathcal H}$ est le faisceau quotient qui est  sans torsion$\,\blacksquare$ 
\end{proof}\vspace{2mm}

\noindent
\cor{}{}\label{c6}{} Supposons le foncteur $\pi^{\#}$ exact à gauche. Alors,  pour tout faisceau cohérent ${\mathcal G}$ de profondeur au moins deux dans $S$,  le morphisme  $\Xi_{\pi}({\mathcal G})$ est un isomorphisme.\rm
\begin{proof}

Considérons, d'abord, le cas d'un faisceau réfléxif  ${\mathcal G}$ et un début de résolution localement libre du ${\mathcal O}_{S}$-dual ${\mathcal G}^{*}$
$${\mathcal L}^{1}\rightarrow{\mathcal L}^{0}\rightarrow {\mathcal G}^{*}\rightarrow 0$$ dont on déduit le début de résolution
$$ 0\rightarrow{\mathcal G}\rightarrow({\mathcal L}^{0})^{*}\rightarrow({\mathcal L}^{1})^{*}$$ et le diagramme commutatif

$$\xymatrix{0\ar[r]&\pi^{\#}({\mathcal G})\ar[r]\ar[d]_{\Xi_{\pi}({\mathcal G})}&\pi^{\#}(({\mathcal L}^{0})^{*})\ar[r]\ar[d]_{\Xi_{{\pi}(({\mathcal L}^{0})^{*})}}&\pi^{\#}(({\mathcal L}^{1})^{*}))\ar[d]_{\Xi_{{\pi}(({\mathcal L}^{1})^{*}))}}\\
0\ar[r]&\pi_{{\mathcal K}}^{!}({\mathcal G})\ar[r]&\pi_{{\mathcal K}}^{!}(({\mathcal L}^{0})^{*}))\ar[r]&\pi_{{\mathcal K}}^{!}(({\mathcal L}^{1})^{*}))}$$
La bijectivité de $\Xi_{\pi({\mathcal G})}$ découlent  de celle des deux autres flèches verticales. \vspace{1mm}

\noindent
Si ${\mathcal G}$ est un faisceau de profondeur au moins deux, on peut l'insérer dans la suite exacte courte
$$0\rightarrow{\mathcal G}\rightarrow{\mathcal G}^{**}\rightarrow {\mathcal H}\rightarrow 0$$
où ${\mathcal H}$ est faisceau sans torsion. Un diagramme commutatif analogue au précédent utilisant la bijectivité de $\Xi_{\pi^{\#}(({\mathcal G}^{**})}$ et l'injectivité de $\Xi_{\pi^{\#}(({\mathcal H})}$ permet de conclure$\,\blacksquare$
\end{proof}
\cor{}{}\label{c7}{} Avec les notations précédentes, on a, pour un morphisme de type $S_2$, les équivalences:\vspace{1mm}

\noindent
{\bf(i)}  $\pi_{{\mathcal K}}^{!}$ est exact à droite, \vspace{1mm}

\noindent
{\bf{(ii)}} ${\mathcal H}^{-n+j}(\pi^{!}{\mathcal G})=0$ pour tout $j>0$.
\rm
\begin{proof}
$\Longrightarrow$: Comme l'exactitude à droite du foncteur  $\pi_{{\mathcal K}}^{!}$ équivaut au fait que le morphisme $\pi$ est de Cohen Macaulay en vertu du  \theoremref{thm3} et que, relativement à toute  paramétrisation locale sur  $X$ donnée par  $\xymatrix{X\ar@/_/[rr]_{\pi}\ar[r]^{f}&Y\ar[r]^{q}&S}$, on a l'isomorphisme
$$f_{*}{\mathcal H}^{-n+j}(\pi^{!}{\mathcal G})\simeq{\mathcal Ext}^{j}(f_{*}{\mathcal O}_{X}, q^{*}{\mathcal G} \otimes_{{\mathcal
O}_{Y}}\Omega^{n}_{Y/S})$$
Or, de la locale liberté de $f_{*}{\mathcal O}_{X}$ résulte l'annulation des faisceaux ${\mathcal Ext}^{j}(f_{*}{\mathcal O}_{X}, q^{*}{\mathcal G} \otimes_{{\mathcal
O}_{Y}}\Omega^{n}_{Y/S})$ et ce, pour tout $j\geq 1$. D'où, 
${\mathcal H}^{-n+j}(\pi^{!}{\mathcal G})=0$.\vspace{1mm}

\noindent
$\Longleftarrow$: 
Le foncteur  ${\mathcal G}\rightarrow{\mathcal H}^{-n+1}(\pi^{!}{\mathcal G})$ est bien évidemment exacte à gauche et, par conséquent, le foncteur $\pi_{{\mathcal K}}^{!}$  est exact à droite.\vspace{1mm}

\noindent
Pour tout diagramme de $S$-plongement local,  $\xymatrix{X\ar@/_/[rr]_{\pi}\ar[r]^{\sigma}&S\times Z\ar[r]^{q}&S}$ où $Z$
 est une variété de Stein de dimension $n+p$, on a 
 $${\mathcal Ext}^{p+j}(\sigma_{*}{\mathcal O}_{X}, \Omega^{n+p}_{S\times Z/S})=0,\,\,\forall\,j\geq 1$$                                                          Comme $$\sigma_{*}{\mathcal H}^{-n+1}(\pi^{!}{\mathcal G})\simeq {\mathcal Ext}^{p+1}(\sigma_{*}{\mathcal O}_{X}, {\mathcal G} \otimes_{{\mathcal
O}_{S\times Z}}\Omega^{n+p}_{S\times Z/S})$$
Un raisonnement par récurrence sur la profondeur de ${\mathcal G}$ et en procédant par dévissage, on obtient l'annulation ${\mathcal Ext}^{p+1}(\sigma_{*}{\mathcal O}_{X}, {\mathcal G} \otimes_{{\mathcal
O}_{S\times Z}}\Omega^{n+p}_{S\times Z/S})=0$
pour tout faisceau cohérent ${\mathcal G}$                                                                                           
\noindent Signalons, au passage que les lemmes {\bf(5.2.2)} et {\bf(5.2.5)} de \cite{A.L} montrent que {\bf(ii)} impose au morphisme $\pi$ d'être de Cohen Macaulay$\,\blacksquare$

\end{proof}
\cor{}{}\label{c8}{} Si $\pi:X\rightarrow S$ est un ${\rm S}_{1}$-morphisme plat d'espaces complexes réduits, le foncteur $\pi_{{\mathcal K}}^{!}$ est exact à droite sur les faisceaux cohérents de profondeur au moins deux sur $S$.\rm\vspace{3mm}

\noindent
\begin{proof} C'est une conséquence du   \corollaryref{c4}$\,\blacksquare$ \end{proof}
\vspace{1mm}

\noindent


\vfill\eject
\section{\color{blue}{La dualité relative pour le foncteur dérivé du produit tensoriel pondéré et le foncteur $\pi^\#$.}}

\subsection{Le foncteur de Flenner.}\vspace{2mm}

\noindent
\begin{notation}
$\bullet$ Etant donné un morphisme $\pi:X\rightarrow S$ d'espaces analytiques complexes et un faisceau cohérent ${\mathcal G}$ plat sur $S$  dont le support est propre sur $S$, on désigne par:\vspace{1mm}

\noindent ${\rm I}\!{\rm L}\pi^{\#}:{\rm D}^{+}_{coh}(S)\rightarrow{\rm D}^{+}_{coh}(X)$ le foncteur covariant défini par ${\rm I}\!{\rm L}\pi^{\#}({\mathcal  N}^{\bullet}):={\mathcal G}\otimes^{{\rm I}\!{\rm L}} {\rm I}\!{\rm L}\pi^{*}{\mathcal N}^{\bullet}$ et par\vspace{1mm}

\noindent ${\mathcal T}_{\mathcal G}:{\rm D}_{coh}(X)\rightarrow {\rm D}_{coh}(X)$ le foncteur contravariant défini par  ${\mathcal A}^{\bullet}\rightarrow {\rm I}\!{\rm R}{\mathcal H}om({\mathcal A}^{\bullet}, {\mathcal G}\otimes^{{\rm I}\!{\rm L}} {\rm I}\!{\rm L}\pi^{*}{\mathcal D}^{\bullet}_{S})$, ${\mathcal D}^{\bullet}_{S}$ étant le complexe dualisant de Ramis-Ruget \cite{RR70} (${\Bbb D}_{S}$ sera le foncteur dualisant associé).\vspace{1mm}

\noindent
On note 

\end{notation}
Nous allons rappeler l'essentiel de \cite{Fle81} ({\bf Satz (2.1),p.178}) dans le: 
\Th{}{}\label{T'2} Soient  $\pi:X\rightarrow S$  un morphisme d'espaces complexes, ${\mathcal G}$ un  faisceau  cohérent   sur $X$, $S$-plat et de support propre sur $S$. Alors, \vspace{1mm}

\noindent
{\bf(i)} il existe un foncteur ${\rm I}\!{\rm L}\pi_{\#}: {\rm D}^{-}_{coh}(X)\rightarrow  {\rm D}^{-}_{coh}(S)$ dépendant de ${\mathcal G}$, covariant, commutant aux changements de base arbitraires et induisant, pour tout ${\mathcal F}^{\bullet}\in  {\rm D}^{-}_{coh}(X)$ et tout ${\mathcal N}^{\bullet}\in  {\rm D}^{+}_{coh}(S)$, un isomorphisme fonctoriel en les arguments,
$${\rm
I}\!{\rm R}\pi_{*}{\rm I}\!{\rm R}{\mathcal H}om({\mathcal F}^{\bullet}, {\rm I}\!{\rm L}\pi^{\#}({\mathcal N}^{\bullet}))\simeq{\rm I}\!{\rm R}{\mathcal H}om({\rm I}\!{\rm L}\pi_{\#}({\mathcal  F}^{\bullet}), {\mathcal N}^{\bullet})$$ 
{\bf(ii)} De plus, pour tout faisceau cohérent donné ${\mathcal F}$ sur $X$, il existe un complexe borné à droite à cohomologie cohérente ${\mathcal M}^{\bullet}$ tel que, pour tout faisceau cohérent ${\mathcal N}$ sur $S$ et tout entier $k$, on ait un isomorphisme
$${\mathcal E}xt^{k}(\pi; {\mathcal F}, {\mathcal G}\otimes \pi^{*}{\mathcal N})\simeq {\mathcal E}xt^{k}({\mathcal M}^{\bullet}, {\mathcal N})$$
\rm\vspace{2mm}

\noindent
\begin{notation} Dans la suite, nous appellerons pondération le faisceau ${\mathcal G}$ servant à définir les foncteurs ${\rm I}\!{\rm L}\pi_{\#}$,   ${\rm I}\!{\rm L}\pi^{\#}$ et nous dirons qu'un morphisme de foncteurs est un ${\mathcal G}$-morphisme s'il est défini relativement à cette ${\mathcal G}$- pondération.
\end{notation}
\noindent
\begin{rem}
{\bf(i)} Rappelons que la commutation aux changements de base plats signifie qu'étant donné un morphisme d'espaces complexes  $\nu:S'\rightarrow S$ de diagramme de changement de base est 
$$\xymatrix{X'\ar[d]_{\pi'}\ar[r]^{\theta}&X\ar[d]^{\pi}\\
 S'\ar[r]_{\nu}&S}$$
on a 
$${\rm I}\!{\rm L}\nu^{*}\circ{\rm I}\!{\rm L}\pi_{\#}={\rm I}\!{\rm L}{\pi'_{\#}}\circ{\rm I}\!{\rm L}\theta^{*}$$
{\bf(ii)} Par construction, on a
$${\rm I}\!{\rm L}\pi_{\#}({\mathcal F}^{\bullet}):= {\Bbb D}_{S}({\rm I}\!{\rm R}\pi_{*}{\mathcal T}_{\mathcal G}({\mathcal  F}^{\bullet}))$$
qui, vu la condition de propreté sur les supports, s'écrit aussi, grâce à la dualité relative de Ramis-Ruget-Verdier \cite{RRV71}, 
$${\rm I}\!{\rm L}\pi_{\#}({\mathcal F}^{\bullet})={\rm I}\!{\rm R}\pi_{*}{\Bbb D}_{X}({\mathcal T}_{\mathcal G}({\mathcal F}^{\bullet})).$$
  \end{rem}\vspace{2mm}
  
  \noindent
    En combinant les deux théorèmes de dualité  \cite{RRV71} et \cite{Fle81}, on peut  tirer ces quelques conséquences:
   \cor{}{} \label{co6} Soient $\pi: X\rightarrow S$ un morphisme propre et plat d'espaces complexes et ${\mathcal G}$ un faisceau cohérent $S$-plat sur $X$.  Alors, on a:
   $$ {\rm I}\!{\rm L}\pi_{\#}({\mathcal A}^{\bullet})={\rm I}\!{\rm R}\pi_{*}({\mathcal A}^{\bullet} \otimes^{{\rm I}\!{\rm L}}
   {\rm I}\!{\rm R}{\mathcal H}om({\mathcal G}, \pi^{!}({\mathcal O}_{S})))$$
 et, en particulier  pour  la ${\mathcal O}_{X}$-pondération, l'isomorphisme canonique
$${\rm I}\!{\rm R}\pi_{*}({\mathcal A}^{\bullet} \otimes^{{\rm I}\!{\rm L}}\pi^{!}{\mathcal B}^{\bullet})
\simeq {\rm I}\!{\rm L}\pi_{\#}({\mathcal A}^{\bullet})\otimes^{{\rm I}\!{\rm L}}{\mathcal B}^{\bullet}$$
\begin{proof} Rappelons que pour tout morphisme propre $\pi:X\rightarrow S$, la dualité relative de \cite{RRV71} induit, en particulier, la formule de projection bien connue
$${\rm I}\!{\rm R}\pi_{*}({\mathcal A}^{\bullet}\otimes^{{\rm I}\!{\rm L}}{\rm I}\!{\rm L}\pi^{*}{\mathcal B}^{\bullet})\simeq {\rm I}\!{\rm R}\pi_{*}({\mathcal A})^{\bullet}\otimes^{{\rm I}\!{\rm L}}{\mathcal B}^{\bullet}$$
Pour montrer la première égalité, on écrit
$$ {\rm I}\!{\rm D}_{S}({\rm I}\!{\rm L}\pi_{\#}({\mathcal A}^{\bullet}))={\rm I}\!{\rm R}\pi_{*}{\rm I}\!{\rm R}{\mathcal H}om({\mathcal A}^{\bullet}, {\mathcal G}\otimes\pi^{*}({\mathcal D}^{\bullet}_{S}))={\rm I}\!{\rm R}\pi_{*}{\rm I}\!{\rm R}{\mathcal H}om({\mathcal A}^{\bullet}\otimes^{{\rm I}\!{\rm L}} {\rm I}\!{\rm D}_{X}(({\mathcal G}\otimes\pi^{*}({\mathcal D}^{\bullet}_{S}))), {\mathcal D}^{\bullet}_{X})$$
qui s'écrit aussi, grâce à la dualité relative de \cite{RRV71},
$$ {\rm I}\!{\rm D}_{S}({\rm I}\!{\rm L}\pi_{\#}({\mathcal A}^{\bullet}))={\rm I}\!{\rm R}{\mathcal H}om({\rm I}\!{\rm R}\pi_{*}({\mathcal A}^{\bullet}\otimes^{{\rm I}\!{\rm L}} {\rm I}\!{\rm D}_{X}({\mathcal G}\otimes\pi^{*}({\mathcal D}^{\bullet}_{S})), {\mathcal D}^{\bullet}_{S})$$
et donc
$${\rm I}\!{\rm L}\pi_{\#}({\mathcal A}^{\bullet})={\rm I}\!{\rm R}\pi_{*}({\mathcal A}^{\bullet}\otimes^{{\rm I}\!{\rm L}} {\rm I}\!{\rm D}_{X}({\mathcal G}\otimes\pi^{*}({\mathcal D}^{\bullet}_{S}))$$
mais comme $\pi^{!}({\mathcal O}_{S})={\rm I}\!{\rm D}_{X}({\rm I}\!{\rm L}\pi^{*}({\mathcal D}^{\bullet}_{S})$, on a 
$${\rm I}\!{\rm D}_{X}({\mathcal G}\otimes\pi^{*}({\mathcal D}^{\bullet}_{S}))={\rm I}\!{\rm R}{\mathcal H}om(({\mathcal G}, \pi^{!}({\mathcal O}_{S})$$
et, par suite, la relation désirée
$$ {\rm I}\!{\rm L}\pi_{\#}({\mathcal A}^{\bullet})={\rm I}\!{\rm R}\pi_{*}({\mathcal A}^{\bullet} \otimes^{{\rm I}\!{\rm L}}
   {\rm I}\!{\rm R}{\mathcal H}om({\mathcal G}, \pi^{!}({\mathcal O}_{S}))$$
On obtient, en particulier, pour ${\mathcal G}:={\mathcal O}_{X}$ qui est un choix légitime puisque  $\pi$ est plat,  
$$ {\rm I}\!{\rm L}\pi_{\#}({\mathcal A}^{\bullet})={\rm I}\!{\rm R}\pi_{*}({\mathcal A}^{\bullet} \otimes^{{\rm I}\!{\rm L}}
   \pi^{!}({\mathcal O}_{S}))$$
   et, donc, 
$$ {\rm I}\!{\rm L}\pi_{\#}({\mathcal A}^{\bullet})\otimes^{{\rm I}\!{\rm L}}{\mathcal B}^{\bullet}={\rm I}\!{\rm R}\pi_{*}({\mathcal A}^{\bullet} \otimes^{{\rm I}\!{\rm L}}
   \pi^{!}({\mathcal O}_{S}))\otimes^{{\rm I}\!{\rm L}}{\mathcal B}^{\bullet}$$
qui, en vertu de la formule de projection, donne
$$ {\rm I}\!{\rm L}\pi_{\#}({\mathcal A}^{\bullet})\otimes^{{\rm I}\!{\rm L}}{\mathcal B}^{\bullet}=
{\rm I}\!{\rm R}\pi_{*}({\mathcal A}^{\bullet} \otimes^{{\rm I}\!{\rm L}}
   \pi^{!}({\mathcal O}_{S})\otimes^{{\rm I}\!{\rm L}}{\rm I}\!{\rm L}\pi^{*}{\mathcal B}^{\bullet})$$
On termine en utilisant le fait que la flèche canonique $$  \pi^{!}({\mathcal O}_{S})\otimes^{{\rm I}\!{\rm L}}{\rm I}\!{\rm L}\pi^{*}{\mathcal B}^{\bullet} \rightarrow\pi^{!}({\mathcal B}^{\bullet}) $$
est, ici,  un isomorphisme puisque $\pi$ est plat$\,\,\,\blacksquare$\end{proof} \vspace{1mm}

\noindent
\begin{rem}Pour montrer l'isomorphisme
$${\rm I}\!{\rm R}\pi_{*}({\mathcal A}^{\bullet} \otimes^{{\rm I}\!{\rm L}}\pi^{!}{\mathcal B}^{\bullet})
\simeq {\rm I}\!{\rm L}\pi_{\#}({\mathcal A}^{\bullet})\otimes^{{\rm I}\!{\rm L}}{\mathcal B}^{\bullet}$$
on peut aussi procéder comme suit en écrivant:
$${\rm I}\!{\rm R}{\mathcal H}om({\rm I}\!{\rm R}\pi_{*}({\mathcal A}^{\bullet} \otimes^{{\rm I}\!{\rm L}}\pi^{!}{\mathcal B}^{\bullet}), {\mathcal D}^{\bullet}_{S})\simeq{\rm I}\!{\rm R}\pi_{*}{\rm I}\!{\rm R}{\mathcal H}om({\mathcal A}^{\bullet} \otimes^{{\rm I}\!{\rm L}}\pi^{!}{\mathcal B}^{\bullet}, {\mathcal D}^{\bullet}_{X})$$
mais, par définition du foncteur $\pi^!$, on obtient
$${\rm I}\!{\rm R}\pi_{*}{\rm I}\!{\rm R}{\mathcal H}om({\mathcal A}^{\bullet} \otimes^{{\rm I}\!{\rm L}}\pi^{!}{\mathcal B}^{\bullet}, {\mathcal D}^{\bullet}_{X})\simeq{\rm I}\!{\rm R}\pi_{*}{\rm I}\!{\rm R}{\mathcal H}om({\mathcal A}^{\bullet}, {\rm I}\!{\rm L}\pi^{*}({\rm I}\!{\rm D}_{S}{\mathcal B}^{\bullet}))$$
qui, en vertu, de la dualité de Flenner, s'écrit aussi
$${\rm I}\!{\rm R}\pi_{*}{\rm I}\!{\rm R}{\mathcal H}om({\mathcal A}^{\bullet} \otimes^{{\rm I}\!{\rm L}}\pi^{!}{\mathcal B}^{\bullet}, {\mathcal D}^{\bullet}_{X})\simeq{\rm I}\!{\rm R}{\mathcal H}om( {\rm I}\!{\rm L}\pi_{\#}({\mathcal A}^{\bullet}), {\rm I}\!{\rm D}_{S}{\mathcal B}^{\bullet})\simeq{\rm I}\!{\rm R}{\mathcal H}om( {\rm I}\!{\rm L}\pi_{\#}({\mathcal A}^{\bullet})\otimes^{{\rm I}\!{\rm L}}{\mathcal B}^{\bullet}, {\mathcal D}^{\bullet}_{S})$$
D'où l'assertion.
\end{rem}
\cor{}{}\label{co7} La paire $({\rm I}\!{\rm L}\pi_{\#}, {\rm I}\!{\rm L}\pi^{\#})$ est adjointe et est  munie d'un morphisme canonique (à isomorphisme près) de foncteurs
  $${\rm I}\!{\rm L}\pi_{\#}{\rm I}\!{\rm L}\pi^{\#}\rightarrow {\rm Id}$$\rm
\begin{proof}
La propriété d'adjonction est bien évidemment une conséquence du \theoremref{T'2}.\vspace{1mm}

\noindent L'existence du morphisme annoncé découle du fait que  la donnée d'un morphisme
$${\rm I}\!{\rm L}\pi_{\#}{\rm I}\!{\rm L}\pi^{\#}({\mathcal B}^{\bullet})\rightarrow {\mathcal B}^{\bullet}$$
est équivalente à celle d'un morphisme
$${\Bbb D}_{S}({\mathcal B}^{\bullet})\rightarrow{\rm
I}\!{\rm R}\pi_{*}{\rm I}\!{\rm R}{\mathcal H}om({\rm I}\!{\rm L}\pi^{\#}({\mathcal B}^{\bullet}), {\rm I}\!{\rm L}\pi^{\#}({\mathcal D}_{S}^{\bullet})) $$
elle même, équivalente,  à celle de 
$${\rm I}\!{\rm L}\pi^{*}({\Bbb D}_{S}({\mathcal B}^{\bullet}))\rightarrow{\rm I}\!{\rm R}{\mathcal H}om({\rm I}\!{\rm L}\pi^{\#}({\mathcal B}^{\bullet}), {\rm I}\!{\rm L}\pi^{\#}({\mathcal D}_{S}^{\bullet})) $$
et,  donc, à celle de 
$${\rm I}\!{\rm L}\pi^{*}({\Bbb D}_{S}({\mathcal B}^{\bullet}))\otimes^{{\rm I}\!{\rm L}}{\rm I}\!{\rm L}\pi^{\#}({\mathcal B}^{\bullet})\rightarrow{\rm I}\!{\rm L}\pi^{\#}({\mathcal D}_{S}^{\bullet}) $$
qu'il est facile de construire en utilisant le morphisme naturel
$${\rm I}\!{\rm L}\pi^{*}({\Bbb D}_{S}({\mathcal B}^{\bullet}))\rightarrow{\rm I}\!{\rm R}{\mathcal H}om({\rm I}\!{\rm L}\pi^{*}({\mathcal B}^{\bullet}), {\rm I}\!{\rm L}\pi^{*}({\mathcal D}_{S}^{\bullet}))$$
et  l'accouplement de Yonéda
$${\rm I}\!{\rm L}\pi^{*}({\mathcal B}^{\bullet})\otimes^{{\rm I}\!{\rm L}} {\rm I}\!{\rm R}{\mathcal H}om({\rm I}\!{\rm L}\pi^{*}({\mathcal B}^{\bullet}), {\rm I}\!{\rm L}\pi^{*}({\mathcal D}_{S}^{\bullet}))\rightarrow {\rm I}\!{\rm L}\pi^{*}({\mathcal D}_{S}^{\bullet})$$
pour avoir
$$\underbrace{ {\rm I}\!{\rm L}\pi^{*}({\mathcal B}^{\bullet})\otimes^{{\rm I}\!{\rm L}} {\rm I}\!{\rm R}{\mathcal H}om({\rm I}\!{\rm L}\pi^{*}({\mathcal B}^{\bullet}), {\rm I}\!{\rm L}\pi^{*}({\mathcal D}_{S}^{\bullet}))}\otimes^{{\rm I}\!{\rm L}}{\mathcal G}\rightarrow{\rm I}\!{\rm L}\pi^{*}({\mathcal D}_{S}^{\bullet})\otimes^{{\rm I}\!{\rm L}}{\mathcal G}$$
On peut remarquer, au passage,  que le problème étant de nature locale et que les objets étant à cohomologie cohérente, on peut procéder à la construction fibre par fibre et, donc, se permettre de supposer le complexe dualisant à composantes injectives et travailler sur la catégorie des complexes $\blacksquare$
\end{proof}
\cor{}{}\label{co8} {\emph{Relation entre les deux dualités.}} \vspace{1mm}

\noindent Soit $\pi:X\rightarrow S$ un morphisme de ${\Bbb E}(S,n)$. Alors, relativement à la $\omega^{n}_{X/S}$-pondération, on a des morphismes canoniques de foncteurs
$${\rm I}\!{\rm R}\pi_{*}[n]\rightarrow{\rm I}\!{\rm L}\pi_{\#},\,\,\,{\rm I}\!{\rm L}\pi^{\#}
\rightarrow\pi^{!}[-n]$$
dont la donnée de l'un est équivalente à la donnée de l'autre. De plus, on a
un diagramme commutatif
$$\xymatrix{{\rm
I}\!{\rm R}\pi_{*}{\rm I}\!{\rm R}{\mathcal H}om({\mathcal F}^{\bullet}, {\mathcal G}\otimes^{{\rm I}\!{\rm L}} {\rm I}\!{\rm L}\pi^{*}{\mathcal N}^{\bullet})\ar[d]\eq[r]&{\rm I}\!{\rm R}{\mathcal H}om({\rm I}\!{\rm L}\pi_{\#}({\mathcal  F}^{\bullet}), {\mathcal N}^{\bullet})\ar[d]\\
{\rm I}\!{\rm R}\pi_{*}{\rm I}\!{\rm R}{\mathcal H}om({\mathcal F}^{\bullet}, \pi^{!}({\mathcal N}^{\bullet})[-n])\eq[r]&{\rm I}\!{\rm R}{\mathcal H}om({\rm I}\!{\rm R}\pi_{*}{\mathcal F}^{\bullet}[n], {\mathcal N}^{\bullet}))}$$
\begin{proof}
La pondération étant donné par le faisceau cohérent $\omega^{n}_{X/S}$ de support $X$, le morphisme est automatiquement plat. Alors, la formule de projection pour un morphisme propre $${\rm I}\!{\rm R}\pi_{*}(\omega^{n}_{X/S}\otimes \pi^{*}{\mathcal D}^{\bullet}_{S})\simeq {\rm I}\!{\rm R}\pi_{*}(\omega^{n}_{X/S})\otimes{\mathcal D}^{\bullet}_{S}$$ et le morphisme de type trace ou d'intégration ${\rm I}\!{\rm R}\pi_{*}\omega^{n}_{X/S}[n]\rightarrow {\mathcal O}_{S}$
nous donnent un morphisme naturel 
$${\rm I}\!{\rm R}\pi_{*}{\mathcal T}_{\mathcal G}({\mathcal F}^{\bullet})\rightarrow {\Bbb D}_{S}({\rm I}\!{\rm R}\pi_{*}({\mathcal F}^{\bullet})[n])$$ et, donc, en dualisant sur $S$, le morphisme canonique
$${\rm I}\!{\rm R}\pi_{*}{\mathcal F}^{\bullet}[n]\rightarrow{\rm I}\!{\rm L}\pi_{\#}({\mathcal F}^{\bullet})$$
et donc le morphisme de foncteurs
$${\rm I}\!{\rm R}\pi_{*}[n]\rightarrow{\rm I}\!{\rm L}\pi_{\#}$$
Le morphisme ${\rm I}\!{\rm L}\pi^{\#}\rightarrow \pi^{!}[-n]$, quant à lui, relève  de la définition $\omega^{n}_{X/S}:={\mathcal H}^{0}(\pi^{!}({\mathcal O}_{S})[-n])$ et de la composition des flèches naturelles suivantes:
$${\rm I}\!{\rm L}\pi^{\#}({\mathcal N}^{\bullet}):=\omega^{n}_{X/S}\otimes^{{\rm I}\!{\rm L}}\pi^{*}({\mathcal N}^{\bullet})\rightarrow\pi^{!}({\mathcal O}_{S})[-n]\otimes^{{\rm I}\!{\rm L}}\pi^{*}({\mathcal N}^{\bullet})\rightarrow\pi^{!}({\mathcal N}^{\bullet})[-n]$$
puisque, pour tout complexe de
faisceaux à cohomologie cohérente ${\mathcal A}^{\bullet}$ et
${\mathcal B}^{\bullet}$, on a un morphisme canonique
$${\rm I}\!{\rm L}\pi^{*}{\mathcal A}^{\bullet}\otimes^{{\rm I}\!{\rm L}} \pi^{!}{\mathcal B}^{\bullet}\rightarrow \pi^{!}({\mathcal A}^{\bullet}\otimes^{{\rm I}\!{\rm L}} {\mathcal B}^{\bullet})$$
qui est un isomorphisme $\pi$ est de {\emph Tordim}-finie ou ${\mathcal A}^{\bullet}$ borné à composantes  plates (ou simplement cohomologiquement
plat!)).\vspace{1mm}

\noindent L'existence du diagramme 
$$\xymatrix{{\rm
I}\!{\rm R}\pi_{*}{\rm I}\!{\rm R}{\mathcal H}om({\mathcal F}^{\bullet}, {\mathcal G}\otimes^{{\rm I}\!{\rm L}} {\rm I}\!{\rm L}\pi^{*}{\mathcal N}^{\bullet})\ar[d]\eq[r]&{\rm I}\!{\rm R}{\mathcal H}om({\rm I}\!{\rm L}\pi_{\#}({\mathcal  F}^{\bullet}), {\mathcal N}^{\bullet})\ar[d]\\
{\rm I}\!{\rm R}\pi_{*}{\rm I}\!{\rm R}{\mathcal H}om({\mathcal F}^{\bullet}, \pi^{!}({\mathcal N}^{\bullet})[-n])\eq[r]&{\rm I}\!{\rm R}{\mathcal H}om({\rm I}\!{\rm R}\pi_{*}{\mathcal F}^{\bullet}[n], {\mathcal N}^{\bullet}))}$$
et sa commutativité découle de ce qui précède et des isomorphismes de dualité de \cite{Fle81} et de \cite{RRV71}

$\,\blacksquare$
\end{proof}

\centerline{\color{blue}{La preuve du \theoremref{thm2}}}
\rm\vspace{2mm}

\noindent
Comme pour le \theoremref{thm1}, nous allons d'abord mettre en évidence l'existence de l'adjoint du foncteur $\pi^{\#}$ et la dualité induite puis nous intéresser aux propriétés fonctorielles.
\Prop{}{}\label{P9}\vspace{1mm}
Soit $\pi:X\rightarrow S$ un morphisme de ${\Bbb E}(S,n)$ pour lequel $\omega^{n}_{X/S}$ est plat sur $S$. Alors,  il existe un unique foncteur(à isomorphisme près) $\pi_{\#}:{\rm Coh}(X)\rightarrow{\rm Coh}(S)$  covariant, exact à droite  tel que la paire $({\pi}_{\#}, \pi^{\#})$  soit  muni d'un morphisme de foncteurs
$${\mathcal T}_{\pi}:{\pi}_{\#}\pi^{\#}\rightarrow {\rm Id}$$ 
et  induise, pour tout faisceaux cohérents ${\mathcal F}$ et  ${\mathcal G}$ sur $X$ et $S$ respectivement, un isomorphisme canonique
$${\rm I}\!{\rm H}om_{X}({\mathcal F}, \pi^{\#}({\mathcal N}))\rightarrow{\rm I}\!{\rm H}om_{S}({\pi}_{\#}{\mathcal F}, {\mathcal N}) $$
 fonctoriel en ces arguments dont la formation commute aux changements de base plats et aux restriction ouvertes sur $X$. De plus, on a le diagramme commutatif
$$\xymatrix{\pi_{*}{\mathcal H}om({\mathcal F}, \omega^{n}_{X/S}\otimes {\mathcal N})\ar[d]\eq[r]&{\mathcal H}om( \pi_{\#}{\mathcal F}, {\mathcal N})\ar[d]\\
\pi_{*}{\mathcal H}om({\mathcal F}, \pi^{!}_{\mathcal K}{\mathcal N})\eq[r]&{\mathcal H}om({\rm I}\!{\rm R}^{n}\pi_{*}{\mathcal F}, {\mathcal N})}$$
 \begin{proof} Comme nous l'avons souligné dans le \corollaryref{co8}, la dualité relative générale de Flenner pour la $\omega^{n}_{X/S}$-pondération 
$${\rm
I}\!{\rm R}\pi_{*}{\rm I}\!{\rm R}{\mathcal H}om({\mathcal F}^{\bullet},  \omega^{n}_{X/S}\otimes^{{\rm I}\!{\rm L}} {\rm I}\!{\rm L}\pi^{*}{\mathcal N}^{\bullet})\simeq{\rm I}\!{\rm R}{\mathcal H}om({\rm I}\!{\rm L}\pi_{\#}({\mathcal  F}^{\bullet}), {\mathcal N}^{\bullet})$$ 
montre que la paire $({\rm I}\!{\rm L}\pi_{\#},{\rm I}\!{\rm L}\pi^{\#})$ est adjointe.  On va donc montrer qu'il en est de même pour les foncteurs $\pi_{\#}:={\mathcal H}^{0}({\rm I}\!{\rm L}\pi_{\#})$ et   $\pi^{\#}:={\mathcal H}^{0}({\rm I}\!{\rm L}\pi^{\#})$ dans la catégorie des faisceaux cohérents.\vspace{1mm}

\noindent
Comme les arguments sont des faisceaux cohérents, les annulations  ${\mathcal H}^{j}({\rm I}\!{\rm L}\pi_{\#}({\mathcal  F}))=0$  (resp. ${\mathcal H}^{j}({\mathcal N})=0$) pour $j>0$ (resp. $j<0$) et les  arguments avancés dans la preuve de la \propositionref{P5},  nous donnent l'isomorphisme  (cf \cite{Fle81}, \S 3)
$$\pi_{*}{\mathcal H}om({\mathcal F}, \omega^{n}_{X/S}\otimes {\mathcal N})\simeq{\mathcal H}om( \pi_{\#}{\mathcal F}, {\mathcal N}),\,\,{\mathcal F}\in {\rm Coh}(X),\,\,\,{\mathcal N}\in {\rm Coh}(S)$$ 
Le  morphisme de foncteurs $\pi_{\#}\pi^{\#}\rightarrow{\rm Id}$
rendant commutatif le diagramme
$$\xymatrix{ {\rm I}\!{\rm L}\pi_{\#}\circ {\rm I}\!{\rm L}\pi^{\#}\ar[rr]^{\beta}\ar[rd]_{\alpha}&&\pi_{\#}\pi^{\#}\ar[ld]^{\gamma}\\ 
&{\rm Id}&}$$
où $\alpha$ est le morphisme défini dans le \corollaryref{co7} et $\beta$ déduit des flèches naturelles
$${\rm I}\!{\rm L}\pi_{\#}\circ {\rm I}\!{\rm L}\pi^{\#}\rightarrow{\rm I}\!{\rm L}\pi_{\#}({\mathcal H}^{0}({\rm I}\!{\rm L}\pi_{\#}))\rightarrow{\mathcal H}^{0}({\rm I}\!{\rm L}\pi_{\#})({\mathcal H}^{0}({\rm I}\!{\rm L}\pi_{\#})=\pi_{\#}\pi^{\#}$$
Pour la construction de  $\gamma$, on procède comme dans le corollaire précédemment cité  mais au niveau des faisceaux.\vspace{1mm}

\noindent
Ainsi, on définira, pour tout faisceau cohérent sur $S$, un morphisme $\pi_{\#}\pi^{\#}{\mathcal F}\rightarrow{\mathcal F}$ si l'on peut définir un morphisme
${\mathcal H}om({\mathcal F}, {\mathcal O}_{S})\rightarrow{\mathcal H}om(\pi_{\#}\pi^{\#}{\mathcal F}, {\mathcal O}_{S}) $. Or, ce dernier revient à se donner un morphisme ${\mathcal H}om({\mathcal F}, {\mathcal O}_{S})\rightarrow\pi_{*}{\mathcal H}om(\pi^{\#}{\mathcal F}, \pi^{\#}{\mathcal O}_{S})$ ou de manière équivalente
$\pi^{*}{\mathcal H}om({\mathcal F}, {\mathcal O}_{S})\rightarrow{\mathcal H}om(\pi^{\#}{\mathcal F},\omega^{n}_{X/S} )$ c'est-à-dire, en raison de la platitude de $\pi$, un morphisme ${\mathcal H}om(\pi^{*}{\mathcal F}, {\mathcal O}_{X})\rightarrow{\mathcal H}om(\pi^{\#}{\mathcal F},\omega^{n}_{X/S} )$
et, donc, un morphisme
$${\mathcal H}om(\pi^{*}{\mathcal F}, {\mathcal O}_{X})\otimes\pi^{*}{\mathcal F}\otimes\omega^{n}_{X/S}\rightarrow\omega^{n}_{X/S} $$
qui sera donné, là encore, par l'accouplement de Yonéda.
 De plus, on a un  diagramme commutatif reliant les deux dualités 
$$\xymatrix{\pi_{*}{\mathcal H}om({\mathcal F}, \omega^{n}_{X/S}\otimes {\mathcal N})\ar[d]\eq[r]&{\mathcal H}om( \pi_{\#}{\mathcal F}, {\mathcal N})\ar[d]\\
\pi_{*}{\mathcal H}om({\mathcal F}, \pi^{!}_{\mathcal K}{\mathcal N})\eq[r]&{\mathcal H}om({\rm I}\!{\rm R}^{n}\pi_{*}{\mathcal F}, {\mathcal N})}$$
En prenant, ${\mathcal N}={\mathcal O}_{S}$, on en déduit que le morphisme canonique $ {\rm I}\!{\rm R}^{n}\pi_{*}\rightarrow \pi_{\#}$ obtenu en prenant la cohomologie de degré $0$ du morphisme ${\rm I}\!{\rm R}\pi_{*}[n]\rightarrow{\rm I}\!{\rm L}\pi_{\#}$  induit l'isomorphisme
$${\mathcal H}om({\rm I}\!{\rm R}^{n}\pi_{*}{\mathcal F}, {\mathcal O}_{S})\simeq {\mathcal H}om(\pi_{\#}{\mathcal F}, {\mathcal O}_{S})$$
\end{proof}
\begin{rem}
On peut remarquer qu'une condition nécéssaire pour que le foncteur ${\mathcal G}\rightarrow\omega^{n}_{X/S}\otimes {\mathcal N}$ admette un adjoint à gauche ${\rm T}$, est que $\omega^{n}_{X/S}$ soit plat sur $S$ (ce qui impose la platitude de $\pi$!). En effet, un isomorphisme du type
$$\pi_{*}{\mathcal H}om({\mathcal F}, \omega^{n}_{X/S}\otimes {\mathcal N})\simeq{\mathcal H}om( {\rm T}({\mathcal F}), {\mathcal N})$$
contraint le foncteur $\pi^{\#}$ d'être exact à gauche; ce qui ne peut avoir lieu que si $\omega^{n}_{X/S}$ est $S$-plat.
\end{rem}

\Prop{}{}\label{P10}\vspace{1mm}
Cette dualité satisfait toutes les propriétés fonctorielles énoncées dans la \propositionref{P9} à savoir
la formation du foncteur $\pi_{\#}$ satisfait les propriétés suivantes: \vspace{2mm}

\noindent
 {\bf(i)  compatibilité avec la
composition des morphismes universellement équidimensionnels:}
  c'est-à-dire que pour tout diagramme commutatif $$\xymatrix{X\ar[rr]^{f}\ar[rd]_{\pi}&&Y\ar[ld]^{g}\\
&S&}$$ de $S$-espaces complexes avec $f\in
{\Bbb E}(Y; m)$ et $g\in{\Bbb E}(S; n)$, on a
$$\pi_{\#}= f_{\#}\circ g_{\#}$$
 {\bf(ii)  compatibilité avec les morphismes de changement de bases:}
Elle commute à tout changement de base  et aux restrictions ouvertes sur $X$ c'est-à-dire que  pour tout morphisme d'espaces complexes  $\nu:S'\rightarrow S$ de  diagramme de changement de base induit 
 $$\xymatrix{X'\ar[d]_{\pi'}\ar[r]^{\theta}&X\ar[d]^{\pi}\\
 S'\ar[r]_{\nu}&S}$$
 on a: $$\Theta^{*}\pi_{\#} =\pi'_{\#}\nu^{*}$$
et pour tout ouvert $U$ de $X$ muni d'un morphisme d'inclusion $j_U$ et de diagramme associé
$$\xymatrix{U\ar[d]_{\pi_U}\ar[r]^{j_U}&X\ar[d]^{\pi}\\
 V:=\pi(U)\ar[r]_{j_V}&S}$$
on a, $$j^{*}_U\circ (\pi_{\#})=({\pi_{U}})_{\#}\circ j^{*}_V$$
{\bf(iii) elle est compatible aux images directes
propres} entre éléments de ${\Bbb E}(S; n)$  c'est-à-dire que, pour
tout diagramme commutatif d'espaces complexes sur $S$,
 $$\xymatrix{{X'}\ar[rr]_{\theta}\ar[rd]_{{\pi'}}&&X\ar[ld]^{\pi}\\
&S&}$$ dans lequel  ${\pi}$ (resp. $\pi'$) est propre à
fibres de
 dimension pure $n$, ${\theta}$ est propre ( donc nécessairement génériquement
 fini), on a un morphisme canonique
 $$\theta_{*}{\tilde{\pi}}_{\#}\rightarrow \pi_{\#}$$ \rm
 \begin{proof}
{\bf (i) la  compatibilité avec la
composition des morphismes universellement équidimensionnels}.\vspace{2mm}

\noindent
Soit  $$\xymatrix{X_{2}\ar[rr]^{f}\ar[rd]_{\pi}&&X_{1}\ar[ld]^{g}\\
&S&}$$ un diagramme commutatif de $S$-espaces complexes avec $f\in
{\Bbb E}(X_{1}; m)$ et $g\in{\Bbb E}(S; n)$. On suppose donné un faisceau cohérent ${\mathcal G}_{1}$ (resp.${\mathcal G}_{2}$) de support $X_{1}$ (resp. $X_{2}$) $S$-plat.  \vspace{1mm}

\noindent
On veut montrer que le morphisme composé $\pi$ (qui est un élément de ${\Bbb E}(S; m+n)$) est naturellement muni d'un faisceau cohérent $S$-plat et satisfait la relation
$${\pi}_{\#}= g_{\#}\circ f_{\#}$$
Pour cela, on va appliquer plusieurs fois la dualité de la  \propositionref{P9}.\vspace{1mm}

\noindent
On a les isomorphismes
$${\rm I}\!{\rm H}om({S}; g_{\#}(f_{\#}{{\mathcal A}}), {\mathcal B})\simeq
{\rm I}\!{\rm H}om(X_{1}; (f_{\#}{{\mathcal A}}), g^{\#}({\mathcal B}))\simeq{\rm I}\!{\rm H}om(X_{2}; {\mathcal A},  f^{\#}(g^{\#}({\mathcal B})))$$
Or,
$$f^{\#}(g^{\#}({\mathcal B}))={\mathcal G}_{2}\otimes f^{*}(g^{\#}({\mathcal B}))={\mathcal G}_{2}\otimes f^{*}({\mathcal G}_{1}\otimes g^{*}({\mathcal B}))\simeq{\mathcal G}_{2}\otimes f^{*}({\mathcal G}_{1})\otimes f^{*}(g^{*}({\mathcal B}))$$
Alors, comme $f$ est plat et que ${\mathcal G}_{1}$ est $S$-plat, $f^{*}({\mathcal G}_{1})$ est aussi $S$-plat.  Il s'en suit que le faisceau cohérent ${\mathcal G}_{2}\otimes f^{*}({\mathcal G}_{1})$ est $S$-plat puisque ${\mathcal G}_{2}$ est $X_{1}$-plat (ces vérifications ne présentent aucune difficultés car il suffit de s'assurer de la préservation de l'injectivité ou de l'exactitude à gauche des suites exactes courtes). Ainsi, il nous est possible de choisir ${\mathcal G}_{2}\otimes f^{*}({\mathcal G}_{1})$ pour  pondération dans la définition de ${\rm I}\!{\rm L }\pi_{\#}$ et donc de $\pi_{\#}$ pour aboutir à la relation voulue.\vspace{1mm}

\noindent
{\bf (ii) Propriétés de compatibilité relativement aux morphismes de changement de bases.} \vspace{1mm}

\noindent
 Soit $\nu:S'\rightarrow S$ un morphisme d'espaces complexes réduits arbitraire et
$$\xymatrix{X'\ar[d]_{\pi'}\ar[r]^{\theta}&X\ar[d]^{\pi}\\
 S'\ar[r]_{\nu}&S}$$
le diagramme de changement de base associé. 
Alors, on a un isomorphisme canonique de foncteurs
$$\nu^{*}\circ\pi_{\#}\rightarrow\pi'_{\#}\circ\theta^{*}$$
Commençons par signaler que si ${\mathcal G}$ est un faisceau $S$-plat, alors $\theta^{*}{{\mathcal G}}$ est naturellement  aussi $S'$-plat. De plus, comme, en ce qui nous concerne, le support des pondérations est toujours l'espace source tout entier, les morphismes sont automatiquement plats. Ainsi, $\pi$ étant plat, $\pi'$ l'est aussi.\vspace{1mm}

\noindent
On alors le diagramme commutatif
$$\xymatrix{{\rm I}\!{\rm H}om({S'}; \nu^{*}\pi_{\#}{{\mathcal A}}, {\mathcal B})\ar[r]^{\alpha}\ar[d]_{\psi}&{\rm I}\!{\rm H}om(S; \pi_{\#}{{\mathcal A}}, \nu_{*}{\mathcal B})\ar[d]_{\gamma }\ar[r]^{\beta}&
{\rm I}\!{\rm H}om({X}; {\mathcal A}, 
{\mathcal G}\otimes\pi^{*}(\nu_{*}{\mathcal B}))\ar[d]_{\delta}\\
\!{\rm H}om(S'; \pi'_{\#}(\theta^{*}{\mathcal A}),{\mathcal B})& 
{\rm I}\!{\rm H}om({X'}; \theta^{*}{\mathcal A}, 
\theta^{*}{\mathcal G}\otimes\pi'^{*}({\mathcal B}))\ar[l]_{\alpha'}&{\rm I}\!{\rm H}om({X}; {\mathcal A}, 
{\mathcal G}\otimes\theta_{*}\pi'^{*}({\mathcal B}))\ar[l]_{\beta'}}$$
dans lequel $\alpha$, $\beta$, $\delta$ et $\alpha'$ sont des isomorphismes avec $\delta$ déduit de la commutation de l'image directe avec le changement de base plat car $\theta$ est déduit de $\nu$ dans le changement de base donné par $\pi$.  Le morphisme $\beta'$ résulte de la composition
$${\mathcal G}\otimes\theta_{*}\pi'^{*}({\mathcal B})\rightarrow\theta_{*}\theta^{*}{\mathcal G}\otimes\theta_{*}\pi'^{*}({\mathcal B})\rightarrow\theta_{*}(\theta^{*}{\mathcal G}\otimes\pi'^{*}({\mathcal B})$$
et de la formule d'adjonction.\vspace{1mm}
\noindent
Cela montre  qu'il existe un morphisme de foncteurs
$$\pi'_{\#}\circ\theta^{*}\rightarrow \nu^{*}\circ\pi_{\#}$$
qui est un isomorphisme si $\nu$ est propre puisque, dans ce cas, $\beta'$ est un isomorphisme (formule de projection et d'adjonction).\vspace{1mm}

\noindent
Pour s'assurer de sa bijectivité dans le cas général, on raisonne localement sur $S$ en décomposant $\nu$ en un plongement suivi d'une projection (cf \cite{Fle81}, \S 3). Ainsi, pour un plongement local  $\nu:=\sigma$, $\theta$ sera aussi un plongement et dans ce cas $\beta'$ est un isomorphisme (c'est la formule de projection!). Il n'est pas plus difficile de l'établir dans  le cas d'une projection. 

\vspace{1mm}

\noindent La compatibilité aux restrictions ouvertes sur:\vspace{1mm}

\indent
$\star$ $S$ se déduit facilement de la propriété précédente de commutation avec les changements de base (d'ailleurs l'inclusion ouverte est plate), \vspace{1mm}

\indent $\star$ $X$, elle découle de l'ouverture de $\pi$ puisque comme  l'image d'un ouvert est un ouvert, on peut recouvrir tout ouvert de $X$ par une réunion localement finie d'ouverts saturés c'est-à-dire du type $\pi^{-1}(\pi(U))$ ce qui nous  ramène au cas précédent  
On a, alors, pour un ouvert $U$  de $X$ muni de l'immersion ouverte  $j_{U}:U\rightarrow X$, un isomorphisme naturel
${j_{U}}^{*}\pi_{\#}() ={\pi_{U}}_{\#}({j_{U}}^{*})$ avec $\pi_{U}=\pi|_{U}$.
 \vspace{1mm}

\noindent
{\bf(iii) la compatibilité aux images directes
propres entre éléments de}  ${\Bbb E}(S; n)$.\vspace{1mm}

\noindent
$\bullet${\bf{Sur la source:} }On se donne un diagramme commutatif d'espaces complexes sur $S$,
 $$\xymatrix{{X'}\ar[rr]^{\theta}\ar[rd]_{{\pi'}}&&X\ar[ld]^{\pi}\\
&S&}$$ dans lequel  ${\pi}$ (resp. $\pi'$) est propre à
fibres de
 dimension pure $n$, ${\theta}$ est propre ( donc nécessairement génériquement
 fini), on a un morphisme canonique
 $${\rm I}\!{\rm L}\pi'_{\#}({\rm I}\!{\rm L}\theta^{*}({\mathcal A}^{\bullet}))\rightarrow {\rm I}\!{\rm L}\pi_{\#}({\mathcal A}^{\bullet})$$ 
 et un isomorphisme
 $$\pi'_{\#}(\theta^{*}({\mathcal A}))\rightarrow\pi_{\#}({\mathcal A})$$
 \indent $\bullet$ Nous allons d'abord montrer que:
 $${\mathcal G'}\,\,\,S-{\rm plat}\Longrightarrow\,\theta_{*}({\mathcal G'})\,\,\,S-{\rm plat}$$\vspace{1mm}
 
 \noindent
 Pour s'en convaincre, il nous suffit de montrer que toute injection  $S$,
 $$0\rightarrow {\mathcal A}\rightarrow {\mathcal B}$$
 de faisceaux cohérents sur $S$ est transformée en une injection
$$0\rightarrow \theta_{*}({\mathcal G'})\otimes\pi^{*}({\mathcal A})\rightarrow  \theta_{*}({\mathcal G'})\otimes\pi^{*}({\mathcal B})$$
 Or, ceci est une conséquence immédiate de la formule de projection et de l'exactitude à gauche du foncteur $\theta_{*}$. En effet, par hypothèse, on a la suite exacte
 $$0\rightarrow {\mathcal G'}\otimes\pi'^{*}({\mathcal A})\rightarrow {\mathcal G'}\otimes\pi'^{*}({\mathcal B})$$
 à laquelle on applique $\theta_{*}$ pour avoir, grâce à la formule de projection et l'exactitude à gauche
 $$0\rightarrow \theta_{*}({\mathcal G'})\otimes\pi^{*}({\mathcal A})\rightarrow  \theta_{*}({\mathcal G'})\otimes\pi^{*}({\mathcal B})$$
 qui prouve la $S$-platitude de $\theta_{*}({\mathcal G'})$.
 \vspace{1mm}
 
 \noindent
Alors, 
 $$ {\rm I}\!{\rm L}\pi'_{\#}({\rm I}\!{\rm L}\theta^{*}{\mathcal A}^{\bullet})={\rm I}\!{\rm R}\pi'_{*}({\rm I}\!{\rm L}\theta^{*}({\mathcal A}^{\bullet}) \otimes^{{\rm I}\!{\rm L}}
   {\rm I}\!{\rm R}{\mathcal H}om({\mathcal G'}, \pi'^{!}({\mathcal O}_{S}))={\rm I}\!{\rm R}\pi_{*}{\rm I}\!{\rm R}\theta_{*}({\rm I}\!{\rm L}\theta^{*}({\mathcal A}^{\bullet}) \otimes^{{\rm I}\!{\rm L}}
   {\rm I}\!{\rm R}{\mathcal H}om({\mathcal G'}, \pi'^{!}({\mathcal O}_{S})))$$
 D'où, en vertu de la formule de projection
 $$ {\rm I}\!{\rm L}\pi'_{\#}({\rm I}\!{\rm L}\theta^{*}{\mathcal A}^{\bullet})={\rm I}\!{\rm R}\pi_{*}({\mathcal A}^{\bullet}\otimes^{{\rm I}\!{\rm L}}{\rm I}\!{\rm R}\theta_{*} {\rm I}\!{\rm R}{\mathcal H}om({\mathcal G'}, \pi'^{!}({\mathcal O}_{S})))$$
 et de la dualité pour un morphisme propre
 $$ {\rm I}\!{\rm L}\pi'_{\#}({\rm I}\!{\rm L}\theta^{*}{\mathcal A}^{\bullet})={\rm I}\!{\rm R}\pi_{*}({\mathcal A}^{\bullet}\otimes^{{\rm I}\!{\rm L}} {\rm I}\!{\rm R}{\mathcal H}om({\rm I}\!{\rm R}\theta_{*}{\mathcal G'}, \pi^{!}({\mathcal O}_{S})))$$
 et, donc, un morphisme canonique
 $$ {\rm I}\!{\rm L}\pi'_{\#}({\rm I}\!{\rm L}\theta^{*}{\mathcal A}^{\bullet})\rightarrow{\rm I}\!{\rm L}\pi_{\#}({\mathcal A}^{\bullet})$$
 ${\rm I}\!{\rm L}\pi_{\#}$ étant défini relativement à la pondération $\theta_{*}{\mathcal G'}$. On remarquera, au passage, que pour $\theta$ fini, ce morphisme est un ismorphisme.\vspace{1mm}
 
 \noindent
 Au niveau des faisceaux, on a un isomorphisme canonique
 $$ \pi'_{\#}(\theta^{*}{\mathcal A})\simeq\pi_{\#}({\mathcal A})$$
 pour tout faisceau cohérent ${\mathcal A}$ sur $X$. \vspace{1mm}
 
 \noindent
 En effet, il nous suffit d'utiliser les  isomorphismes de dualité donnés par la \propositionref{P9} et l'adjonction naturelles pour avoir:
 $$\xymatrix{{\rm I}\!{\rm H}om(S; {\pi'_{\#}}(\theta^{*}{\mathcal A}), {\mathcal B})\ar[d]_{\alpha_{0}}\eq[r]^{\alpha_{1}}&{\rm I}\!{\rm H}om({X'}; \theta^{*}{\mathcal A},
 {\pi'}^{\#}({\mathcal B}))\eq[r]^{\alpha_{2}}
 &{\rm I}\!{\rm H}om(X'; {\theta}^{*}
 {\mathcal A}, {\mathcal G'}\otimes \pi'^{*}{\mathcal B})\eq[d]^{\alpha_{3}}\\
{\rm I}\!{\rm H}om(S; \pi_{\#}({\mathcal A}), {\mathcal B})&
{\rm I}\!{\rm H}om(X; {\mathcal A},  \theta_{*}{\mathcal G'}\otimes \pi^{*}{\mathcal B})\eq[l]^{\alpha_{5}}&{\rm I}\!{\rm H}om(X; {\mathcal A}, \theta_{*}({\mathcal G'}\otimes \theta^{*}(\pi^{*}{\mathcal B})))\eq[l]^{\alpha_{4}}}$$ 
et donc l'isomorphisme voulue.\vspace{2mm}

\noindent
$\bullet$ {\bf{Sur la base}}:\vspace{2mm}

\noindent
On se donne un diagramme commutatif d'espaces complexes 
$$\xymatrix{X\ar[d]_{\pi}\ar[rd]^{\pi'}&\\
 S\ar[r]_{\nu}&S'}$$
 dans lequel $\pi$ est un  morphisme de ${\Bbb E}(S; n)$  muni d'un faisceau $S$-plat ${\mathcal G}$ et $\nu$ plat. Alors, ${\mathcal G}$ est $S'$-plat et on a un morphisme naturel $${\pi}_{\#}({\mathcal A})\rightarrow\nu^{*}({\pi'}_{\#}({\mathcal
A}))$$
dont la construction résulte du diagramme commutatif
$$\xymatrix{{\rm I}\!{\rm H}om(S; \nu^{*}({\pi'}_{\#}({\mathcal
A})), {\mathcal B}) \ar[d]\eq[r]&{\rm I}\!{\rm H}om({S'}; {\pi'}_{\#}({\mathcal A}), \nu_{*}{\mathcal B})
)\eq[r]&{\rm I}\!{\rm H}om({S'};\nu_{\#}({\pi}_{\#}({\mathcal
A})), \nu_{*}{\mathcal B}))\eq[d]\\
{\rm I}\!{\rm H}om(S;{\pi}_{\#}({\mathcal
A})), {\mathcal B}))&&
{\rm I}\!{\rm H}om(S;{\pi}_{\#}({\mathcal
A})), \nu^{\#}(\nu_{*}{\mathcal B}))\ar[ll]}$$
dans lequel les flèches sont données par les isomorphismes d'adjonction, de dualité de la \propositionref{P9}, de l'égalité $\nu^{\#}=\nu^{*}$ pour la pondération triviale pour $\nu$ et du moprhisme canonique $\nu^{*}\nu_{*}\rightarrow {\rm Id}$.
$\blacksquare$\end{proof}
\vspace{2mm}
\noindent

\section{\color{blue}{ Propriétés du faisceau $\pi^{\#}({\mathcal G}):=\pi^{*}({\mathcal G})\otimes \omega^{n}_{X/S}$}}
\rm\vspace{4mm}

\noindent
\Prop{}{}\label{P11} Soit $\pi:X\rightarrow S$ un morphisme propre
d'espaces complexes dont les fibres sont de dimension (pure) $n$. Alors,  le foncteur $\pi^{\#}:{\mathcal G}\rightarrow {\mathcal G}\otimes \omega^{n}_{X/S}$ défini de ${\rm Coh}(S)$ sur  ${\rm Coh}(X)$ est covariant, exact à droite et muni d'un morphisme de foncteurs ${\rm I}\!{\rm R}^{n}\pi_{*}\pi^{\#}\rightarrow {\rm Id}$ induisant pour tout faisceau cohérent ${\mathcal G}$, un morphisme canonique
$${\rm I}\!{\rm R}^{n}\pi_{*}\pi^{\#}{\mathcal G}\rightarrow {\mathcal G}$$
de formation compatible aux changements de bases plats.\rm
\begin{proof}
Il est bien connu que le produit tensoriel est covariant et exact à droite en ces arguments. Il nous reste seulement à prouver l'existence du morphisme de foncteurs ${\rm I}\!{\rm R}^{n}\pi_{*}\pi^{\#}\rightarrow {\rm Id}$. Pour ce faire, nous allons montrer que le morphisme canonique $${\mathcal G}\otimes{\rm I}\!{\rm R}^{n}\pi_{*}\omega^{n}_{X/S}\rightarrow{\rm I}\!{\rm R}^{n}\pi_{*}(\pi^{*}({\mathcal G})\otimes \omega^{n}_{X/S})$$
est un isomorphisme  naturellement fonctoriel en ${\mathcal G}$  \vspace{2mm}
 
 \noindent
 $\bullet$ Cette flèche naturelle  est  la composée du morphisme naturel ${\mathcal
G}\rightarrow \pi_{*}\pi^{*}{\mathcal G}$ et du morphisme
``cup-produit''
$${\mathcal G}\otimes{\rm I}\!{\rm R}^{n}\pi_{*}\omega^{n}_{X/S}\rightarrow\pi_{*}\pi^{*}{\mathcal G}\otimes{\rm I}\!{\rm R}^{n}\pi_{*}\omega^{n}_{X/S}\rightarrow
{\rm I}\!{\rm R}^{n}\pi_{*}(\pi^{*}({\mathcal G})\otimes
\omega^{n}_{X/S})$$ 
Pour montrer que c'est un isomorphisme, il nous suffit de voir que les foncteurs
${\mathcal G}\rightarrow{\mathcal G}\otimes{\rm I}\!{\rm
R}^{n}\pi_{*}\omega^{n}_{X/S}$ et ${\mathcal G}\rightarrow{\rm I}\!{\rm
R}^{n}\pi_{*}(\pi^{*}({\mathcal G})\otimes \omega^{n}_{X/S})$ sont isomorphes.\vspace{1mm}

\noindent
Or, ils sont exactes à droites en raison de la propriété du produit tensoriel et de l'annulation des images
directes supérieures ${\rm I}\!{\rm R}^{k}\pi_{*}{\mathcal F}$ pour tout $k>n$ et tout faisceaux cohérents sur
$X$ (cf \lemmaref{L4}).  De plus, ils coincident naturellement sur les faisceaux localement libres.  Comme  le problème est de nature locale sur
$S$, un raisonnement par récurrence descendante sur la longueur de la
résolution localement libre ou la profondeur de ${\mathcal G}$ nous ramène essentiellement à considérer le cas d'un début de résolution à deux termes
$${\mathcal L}^{1}\rightarrow{\mathcal L}^{0} \rightarrow {\mathcal G}\rightarrow 0$$
dont on déduit (cf annexe) le diagramme commutatif
$$\xymatrix{{\mathcal L}^{1}\otimes{\rm I}\!{\rm R}^{n}\pi_{*}\omega^{n}_{X/S}\ar[r]\ar[d]& {\mathcal L}^{0}\otimes{\rm I}\!{\rm R}^{n}\pi_{*}\omega^{n}_{X/S}\ar[r]\ar[d]&{\mathcal G}\otimes{\rm I}\!{\rm R}^{n}\pi_{*}\omega^{n}_{X/S}\ar[r]\ar[d]&0\\
{\rm I}\!{\rm R}^{n}\pi_{*}(\pi^{*}({\mathcal L}^{1})\otimes
\omega^{n}_{X/S})\ar[r]& {\rm I}\!{\rm R}^{n}\pi_{*}(\pi^{*}({\mathcal L}^{0})\otimes \omega^{n}_{X/S})\ar[r]_{\alpha}& {\rm I}\!{\rm
R}^{n}\pi_{*}(\pi^{*}({\mathcal G})\otimes \omega^{n}_{X/S})\ar[r]&0}$$
avec un morphisme surjectif $${\rm I}\!{\rm R}^{n}\pi_{*}(\pi^{*}({\mathcal L}^{1})\otimes
\omega^{n}_{X/S})\rightarrow {\rm Ker}\alpha$$
Il est, alors, facile de conclure en vertu du cas localement libre.
\vspace{1mm}

\noindent Si $\omega^{n}_{X/S}$ est plat sur $S$, cette flèche se déduit du \corollaryref{C2} et du \corollaryref{C3}.
En effet, on a  la composée des morphismes
$${\rm I}\!{\rm
R}\pi_{*}({\mathcal H}^{0}({\rm I}\!{\rm L}\pi^{\#}))[n]\rightarrow{\rm I}\!{\rm L}\pi_{\#}({\mathcal H}^{0}({\rm I}\!{\rm L}\pi^{\#}))\rightarrow{\mathcal H}^{0}({\rm I}\!{\rm L}\pi_{\#})({\mathcal H}^{0}({\rm I}\!{\rm L}\pi^{\#}))=\pi_{\#}\circ \pi^{\#}\rightarrow{\rm Id}$$
et, par suite, le morphisme
$${\rm I}\!{\rm
R}^{n}\pi_{*}\circ \pi^{\#}\rightarrow{\rm Id}$$
La compatibilité aux changements de base arbitraires est une conséquence naturelle de cette propriété qu'ont ces foncteurs$\,\blacksquare$
\end{proof}

 \begin{rem} On peut remarquer que l'existence d'un morphisme canonique  de foncteurs ${\rm I}\!{\rm R}^{n}\pi_{*}\pi^{\#}\rightarrow {\rm Id}$ ne suffit pas pour dire que la paire $({\rm I}\!{\rm R}^{n}\pi_{*}$, $\pi^{\#})$ est adjointe comme il est facile de s'en convaincre en utilisant les propriétés classiques de l'adjonction.
 \end{rem}
\Prop{}{}\label{prop12} Soit
$\pi:X\rightarrow S$ un morphisme propre universellement $n$-équidimensionnel  d'espaces
complexes réduits de dimension  finie à  fibres de même dimension pure $n$ . Alors, pour tout faisceau
cohérent ${\mathcal G}$ sans torsion sur $S$, les propriétés
suivantes sont équivalentes:\vspace{2mm}

\noindent {\bf(i)} le
faisceau $\omega^{n}_{X/S}$ est plat sur $S$,\vspace{1mm}

\noindent {\bf(ii)}
le foncteur $\pi^{\#}:{\mathcal G}\rightarrow \pi^{*}{\mathcal G}\otimes
\omega^{n}_{X/S}$ est exact,\vspace{1mm}

\noindent {\bf(iii)} le faisceau
$\pi^{*}{\mathcal G}\otimes\omega^{n}_{X/S}$ est sans
torsion sur $X$.\vspace{1mm}\noindent\rm \rm \begin{proof}\vspace{1mm}

\noindent Toutes les assertions sont de nature locales et l'hypothèse de propreté du morphisme n'est là que pour garantir l'existence du faisceau $\omega^{n}_{X/S}$. De plus, étant de support $X$, sa $S$-platitude entraîne automatiquement la platitude de $\pi$.\vspace{1mm}

\noindent
{\bf(i)}$\Longrightarrow${\bf(ii)}: Evident puisque, par platitude de $\pi$,  toute suite exacte courte sur $S$ est transformée par $\pi^{*}$ en une suite exacte courte qui reste exacte par tensorisation  par  $\omega^{n}_{X/S}$ en raison de sa $S$-platitude. Remarquons qu'il nous suffit  simplement d'avoir l'exactitude à gauche puisque le foncteur $\pi^{\#}$ est déjà  exacte à droite.\vspace{2mm}

\noindent {\bf(ii)}$\Longrightarrow${\bf(iii)}: En fait, la condition {\bf(ii)} entraine le fait que le faisceau $\pi^{*}{\mathcal G}\otimes
\omega^{n}_{X/S}$ est sans ${\mathcal O}_{X}$-torsion. En effet, ${\mathcal G}$ étant sans torsion s'injecte dans un faisceau localement libre de rang fini  ${\mathcal L}$. Comme la platitude de $\pi$ préserve cette injection par image réciproque,  la tensorisation par $\omega^{n}_{X/S}$ exhibe, alors, grâce à  {\bf(ii)},  $\pi^{*}({\mathcal G})\otimes_{{\mathcal O}_{X}}
\omega^{n}_{X/S}$ comme un sous faisceau d'un faisceau sans torsion sur $X$ puisque $\omega^{n}_{X/S}$ l'est d'après la \corollaryref{c3}. D'où {\bf(iii)}.\vspace{1mm}

\noindent
 {\bf(iii)}$\Longrightarrow${\bf(i)}:  La nature locale du problème permet de travailler dans une une paramétrisation locale $\xymatrix{X\ar@/_/[rr]_{\pi}\ar[r]^{f}&S\times
U\ar[r]^{q}&S}$, 
pour laquelle on a
$$f_{*}(\pi^{*}{\mathcal G}\otimes_{{\mathcal O}_{X}}
\omega^{n}_{X/S})\simeq q^{*}{\mathcal G}\otimes_{{\mathcal O}_{S\times U}}f_{*}\omega^{n}_{X/S}\simeq {\mathcal H}om(f_{*}{\mathcal O}_{X}, \Omega^{n}_{S\times U/S})\otimes q^{*}{\mathcal G}$$
Notons ${\mathcal A}$ le faisceau ${\mathcal H}om(f_{*}{\mathcal O}_{X}, \Omega^{n}_{S\times U/S})$. Alors, ${\mathcal G}$ étant sans torsion, peut être installé (localement) dans une suite exacte courte
$$0\rightarrow{\mathcal G}\rightarrow{\mathcal L}\rightarrow{\mathcal K}\rightarrow 0$$
avec ${\mathcal L}$ localement libre de rang fini. On en déduit la suite exacte 
$$\xymatrix{0\rightarrow{\rm Tor}^{{\mathcal O}_{S\times U}}_{1}({\mathcal A}, q^{*}({\mathcal K}))\ar[r]&{\mathcal A}\otimes q^{*}({\mathcal G})\ar[r]&{\mathcal A}\otimes q^{*}({\mathcal L})\ar[r]&{\mathcal A}\otimes q^{*}({\mathcal K})\ar[r]&0 }$$
Comme le morphisme $${\mathcal A}\otimes q^{*}({\mathcal G})\rightarrow{\mathcal A}\otimes q^{*}({\mathcal L})$$
est un isomorphisme générique et que ces faisceaux sont sans torsion, il est injectif et, par conséquent 
$${\rm Tor}^{{\mathcal O}_{S\times U}}_{1}({\mathcal A}, q^{*}({\mathcal K}))=0$$
Ceci étant vrai pour tout faisceau cohérent ${\mathcal K}$, cela impose à ${\mathcal A}$ d'être $S$-plat (cf l'argument d'algèbre commutative classique ci-dessous). Mais $f$ étant fini, la donnée d'un  morphisme injectif 
$$0\rightarrow f_{*}(\pi^{*}{\mathcal G}_{1}\otimes_{{\mathcal O}_{X}}
\omega^{n}_{X/S})\rightarrow f_{*}(\pi^{*}{\mathcal G}_{2}\otimes_{{\mathcal O}_{X}}
\omega^{n}_{X/S})$$
est équivalente à la donnée d'un morphisme injectif
$$0\rightarrow \pi^{*}{\mathcal G}_{1}\otimes_{{\mathcal O}_{X}}
\omega^{n}_{X/S}\rightarrow \pi^{*}{\mathcal G}_{2}\otimes_{{\mathcal O}_{X}}
\omega^{n}_{X/S}$$
D'où la $S$-platitude du faisceau $\omega^{n}_{X/S}$.\vspace{1mm}

\noindent

 Le problème  d'algèbre commutative  classique auquel on faisait allusion peut s'énoncer sous la forme:\vspace{2mm}
 
\noindent
{\emph Soit ${A}$ et ${B}$ deux anneaux localement noethériens avec ${A}$ réduit et ${ B}$ plat sur ${A}$. Soit ${M}$ un ${B}$-module de type fini. Alors, on a:}\vspace{3mm}

\noindent
$ {N}\otimes_{{A}} { M}$ est sans ${A}$- torsion  pour tout  ${A}$-module ${N}$ de type fini et sans torsion  $\Longleftrightarrow$ ${ M}$ est ${A}$-plat.\vspace{1mm}

\noindent
On se convainc aisément que la question revient à prouver,  pour l'idéal maximal  ${\mathfrak m}_{A}$ de $A$ de corps résiduel $k$, l'équivalence\vspace{2mm}

\noindent
\centerline{{\emph{${{\mathfrak m}_{A}} \otimes_{A} M$ est sans $A$-torsion $\Longleftrightarrow\,$ $M$ est $A$-plat}}}\vspace{2mm}

\noindent
Si $M\,$ est $A$-plat, il est clair que ${{\mathfrak m}_{A}}\otimes_{A} M \,$ est sans torsion sur $A$  puisque ${\mathfrak m}_{A}$ l'est naturellement  (cf \cite{EGA4} \S. 6). \vspace{1mm}

\noindent
L'autre implication découle de la suite exacte courte
$$0\rightarrow {\mathfrak m}_{A}\rightarrow A\rightarrow k\rightarrow 0$$
et de la suite d'homologie qui s'en déduit
$${\rm Tor}^{A}_{1}(k, M)\rightarrow {{\mathfrak m}_{A}}
\otimes_{A} M\rightarrow M\rightarrow k\otimes_{A} M\rightarrow 0$$
Alors, $A$ étant réduit, l'idéal maximal ${\mathfrak m}_{A}$ n'est
pas un idéal associé de $A$ ( dans le cas contraire, il existerait $x$ non nul dans $A$ tel que ${\mathfrak m}_{A}.x=0$ et donc, en particulier, $x^{2}=0$),  ${\rm Tor}^{A}_{1}(k, M)$, étant
annulé par ${\mathfrak m}_{A}$,  s'identifie au sous faisceau de torsion
 ( sur $A$ ) de  ${{\mathfrak
m}_{A}}\otimes_{A} M$. Or, ce dernier étant sans $A$-torsion, par
hypothèse, il en résulte que ${\rm Tor}^{A}_{1}(k, M)=0$. Mais
d'après  la caractérisation des modules plats (cf \cite{Mat},
Thm p.17), cela revient à dire que $M$ est
$A$-plat$\,\blacksquare$\end{proof}\vspace{2mm}
 
 \noindent

\begin{rem}\vspace{1mm}

 \noindent {\bf(i)} Avec les notations précédentes, considérons à cet égard, un anneau $A$ de profondeur au moins deux d'idéal maximal ${\mathfrak m}_{A}$ de corps résiduel $k$. Si $M$ est la première syzygie de ${\mathfrak m}_{A}$ (donc de profondeur au moins deux!)  et $N$ un $A$-module sans torsion, alors, 
$M\otimes_{A} N$ est sans torsion si et seulement si $N$ est localement libre.\vspace{1mm}

\noindent
{\bf(ii)} L'hypothèse de
 $A$-platitude de $B$ peut-être éviter si l'on a ${\mathfrak m}_{A}.B\subset {\rm Rad}(B)$ puisque , pour tout anneau local noetherien
$(A,{\mathfrak m}_{A},k)$ réduit et toute  $A$-algèbre  $B$
 satisfaisant cette condition et tout  $B$-module $M$ de
 type fini et séparé pour la topologie ${\mathfrak
m}_{A}$-adic, on a\vspace{2mm}

\noindent
\centerline{$M$ est  $A$-plat si et seulement si
$M\otimes_{A}{\mathfrak m}_{A}$ est sans $A$-torsion.}
\vspace{2mm}
En effet, la suite exacte 
tensorisée avec $N$ donne la suite exacte
$$0\rightarrow {\rm Tor}^{A}_{1}({\mathfrak m}_{A}, N)\rightarrow M\otimes_{A} N\rightarrow F\otimes_{A} N$$
Ainsi,  $M\otimes_{A} N$ est sans torsion si et seulement si ${\rm Tor}^{A}_{1}({\mathfrak m}_{A}, N)={\rm Tor}^{A}_{2}(k, N)=0$, c'est-à-dire $N$ de profondeur homologique inférieure ou égale à 1.\vspace{1mm}

\noindent
{\bf(iii)} On  peut signaler un autre résultat assez classique (par exemple \cite{Mat},
Thm (50),p.153 seconde édition) disant que si  $(A,{\mathfrak m}_{A},k)$ et $(B,{\mathfrak m}_{B},k')$ sont deux anneaux locaux noethériens muni d'un morphisme local $\phi: A\rightarrow B$ et si $M$ et $N$ sont deux $A$-modules et $B$-modules finis avec $N$ plat sur $A$, alors
$${\rm Prof}_{B}(M\otimes_{A} N)={\rm Prof}_{A}(M) + {\rm Prof}_{B\times k}(N\otimes k)$$
\end{rem}\vspace{1mm}

\noindent 
\begin{obs} On peut signaler que le produit tensoriel de deux faisceaux  sans torsion est généralement de torsion. D'ailleurs, même si l'un deux est de profondeur supérieure à deux, cela ne garantit pas l'absence de torsion dans leur   produit tensoriel.\vspace{1mm}

\noindent
 Dans (\cite{A.L}, Remarque.(5.2.4),p.47) est donné un exemple simple où ${\rm E}xt_{R}^{m-n}(M, R)\otimes_{R} M$
 est de torsion  alors
  que l'un des deux modules est de profondeur au moins deux et l'autre sans torsion.\vspace{1mm}
  
  \noindent
En fait, on peut facilement voir que l'on ne peut faire guère mieux que \cite{EGA3}, $\S 6$:\vspace{1mm}

\noindent
 Soit $A$ un anneau noetherien de profondeur au moins deux. Soient
$\mathfrak {m}$ l'idéal maximal,  $k$  le corps résiduel et  $M$ la première (resp. seconde) syzygie de
$\mathfrak m$ ( resp. $k$). Alors ${\rm Prof}_{A}(M)=2$  et $M$ est sans
torsion puisque  sous module d'un module libre. Soit  $N$ un
$A$- module de type fini et sans $A$-torsion.\vspace{1mm}

\noindent   Alors $
M\otimes_{A}N$ est sans torsion si et seulement si ${\rm
pd}_{A}(N)\leq 1$.\vspace{2mm}

\noindent  En effet, en tensorisant par
$N$  la suite exacte courte  $0\rightarrow M\rightarrow F\rightarrow
{\mathfrak m}\rightarrow 0$ nous obtenons la suite exacte
$$0\rightarrow
{\rm Tor}^{A}_{1}({\mathfrak m}, N)\rightarrow M\otimes N\rightarrow
F\otimes N\rightarrow {\mathfrak m}\otimes N\rightarrow 0$$
Alors   ${\rm Tor}^{A}_{1}(\mathfrak{m}, N)={\rm
Tor}^{A}_{2}(k, N)$ est le sous module
de torsion de $M\otimes N$ puisqu'il est annulé par $\mathfrak m$. De plus, le quotient  $(M\otimes N)/{\rm Tor}^{A}_{1}(\mathfrak{m}, N)$  est
plongé dans le module sans torsion $F\otimes N$. Ainsi, le
produit tensoriel  $M\otimes N$ est sans torsion si et seulement si ${\rm
Tor}^{A}_{2}(k, N)=0$ et donc ${\rm pd}_{A}(N)\leq 1$.\vspace{1mm}

\noindent
Cela montre que dans la situation où  $f:A\rightarrow B$ est un morphisme d'anneaux localement noethériens avec $B$ plat sur $A$, $M$ (resp. $N$) un $B$ (resp. $A$)-module de type fini tel que ${\rm Prof}_{A}(M)\geq 2$ et sans $B$-torsion (resp. sans $A$-torsion), le produit tensoriel $M\otimes_{A} N$ présente de la torsion; ce qui nous incite à dire que l'on ne peut faire mieux que  \cite{EGA3}, \S 6 puisque cette condition sur la dimension
projective implique la liberté du module canonique.
\end{obs}

\vspace{2mm}

\noindent
\Prop{}{}\label{P13} Soit $\pi:X\rightarrow S$ un morphisme universellement $n$-équidimensionnel propre. On suppose $\omega^{n}_{X/S}$ plat sur $S$. Alors, le faisceau  $\pi^{*}{\mathcal G}\otimes\omega^{n}_{X/S}$ est sans
torsion fibre par fibre sur $X$ pour tout faisceau cohérent ${\mathcal G}$ sur $S$.
\begin{proof}
On veut montrer que pour tout fermé $\Sigma$ d'intérieur vide dans $X$ tel que, pour tout $s\in S$,  $\Sigma\cap \pi^{-1}(s)$ soit d'intérieur vide dans la fibre $\pi^{-1}(s)$, on a 
$${\mathcal H}^{0}_{\Sigma}(\pi^{*}{\mathcal G}\otimes\omega^{n}_{X/S})=0$$
 Alors, si ${\mathcal G}$ est sans torsion, ce n'est qu'un cas particulier de la \propositionref{prop12}. Dans le cas général, on considère la suite exacte courte 
$$0\rightarrow {\mathcal T}\rightarrow{\mathcal G}\rightarrow{\mathcal G}/{\mathcal T}\rightarrow0$$
où ${\mathcal T}$ est le sous faisceau de torsion de ${\mathcal G}$.\vspace{1mm}

\noindent
On en déduit, alors, la suite exacte courte 
$$0\rightarrow \pi^{*}({\mathcal T})\otimes_{{\mathcal O}_{X}}
\omega^{n}_{X/S}\rightarrow\pi^{*}({\mathcal G})\otimes_{{\mathcal O}_{X}}
\omega^{n}_{X/S}\rightarrow\pi^{*}({\mathcal G}/{\mathcal T})\otimes_{{\mathcal O}_{X}}
\omega^{n}_{X/S}\rightarrow 0$$
Pour $\Sigma$ d'intérieur vide fibre par fibre dans $X$, il s'en suit que
$${\mathcal H}^{0}_{\Sigma}(\pi^{*}({\mathcal T})\otimes_{{\mathcal O}_{X}}
\omega^{n}_{X/S})\simeq{\mathcal H}^{0}_{\Sigma}(\pi^{*}({\mathcal G})\otimes_{{\mathcal O}_{X}}\omega^{n}_{X/S})$$
Comme $\omega^{n}_{X/S}$ commute aux changements de base et que ${\rm Supp}({\mathcal T})$ est un  fermé strict $S_{0}$ de $S$, on a, en notant $\pi_{0}:X_{0}\rightarrow S_{0}$ le morphisme déduit de $\pi$ par changement de base, 
$${\mathcal H}^{0}_{\Sigma}(\pi^{*}({\mathcal T})\otimes_{{\mathcal O}_{X}}
\omega^{n}_{X/S})\simeq {\mathcal H}^{0}_{\Sigma_{0}}({\pi_{0}}^{*}({\mathcal T}_{0})\otimes_{{\mathcal O}_{X_{0}}}
\omega^{n}_{X_{0}/S_{0}})$$
avec  ${\mathcal T}_{0}:={\mathcal T}$ et $\Sigma_{0}:=\Sigma\cap X_{0}$ qui est encore d'intérieur vide fibre par fibre.  Si ${\mathcal T}_{0}$  est sans torsion sur $S_{0}$, le résultat est prouvé sinon on reprend le dévissage en remplaçant ${\mathcal G}$ par ${\mathcal T}$ jusqu'à obtenir un indice $k$ pour lequel  ${\mathcal T}_{k}$ soit sans torsion sur $S_k\,\blacksquare$
\end{proof}
\vspace{1mm}

\noindent

\section{\color{blue}{La preuve du \theoremref{thm3}: Relation entre les foncteurs $\pi^{!}_{\mathcal K}$ et $\pi^{\#}$.}}
 On commence par mettre en évidence, pour tout  élément, $\pi$, de ${\Bbb E}(S;n)$,  un morphisme canonique de foncteurs 
$$\Xi_{\pi}:{\pi}^{\#}\rightarrow
{\pi}^{!}_{{\mathcal K}}$$  ayant les propriétés suivantes:\vspace{3mm}

\noindent
{\bf(i)} {\emph{il est de formation compatible aux changements de base plats
et aux restrictions ouvertes sur $X$}}. \vspace{2mm}

\noindent
{\bf(ii)} il est  génériquement bijectif relativement à $\pi$. De plus, pour tout faisceau cohérent ${\mathcal
G}$ de support  $S$,  $\Xi_{\pi}({\mathcal G})$ est un isomorphisme
générique.\vspace{2mm}

\noindent
{\bf(iii)}  il est injectif si et seulement si ${\pi}^{\#}$ est exact.\vspace{2mm}

\noindent{\bf(iv)} il est surjectif si seulement si ${\pi}^{!}_{{\mathcal K}}$ est exact.
\vspace{2mm}

\noindent
{\bf(v)} ${\pi}^{\#}$ est exact si et seulement si ${\pi}^{!}_{{\mathcal K}}$ est exact.\vspace{2mm}

\noindent
{\bf(vi)} ${\pi}^{!}_{{\mathcal K}}$ est exact si et seulement si $\pi$ est de Cohen Macaulay.
\rm\vspace{2mm}

\noindent
{\bf(a) La construction de $\Xi_{\pi}$}: 

On peut la déduire  de  la dualité du  \theoremref{thm1} ou de la \propositionref{P6}\footnote{Elle émane aussi de la flèche $${\rm I}\!{\rm L}\pi^{*}({\mathcal G})\otimes^{{\rm I}\!{\rm L}} \pi^{!}{\mathcal O}_{S}[-n]\rightarrow \pi^{!}({\mathcal G})[-n]$$
dont on  prend la cohomologie de degré $0$ de sorte à avoir la composée
$${\mathcal H}^{0}({\rm I}\!{\rm L}\pi^{*}{\mathcal G})\otimes {\mathcal H}^{-n}(\pi^{!}{\mathcal O}_{S})\rightarrow {\mathcal H}^{-n}({\rm I}\!{\rm L}\pi^{*}{\mathcal G}\otimes \pi^{!}{\mathcal O}_{S})\rightarrow{\mathcal H}^{-n}(\pi^{!}({\mathcal G}))$$}. En effet, la donnée du morphisme $\int_{\pi,\mathcal G}:{\rm I}\!{\rm
R}^{n}\pi_{*}(\pi^{*}({\mathcal G})\otimes
\omega^{n}_{X/S})\rightarrow {\mathcal G}$ donné par la \propositionref{P11} est strictement
équivalente à la donnée d'un morphisme $\Xi_{\pi}:\pi^{*}{\mathcal
G}\otimes \omega^{n}_{X/S}\rightarrow {\mathcal H}^{-n}(\pi^{!}{\mathcal
G})$ rendant commutatif le diagramme
$$\xymatrix{{\rm I}\!{\rm R}^{n}\pi_{*}(\pi^{*}{\mathcal G}\otimes
\omega^{n}_{X/S})\ar[rr]\ar[rd]&&{\rm I}\!{\rm R}^{n}\pi_{*}{\mathcal
H}^{-n}(\pi^{!}{\mathcal G})\ar[ld]\\&{\mathcal G}&}$$
\vspace{1mm}

\noindent
Remarquons que la construction peut se faire localement en vertu de  la nature locale du problème ( sur $X$ (et $S$) en con sidérant une paramétrisation locale (arbitraire)  $\xymatrix{X\ar@/_/[rr]_{\pi}\ar[r]^{f}&S\times
U\ar[r]^{q}&S}$ 
Alors,

$$f_{*}(\pi^{*}({\mathcal G})\otimes_{{\mathcal O}_{X}}
\omega^{n}_{X/S})\simeq q^{*}({\mathcal G})\otimes_{{\mathcal O}_{S\times U}}f_{*}\omega^{n}_{X/S}\simeq {\mathcal H}om(f_{*}{\mathcal O}_{X}, \Omega^{n}_{S\times U/S})\otimes q^{*}({\mathcal G})$$
et 
$$f_{*}{\mathcal H}^{-n}(\pi^{!}({\mathcal G}))\simeq{\mathcal H}om(f_{*}{\mathcal O}_{X}, \Omega^{n}_{S\times U/S}\otimes q^{*}({\mathcal G}))$$ 
D'où, le morphisme naturel
$$f_{*}(\pi^{\#}({\mathcal G}))\rightarrow f_{*}(\pi^{!}_{{\mathcal K}}({\mathcal G}))$$
Comme la paramétrisation est arbitraire, on en déduit aisément le morphisme
$$\Xi_{\pi}({\mathcal G}):\pi^{\#}({\mathcal G})\rightarrow \pi^{!}_{{\mathcal K}}({\mathcal G})$$
pour tout faisceau cohérent ${\mathcal G}$.\vspace{2mm}

\noindent

{\bf(b) Propriétés de $\Xi_{\pi}$.}\vspace{2mm}

\noindent
 {\bf(i)} La compatibilité aux restrictions ouvertes sur $X$ et
aux changements de base plats est une conséquence du bon
comportement de ces foncteurs vis-à-vis de ces opérations comme il a été vu dans  les  \propositionref{P6} et \propositionref{P7}.\vspace{1mm}

\noindent
{\bf(ii)} La  bijectivité générique peut être induite par la nature du morphisme $\pi$ ou par celle du faisceau ${\mathcal G}$.\vspace{1mm}

\indent $\bullet$ Le morphisme $\Xi_{\pi}({\mathcal G})$ est un isomorphisme pour $\pi$ de Gorenstein (donc, en particulier, lisse) puisqu'en vertu du \corollaryref{C3}  les foncteurs ${\rm I}\!{\rm R}^{n}\pi_{*}$ et $\pi_{\#}$ sont isomorphes et, par suite, leurs adjoints  $\pi^{\#}\simeq\pi^{!}_{\mathcal K}$ aussi comme le précise le \corollaryref{C4}\footnote{si $\pi$ est lisse et propre, la bijectivité de  $\Xi_{\pi}$ découle immédiatement de la platitude de $\pi$ garantissant l'isomorphisme 
$${\rm I}\!{\rm L}\pi^{*}{\mathcal G}\otimes \pi^{!}({\mathcal O}_{S}[-n])\simeq\pi^{!}({\mathcal G})[-n]$$
et des annulations ${\mathcal H}^{j}(\pi^{!}({\mathcal O}_{S}))=0$ pour tout $j\not=-n$. On a, d'ailleurs, l'isomorphisme de Verdier 
$$\pi^{!}({\mathcal
G})\simeq\pi^{*}{\mathcal
G}\otimes \Omega^{n}_{X/S}[n]$$}. Comme un morphisme d'espaces complexes surjectif est toujours génériquement lisse ou de Gorenstein, $\Xi_{\pi}({\mathcal G})$ est génériquement bijectif. 
 \vspace{1mm}

\indent
$\bullet$ Comme tout faisceau cohérent arbitraire est génériquement localement libre sur un espace réduit et que les faisceaux ${\mathcal H}^{-n}(\pi^{!}{\mathcal G})$ et $\pi^{*}{\mathcal
G}\otimes_{{\mathcal O}_{X}}\omega^{n}_{X/S}$ sont canoniquement isomorphes sur les faisceaux localement libres, il s'en suit que  $\Xi_{\pi}({\mathcal G})$ génériquement bijectif en vertu de la compatibilité de la construction avec les rectrictions ouvertes. \vspace{1mm}

\noindent
Pour être précis, cette propriété de bijectivité générique  résultera de:\vspace{1mm}

$\bullet$ la compatibilité aux restrictions ouvertes sur $X$ et $S$ de ces 
foncteurs, de l'ouverture et de la surjectivité de $\pi$ et de ce
que l'ensemble des points en lesquels un morphisme  d'espaces
complexes est Gorenstein (ou lisse) est un
ouvert $U$ de $X$ dont le complémentaire $F$ est un fermé
analytique (cf [Ba1]). $X$ étant réduit ( ou à
partie régulière dense) , $U$ est dense dans $X$ ainsi que $\pi(U)$ puisque 
$\pi$ est ouvert et la partie régulière de $S$ dense dans
$S$ (le théorème de Remmert \cite{Rem1} nous dit que   $\pi(F)$ un sous ensemble analytique
strict de $S$).\vspace{1mm}

\noindent

\vspace{3mm}

\noindent
 {\bf(iii) $\pi^{\#}$  exact $\Longleftrightarrow$   $\Xi_{\pi}$ injectif:}
\vspace{3mm}

\noindent 
$\Longrightarrow$ : Nous supposerons dans toute la suite $S$ réduit.
Commençons par le cas où ${\mathcal G}$ est sans torsion. Alors, il s'injecte dans 
un faisceau cohérent localement libre ${\mathcal L}$ sur $S$. Mais
$\pi^{\#}$ étant supposé exact donc exact à gauche, on a une
injection $\pi^{\#}({\mathcal G})\rightarrow \pi^{\#}({\mathcal L})$ montrant que
 $\pi^{\#}({\mathcal G})$ est sans torsion puisque $\pi^{\#}({\mathcal L})$ l'est (cf \propositionref{prop12}. Comme $\Xi_{\pi}({\mathcal G})$ est un isomorphisme générique et $\pi_{{\mathcal K}}^{!}({\mathcal G})$ est sans
 torsion sur $X$ d'après la \propositionref{P6}), on en
 déduit que $\Xi_{\pi}({\mathcal G} )$ est injectif.\vspace{1mm}
 
 \noindent
 Supposons, à présent, que ${\mathcal G}$ admette de la torsion dont on notera le
 faisceau associé par ${\mathcal 
 T}$. Considérons, alors, la suite
 exacte courte
 $$0\rightarrow {\mathcal T}\rightarrow {\mathcal G}\rightarrow {\mathcal G}/{\mathcal T}\rightarrow 0$$
 dont on déduit, par fonctorialité, le diagramme commutatif à
 lignes exactes
$$\xymatrix{0\ar[r]&\pi^{\#}({\mathcal T})\ar[r]\ar[d]_{u}&\pi^{\#}({\mathcal G})\ar[d]_{v}\ar[r]
&\pi^{\#}({\mathcal G}/{\mathcal T})\ar[d]_{w}\ar[r]&0\\
0\ar[r]&\pi_{{\mathcal K}}^{!}({\mathcal T})\ar[r]&\pi_{{\mathcal K}}^{!}({\mathcal G})\ar[r]&
\pi_{{\mathcal K}}^{!}({\mathcal G}/{\mathcal T})&}$$ 
Comme $w$ est injectif d'après le cas traité précédemment, il est facile de voir que l'injectivité de $u$
est équivalente à celle de $v$. En désignant le support de ${\mathcal T}$ par $S'$ muni d'un plongement (local) $i$ dans $S$, on a un diagramme naturel de changement de base
$$\xymatrix{X'\ar[r]^{i'}\ar[d]^{\pi'}&X\ar[d]^{\pi}\\
S'\ar[r]_{i}&S}$$
avec $\pi'$ propre à fibres de dimension au plus $n$. L'assertion se prouve, alors, par dévissage ou récurrence sur la dimension de la base car $\pi^{\#}({\mathcal T})={\pi'}^{\#}({\mathcal T})$ et $\pi_{{\mathcal K}}^{!}({\mathcal T})={\pi'}_{{\mathcal K}}^{!}({\mathcal T})$. \vspace{2mm} 

\noindent

$\Longleftarrow$:\vspace{2mm}

\noindent Soit
$0\rightarrow {\mathcal G}_{1}\rightarrow {\mathcal G}_{1}$ une suite exacte de faisceaux
cohérents sur $S$. Alors, par fonctorialité et exactitude à
gauche du foncteur ${\mathcal G}\rightarrow {\mathcal H}^{-n}(\pi^{!}{\mathcal G})$, nous
obtenons  le diagramme commutatif à lignes exactes
$$\xymatrix{&0\ar[d]&0\ar[d]\\
&{\mathcal G}_{1}\otimes \omega^{n}_{X/S}\ar[d]\ar[r]&{\mathcal G}_{2}\otimes \omega^{n}_{X/S}\ar[d]\\
0\ar[r]&{\mathcal H}^{-n}(\pi^{!}{\mathcal G}_{1})\ar[r]&{\mathcal
H}^{-n}(\pi^{!}{\mathcal G}_{2})}$$
dont on déduit   sans peine l'exactitude à gauche du foncteur $\pi^{\#}$ et, par suite, son
exactitude puisqu'il est déjà exact à droite.
\vspace{3mm}

\noindent
{\bf(iv) $\Xi_{\pi}$ surjectif $\Longleftrightarrow$
$\pi_{{\mathcal K}}^{!}$ exact à droite.}\vspace{2mm}
 
 \noindent
$\Longrightarrow$:  Considérons la suite exacte de faisceaux cohérents
 ${\mathcal G}_{1}\rightarrow {\mathcal G}_{2}\rightarrow 0$. Alors, par fonctorialité,
 exactitude à droite du foncteur $\pi^{\#}$ et hypothèse, on a le diagramme
 commutatif $$\xymatrix{
{\mathcal G}_{1}\otimes \omega^{n}_{X/S}\ar[d]\ar[r]&{\mathcal G}_{2}\otimes
\omega^{n}_{X/S}
\ar[d]\ar[r]&0\\
{\mathcal H}^{-n}(\pi^{!}{\mathcal G}_{1})\ar[r]\ar[d]&{\mathcal
H}^{-n}(\pi^{!}{\mathcal G}_{2})\ar[d]\\
0&0&}$$
 La surjectivité de la flèche
 horizontale inférieure en résulte. D'où l'exactitude à droite du
 foncteur $\pi^{!}_{{\mathcal K}}$
  et par conséquent son exactitude puisqu'il est déjà exact à gauche.
  \vspace{1mm}

  \noindent On peut voir que cette surjectivité entraîne l'identification des deux foncteurs $\pi^{!}_{{\mathcal K}}$ et $\pi^{\#}$ puisqu'ils le sont sur les faisceaux localement libres. En effet, en prenant un début de résolution à deux termes 
  $${\mathcal L}^{1}\rightarrow{\mathcal L}^{0}\rightarrow {\mathcal G}\rightarrow 0$$
  de laquelle résulte le diagramme commutatif à lignes exactes et colonnes surjectives
 $$\xymatrix{\pi^{\#}({\mathcal L}^{1})\ar[r]\eq[d]_{\Xi_{\pi}({\mathcal L}^{1})}&\pi^{\#}({\mathcal L}^{0})\ar[r]\eq[d]_{\Xi_{\pi}({\mathcal L}^{0})}&\pi^{\#}({\mathcal G})\ar[r]\ar[d]_{\Xi_{\pi}({\mathcal G})}&0\\
\pi_{{\mathcal K}}^{!}({\mathcal L}^{1})\ar[r]\ar[d]&\pi_{{\mathcal K}}^{!}({\mathcal L}^{0})\ar[r]\ar[d]&\pi_{{\mathcal K}}^{!}({\mathcal G})\ar[d]\ar[r]&\cdots\\
0&0&0&}$$ 
En chassant dans ce diagramme commutatif, on établi facilement l'injectivité de $\Xi_{\pi}({\mathcal G})$ et, donc, sa bijectivité. D'où, en particulier, l'exactitude à droite du foncteur $\pi^{!}_{{\mathcal K}}$.\vspace{1mm}

\noindent $\Longleftarrow$:\vspace{1mm}
 
 \noindent
 
Comme on vient de le voir, l'exactitude à droite du foncteur $\pi^{!}_{{\mathcal K}}$ contraint ces deux foncteurs à être égaux puisqu'ils s'identifient sur les faisceaux localement libres et sont tous deux exacts à droite. On peut raisonner par récurrence sur la profondeur de
${\mathcal G}$ en initialisant  le processus  par le cas localement libre mais, ici, il est plus simple de considérer, comme ci-dessus,   un début de résolution localement libre 
$${\mathcal L}^{1}\rightarrow{\mathcal L}^{0}\rightarrow {\mathcal G}\rightarrow 0$$
à laquelle on applique les foncteurs $\pi_{{\mathcal K}}^{!}$ et $\pi^{\#}$ pour obtenir le  diagramme commutatif
$$\xymatrix{\pi^{\#}({\mathcal L}^{1})\ar[r]\ar[d]_{\Xi_{\pi}({\mathcal L}^{1})}&\pi^{\#}({\mathcal L}^{0})\ar[r]\ar[d]_{\Xi_{\pi}({\mathcal L}^{0})}&\pi^{\#}({\mathcal G})\ar[r]\ar[d]_{\Xi_{\pi}({\mathcal G})}&0\\
\pi_{{\mathcal K}}^{!}({\mathcal L}^{1})\ar[r]&\pi_{{\mathcal K}}^{!}({\mathcal L}^{0})\ar[r]&\pi_{{\mathcal K}}^{!}({\mathcal G})\ar[r]&0}$$
dont la bijectivité des deux premières flèches verticales entraine la surjectivité de la dernière et même sa bijectivité.\vspace{1mm}

\noindent
{\emph{En fait, ce raisonnement montre que si $\phi:T\rightarrow T'$ est un morphisme de foncteurs exact à droites qui est un isomorphisme sur les faisceaux localement libres alors c'est un isomorphisme}}.
\vspace{3mm}

\noindent
{\bf(v) $\pi^{\#}$ est  exact  si et seulement si $\pi_{{\mathcal K}}^{!}$ est exact.}\vspace{1mm}

\noindent 
$\bullet$ Si  $\pi_{{\mathcal K}}^{!}$ est exact à droite et, donc, exact, il est clair que $\pi^{\#}$ l'est aussi puisque ceux sont deux foncteurs exacts à droite coïncidant sur les faisceaux localement libres.\vspace{1mm}

\noindent
Réciproquement, si $\pi^{\#}$ est exact à gauche donc exact, le faisceau  $\omega^{n}_{X/S}$ est nécessairement  $S$-plat et  de support $X$ entrainant, en particulier  la platitude de $\pi$.\vspace{1mm}

\noindent 
Pour atteindre notre objectif, nous allons commencer par le \vspace{1mm}

\indent 
{{\bf(a) Cas ${\mathcal G}$ réfléxif.} }\vspace{1mm}

\noindent Dans ce cas, en considérant on déduit d'un début de résolution localement libre à deux termes, la suite exacte courte
$$0\rightarrow{\mathcal G}\rightarrow {\mathcal L}^{0}\rightarrow{\mathcal L}^{1}$$
conduisant, grâce à l'exactitude de $\pi^{\#}$ et l'exactitude à gauche de $\pi_{{\mathcal K}}^{!}$, au diagramme commutatif
$$\xymatrix{0\ar[r]&\pi^{\#}({\mathcal G})\ar[d]_{\Xi_{\pi}({\mathcal G})}\ar[r]& \pi^{\#}({\mathcal L}^{0})\ar[r]\ar[d]_{\Xi_{\pi}({\mathcal L}^{0})}&\pi^{\#}({\mathcal L}^{1})\ar[r]\ar[d]_{\Xi_{\pi}({\mathcal L}^{1})}\ar[r]&0\\
0\ar[r]&\pi_{{\mathcal K}}^{!}({\mathcal G})\ar[r]&\pi_{{\mathcal K}}^{!}({\mathcal L}^{0})\ar[r]&\pi_{{\mathcal K}}^{!}({\mathcal L}^{1})\ar[r]&\cdots}$$
montrant que $\Xi_{\pi}({\mathcal G})$ est un isomorphisme pour tout faisceau cohérent ${\mathcal G}$ réfléxif sur $S$.\vspace{1mm}

\noindent
{\bf(b) Cas ${\mathcal G}$ sans torsion.}\vspace{1mm}

\indent $\bullet$ On voit facilement que $\Xi_{\pi}({\mathcal G})$ est injectif pour tout faisceau cohérent ${\mathcal G}$ sans torsion puisqu'il  suffit de le voir comme sous faisceau d'un faisceau localement libre et d'utiliser l'exactitude à gauche de ces foncteurs et la bijectivité de $\Xi_{\pi}({\mathcal L})$ pour ${\mathcal L}$ localement libre. \vspace{1mm}

\indent $\bullet$ On en déduit, alors, qu'il est injectif pour tout faisceau cohérent ${\mathcal G}$. En effet, si ${\mathcal 
 T}$ est le sous faisceau de torsion, la suite
 exacte courte
 $$0\rightarrow {\mathcal T}\rightarrow {\mathcal G}\rightarrow {\mathcal G}/{\mathcal T}\rightarrow 0$$
 nous donne, par fonctorialité, le diagramme commutatif à
 lignes exactes
$$\xymatrix{0\ar[r]&\pi^{\#}({\mathcal T})\ar[r]\ar[d]_{u}&\pi^{\#}({\mathcal G})\ar[d]_{v}\ar[r]
&\pi^{\#}({\mathcal G}/{\mathcal T})\ar[d]_{w}\ar[r]&0\\
0\ar[r]&\pi_{{\mathcal K}}^{!}({\mathcal T})\ar[r]&\pi_{{\mathcal K}}^{!}({\mathcal G})\ar[r]&
\pi_{{\mathcal K}}^{!}({\mathcal G}/{\mathcal T})&}$$ 
Comme $w$ est injectif, il apparait que $v$ est injectif si et seulement si $u$ l'est. On procède, alors, par dévissage sachant que ces foncteurs commutent au changements de bases. En notant $S_{0}:={\rm Supp}({\mathcal T})$, ${\mathcal T}_{0}:={\mathcal T}$, $u_{0}:=u$, $\pi_{0}:X_{0}\rightarrow S_{0}$ le morphisme déduit de $\pi$ par changement de base, on a
$$\pi^{\#}({\mathcal T})={\pi_{0}}^{\#}({\mathcal T}_{0}),\,\,\, \pi_{{\mathcal K}}^{!}({\mathcal T})={\pi_{0}}_{{\mathcal K}}^{!}({\mathcal T}_{0})$$
Si ${\mathcal T}_{0}$ est sans torsion sur $S_0$, on a directement l'injectivité de $v$. Sinon, on dévisse ${\mathcal T}_0$ jusquà obtenir un indice $m$ tel que ${\mathcal T}_m$ soit sans torsion sur $S_m$ donnant l'injectivité de $u_m$ dont on déduira successivement l'injectivité de $u_{m+1}$, ..., $u_{0}$ et celle de $v$. \vspace{1mm}

\noindent
En résulte, alors, la bijectivité désirée. En effet, en considérant la suite exacte courte
$$0\rightarrow{\mathcal G}\rightarrow{\mathcal G}^{**}\rightarrow {\mathcal H}\rightarrow 0$$
on obtient le diagramme commutatif
$$\xymatrix{0\ar[r]&\pi^{\#}({\mathcal G})\ar[r]\ar[d]&\pi^{\#}({\mathcal G}^{**})\eq[d]\ar[r]
&\pi^{\#}({\mathcal H})\ar@{^{(}->}[d]\ar[r]&0\\
0\ar[r]&\pi_{{\mathcal K}}^{!}({\mathcal G})\ar[r]&\pi_{{\mathcal K}}^{!}({\mathcal G}^{**})\ar[r]&
\pi_{{\mathcal K}}^{!}({\mathcal H})&}$$ 
duquel ressort la bijectivité de la première flèche verticale qui est $\Xi_{\pi}({\mathcal G})$.\vspace{1mm}

\noindent
{\bf(c) Cas ${\mathcal G}$ cohérent quelconque.}\vspace{1mm}
Il découle immédiatement de ce qui précède car en reprenant la suite exacte courte de {\bf(b)} et le diagramme induit, il est facile de voir que $\Xi_{\pi}({\mathcal G})$ est surjective. D'où sa bijectivité$\blacksquare$\vspace{2mm}

\noindent
{\bf(vi)$\Xi_{\pi}$ est un isomorphisme si et seulement si  $\pi$ est un morphisme de Cohen Macaulay}.\rm
\vspace{3mm}
 
 \noindent
$\Longleftarrow$ : Si $\pi$ est de Cohen Macaulay, on a $\pi^{!}({\mathcal O}_{S}[-n]\simeq \omega^{n}_{X/S}$  et, par conséquent, l'isomorphisme
$$\pi^{*}({\mathcal G})\otimes \omega^{n}_{X/S}\simeq{\rm I}\!{\rm L}\pi^{*}({\mathcal G})\otimes^{{\rm I}\!{\rm L}} \pi^{!}{\mathcal O}_{S}[-n]\rightarrow \pi^{!}({\mathcal G})[-n]$$
qui s'écrit aussi
$$\pi^{*}({\mathcal G})\otimes \omega^{n}_{X/S}\simeq{\mathcal H}^{-n}(\pi^{!}{\mathcal G})$$
prouvant l'implication désirée.\vspace{1mm}

\noindent
$\Longrightarrow$ : L'exactitude à gauche de $\pi^{\#}$ entraîne la $S$-platitude de $\omega^{n}_{X/S}$ qui, étant de support $X$, entraîne la platitude de $\pi$. l'exactitude à droite de $\pi_{{\mathcal K}}^{!}$, quant à elle, implique que sa formation commute à tout changement de base et que, par conséquent, $\omega^{n}_{X/S}$ commute aux changements de base arbitraire. L'isomorphisme
$$\pi^{*}({\mathcal G})\otimes^{{\rm I}\!{\rm L}} \pi^{!}{\mathcal O}_{S}[-n]\rightarrow \pi^{!}({\mathcal G})[-n]$$
ne garantit pas  les annulations ${\mathcal H}^{j}(\pi^{!}{\mathcal O}_{S})=0$ pour tout $j\not= -n$; dans un plongement donné, cela nécessite d'annuler tous les faisceaux "Ext" en degré supérieur à ${\rm codim}(X) +1$.\vspace{2mm}

\noindent
On va adopter une autre stratégie consistant à raisonner localement sachant  que la propriété d'être de Cohen Macaulay pour un morphisme est  une propriété locale sur la source. On considère, alors, une paramétrisation locale en un certain point de $X$ donnée par  $\xymatrix{X\ar@/_/[rr]_{\pi}\ar[r]^{f}&Y\ar[r]^{q}&S}$. Alors, le morphisme $\Xi_{\pi}$ étant compatible aux localisations (puisque la formation de ces foncteurs commute à tous changements de base plats),  on a, relativement à cette paramétrisation locale, le morphisme  
$$f_{*}(\pi^{*}{\mathcal G}\otimes_{{\mathcal
O}_{X}}\omega^{n}_{X/S})\rightarrow f_{*}{\mathcal H}^{-n}(\pi^{!}{\mathcal G})$$
que l'on peut aussi écrire
$$q^{*}{\mathcal G}\otimes
{\mathcal H}om(f_{*}{\mathcal O}_{X}, \Omega^{n}_{Y/S})\rightarrow{\mathcal
H}om(f_{*}{\mathcal O}_{X}, q^{*}{\mathcal G} \otimes_{{\mathcal
O}_{Y}}\Omega^{n}_{Y/S})$$
On va montrer que cette flèche est  bijective pour tout faisceau cohérent ${\mathcal G}$ si et seulement si  le faisceau   $f_{*}{\mathcal O}_{X}$   est localement libre de rang fini. \vspace{1mm}

\noindent On terminera, alors, en utilisant le fait que cette condition est équivalente à la platitude de $f$ (choisi arbitrairement) qui, au vu du critère de platitude par fibre impose aux fibres d'être de Cohen macaulay.\vspace{2mm}

\noindent

$\bullet$ Si $f_{*}{\mathcal O}_{X}$   est localement libre de rang fini, le résultat est évident.\vspace{1mm}

$\bullet$ Réciproquement, $\bar\pi$ est exact à droite, alors il coincide avec le foncteur $\pi^{\#}$ et cette flèche est un isomorphisme. Alors, comme on l'a déjà vu, cela impose la $S$-platitude du faisceau $\omega^{n}_{X/S} $ et, donc, la platitude de $\pi$. D'ailleurs,  de l'isomorphisme précédent, il résulte  que, pour tout faisceau cohérent sans torsion ${\mathcal G}$,  le faisceau cohérent ${\mathcal
H}om(f_{*}{\mathcal O}_{X}, q^{*}{\mathcal G} \otimes_{{\mathcal
O}_{Y}}\Omega^{n}_{Y/S})$ étant  sans torsion, il en sera de même  pour le produit tensoriel $q^{*}{\mathcal G}\otimes
{\mathcal H}om(f_{*}{\mathcal O}_{X}, \Omega^{n}_{Y/S})$\footnote{On voit bien, d'après nos discussions précédentes que cela n'a aucune de chance de se réaliser sans conditions fortes sur $f_{*}{\mathcal O}_{X}$. Signalons, au passage, que même pour  ${\mathcal G}$ sans torsion, le produit tensoriel aura de la torsion alors que ${\mathcal
H}om(f_{*}{\mathcal O}_{X}, q^{*}{\mathcal G} \otimes_{{\mathcal
O}_{Y}}\Omega^{n}_{Y/S})$ en est dépourvu!}
Or ceci n'est réalisable que si et seulement si ${\mathcal H}om(f_{*}{\mathcal O}_{X},\Omega^{n}_{Y/S})=f_{*}\omega^{n}_{X/S} $ est localement libre. En effet, prenons ${\mathcal G}:={\mathcal H}om({\mathcal G'}, {\mathcal O}_{Y})$. Alors, on doit 
$$ {\mathcal H}om(q^{*}{\mathcal G'}, {\mathcal O}_{Y})\otimes f_{*}\omega^{n}_{X/S}\simeq {\mathcal H}om(f_{*}({\mathcal O}_{X})\otimes {\mathcal G'}, \Omega^{n}_{Y/S})\simeq {\mathcal H}om(f_{*}({\mathcal G'}, f_{*}\omega^{n}_{X/S})$$
et ce quelque soit le faisceau cohérent ${\mathcal G'}$. Mais cela impose à $f_{*}\omega^{n}_{X/S}$ d'être localement libre et, donc, à $f_{*}{\mathcal O}_{X}$ de l'être aussi$\,\blacksquare$\vspace{2mm}

\noindent 
\cor{}{}\label{cor8}{} Supposons le foncteur $\pi^{\#}$ exact à gauche. Alors,  pour tout faisceau cohérent ${\mathcal G}$ de profondeur au moins deux dans $S$,  le morphisme  $\Xi_{\pi}({\mathcal G})$ est un isomorphisme.\rm
\begin{proof}

Considérons, d'abord, le cas d'un faisceau réfléxif  ${\mathcal G}$ et un début de résolution localement libre du ${\mathcal O}_{S}$-dual ${\mathcal G}^{*}$
$${\mathcal L}^{1}\rightarrow{\mathcal L}^{0}\rightarrow {\mathcal G}^{*}\rightarrow 0$$ dont on déduit le début de résolution
$$ 0\rightarrow{\mathcal G}\rightarrow({\mathcal L}^{0})^{*}\rightarrow({\mathcal L}^{1})^{*}$$ et le diagramme commutatif

$$\xymatrix{0\ar[r]&\pi^{\#}({\mathcal G})\ar[r]\ar[d]_{\Xi_{\pi}({\mathcal G})}&\pi^{\#}(({\mathcal L}^{0})^{*})\ar[r]\ar[d]_{\Xi_{{\pi}(({\mathcal L}^{0})^{*})}}&\pi^{\#}(({\mathcal L}^{1})^{*}))\ar[d]_{\Xi_{{\pi}(({\mathcal L}^{1})^{*}))}}\\
0\ar[r]&\pi_{{\mathcal K}}^{!}({\mathcal G})\ar[r]&\pi_{{\mathcal K}}^{!}(({\mathcal L}^{0})^{*}))\ar[r]&\pi_{{\mathcal K}}^{!}(({\mathcal L}^{1})^{*}))}$$
La bijectivité de $\Xi_{\pi({\mathcal G})}$ découlent  de celle des deux autres flèches verticales. \vspace{1mm}

\noindent
Si ${\mathcal G}$ est un faisceau de profondeur au moins deux, on peut l'insérer dans la suite exacte courte
$$0\rightarrow{\mathcal G}\rightarrow{\mathcal G}^{**}\rightarrow {\mathcal H}\rightarrow 0$$
où ${\mathcal H}$ est faisceau sans torsion. Un diagramme commutatif analogue au précédent utilisant la bijectivité de $\Xi_{\pi^{\#}(({\mathcal G}^{**})}$ et l'injectivité de $\Xi_{\pi^{\#}(({\mathcal H})}$ permet de conclure$\,\blacksquare$
\end{proof}
\cor{}{}\label{cor9}{} Avec les notations précédentes, on a, pour un morphisme de type $S_2$, les équivalences:\vspace{1mm}

\noindent
{\bf(i)}  $\pi_{{\mathcal K}}^{!}$ est exact à droite, \vspace{1mm}

\noindent
{\bf{(ii)}} ${\mathcal H}^{-n+j}(\pi^{!}{\mathcal G})=0$ pour tout $j>0$.
\rm
\begin{proof}
$\Longrightarrow$: Comme l'exactitude à droite du foncteur  $\pi_{{\mathcal K}}^{!}$ équivaut au fait que le morphisme $\pi$ est de Cohen Macaulay en vertu du  \theoremref{thm3} et que, relativement à toute  paramétrisation locale sur  $X$ donnée par  $\xymatrix{X\ar@/_/[rr]_{\pi}\ar[r]^{f}&Y\ar[r]^{q}&S}$, on a l'isomorphisme
$$f_{*}{\mathcal H}^{-n+j}(\pi^{!}{\mathcal G})\simeq{\mathcal Ext}^{j}(f_{*}{\mathcal O}_{X}, q^{*}{\mathcal G} \otimes_{{\mathcal
O}_{Y}}\Omega^{n}_{Y/S})$$
Or, de la locale liberté de $f_{*}{\mathcal O}_{X}$ résulte l'annulation des faisceaux ${\mathcal Ext}^{j}(f_{*}{\mathcal O}_{X}, q^{*}{\mathcal G} \otimes_{{\mathcal
O}_{Y}}\Omega^{n}_{Y/S})$ et ce, pour tout $j\geq 1$. D'où, 
${\mathcal H}^{-n+j}(\pi^{!}{\mathcal G})=0$.\vspace{1mm}

\noindent
$\Longleftarrow$: 
Le foncteur  ${\mathcal G}\rightarrow{\mathcal H}^{-n+1}(\pi^{!}{\mathcal G})$ est bien évidemment exacte à gauche et, par conséquent, le foncteur $\pi_{{\mathcal K}}^{!}$  est exact à droite.\vspace{1mm}

\noindent
Pour tout diagramme de $S$-plongement local,  $\xymatrix{X\ar@/_/[rr]_{\pi}\ar[r]^{\sigma}&S\times Z\ar[r]^{q}&S}$ où $Z$
 est une variété de Stein de dimension $n+p$, on a 
 $${\mathcal Ext}^{p+j}(\sigma_{*}{\mathcal O}_{X}, \Omega^{n+p}_{S\times Z/S})=0,\,\,\forall\,j\geq 1$$                                                          Comme $$\sigma_{*}{\mathcal H}^{-n+1}(\pi^{!}{\mathcal G})\simeq {\mathcal Ext}^{p+1}(\sigma_{*}{\mathcal O}_{X}, {\mathcal G} \otimes_{{\mathcal
O}_{S\times Z}}\Omega^{n+p}_{S\times Z/S})$$
Un raisonnement par récurrence sur la profondeur de ${\mathcal G}$ et en procédant par dévissage, on obtient l'annulation ${\mathcal Ext}^{p+1}(\sigma_{*}{\mathcal O}_{X}, {\mathcal G} \otimes_{{\mathcal
O}_{S\times Z}}\Omega^{n+p}_{S\times Z/S})=0$
pour tout faisceau cohérent ${\mathcal G}$                                                                                           
\noindent Signalons, au passage que les lemmes {\bf(5.2.2)} et {\bf(5.2.5)} de \cite{A.L} montrent que {\bf(ii)} impose au morphisme $\pi$ d'être de Cohen Macaulay$\,\blacksquare$

\end{proof}
\cor{}{}\label{cor10}{} Si $\pi:X\rightarrow S$ est un ${\rm S}_{2}$-morphisme plat d'espaces complexes réduits, le foncteur $\pi_{{\mathcal K}}^{!}$ est exact à droite sur les faisceaux cohérents de profondeur au moins deux sur $S$.\rm\vspace{3mm}

\noindent
\begin{proof} Il suffit d'appliquer le \corollaryref{c4} qui est une adaptation, à notre cadre, des lemmes {\bf(5.2.2)} et {\bf(5.2.5)} de \cite{A.L}$\blacksquare$ \end{proof} 
\vspace{1mm}

\noindent
\section{\color{blue}{Les foncteurs $\pi_{{\mathcal K}}^{!}$ et $\pi^{\#}$ sur des exemples.}}\vspace{2mm}

\noindent Nous proposons dans ce qui suit
un exemple simple dans  lequel $\Xi_{\pi}$ n'est pas un
isomorphisme.\vspace{1mm}

\noindent {\bf(i) Exemple avec une base non normale.} Soient $X={\Bbb C}^{2}$,
$S=\{(x,y,z)\in {\Bbb C}^{3}/ xy^{2}=z^{2}\}$ et $\pi:X\rightarrow S$ le morphisme de paramétrisation  défini par $(u,v)\rightarrow (u^{2}, v^{2},
uv)$ qui est aussi une désingularisation de $S$. \vspace{1mm}

\noindent Alors, l'égalité des foncteurs $\pi^{!}_{{\mathcal K}}$ et $\pi^{\#}$ s'exprime par l'isomorphisme
$${\mathcal G}\otimes\pi_{*}{\omega}^{0}_{X/S}\simeq {\mathcal H}om(\pi_{*}{\mathcal O}_{X}, {\mathcal G})$$ 
pour tout faisceau cohérent ${\mathcal G}$ sur $S$.\vspace{1mm}

\noindent Comme $$\pi_{*}{\omega}^{0}_{X/S}\simeq{\mathcal H}om(\pi_{*}{\mathcal O}_{X}, {\mathcal O}_{S}):={\mathcal H}om({\mathcal L}^{0}_{S}, {\mathcal O}_{S})$$
et que ${\mathcal L}^{0}_{S}$ est réfléxif puisque
$${\mathcal L}^{0}_{S}:=\pi_{*}({\mathcal O}_{X})\simeq {\mathcal H}om({\mathcal L}^{2}_{S}, \omega^{2}_{S})$$
D'ailleurs, en utilisant le fait que $X$ et $S$ sont de Gorenstein, on voit facilement  que 
$${\mathcal L}^{0}_{S}\simeq {\mathcal H}om({\mathcal H}om({\mathcal L}^{0}_{S}, {\mathcal O}_{S}), {\mathcal O}_{S})$$
Ainsi, l'isomorphisme déduit de l'égalité de nos foncteurs s'écrit aussi
$${\mathcal G}\otimes{\mathcal L}^{0}_{S} \simeq{\mathcal H}om({\mathcal L}^{0}_{S}, {\mathcal G})$$ 
ce qui est faux comme le montre l'exemple simple en prenant ${\mathcal G}:=\omega^{1}_{S}$. Si tel était le cas, on aurait l'injection
$$\omega^{1}_{S}\otimes{\mathcal L}^{0}_{S}\subset \omega^{1}_{S}$$
déduite de l'injection de ${\mathcal O}_{S}$ dans ${\mathcal L}^{0}_{S}$ et de l'absence de torsion dans $\omega^{1}_{S}$. mais les  sections $\frac{dz}{y}$ et $\frac{z}{y}$ de $\omega^{1}_{S}$ et ${\mathcal L}^{0}_{S}$ respectivement donnent par cup-produit la forme $\frac{zdz}{y^2}$ qui n'est pas du tout une section de $\omega^{1}_{S}$.\vspace{2mm}

\noindent
{\bf(ii) Exemple avec une base normale.}
Soient $X={\Bbb C}^{2}$,
$S=\{(x,y,z)\in {\Bbb C}^{3}/ xy=z^{2}\}$ ( qui présente une
singularité quotient)  et $\pi:X\rightarrow S$ le morphisme de
paramétrisation  défini par $(u,v)\rightarrow (u^{2}, v^{2},
uv)$.\vspace{1mm}

\noindent Alors, $\pi$ n'est ni plat, ni même de {\emph
Tor-dimension} finie (la multiplicité pontuelle en un point autre que l'origine est $2$  alors qu'en l'origine elle vaut $3$) mais, toutefois,  analytiquement géométriquement plat
puisque équidimensionnel sur une base normale (c-à-d il définit une famille analytique de points). 
Là encore, l'égalité des foncteurs $\pi^{!}_{{\mathcal K}}$ et $\pi^{\#}$ revient exactement à imposer l'isomorphisme
$$\pi^{*}({\mathcal G})\otimes\pi^{!}({\mathcal O}_{S})\simeq\pi^{!}({\mathcal G})  $$
qui s'exprime aussi sous la forme
$${\mathcal G}\otimes\pi_{*}{\omega}^{0}_{X/S}\simeq {\mathcal H}om(\pi_{*}{\mathcal O}_{X}, {\mathcal G})$$ 
Comme $X$ et $S$ sont des espaces de Gorenstein\footnote{ On est dans la situation
où le faisceau des fonctions méromorphes régulières de
Kunz-Waldi-Kersken (qui s'identifieà l'idéal conducteur)
coincide avec le faisceau canonique relatif.} , on a
$$\pi_{*}{\omega}^{0}_{X/S}\simeq{\mathcal H}om(\pi_{*}{\mathcal O}_{X}, {\mathcal O}_{S})\simeq \pi_{*}{\mathcal O}_{X}$$
D'où $${\mathcal G}\otimes\pi_{*}{\omega}^{0}_{X/S}\simeq \pi_{*}\pi^{*}({\mathcal G})$$
Ainsi, notre isomorphisme de foncteurs, donne, pour tout faisceau cohérent ${\mathcal G}$
$$\pi_{*}\pi^{*}({\mathcal G})\simeq{\mathcal H}om(\pi_{*}{\mathcal O}_{X}, {\mathcal G})$$
 Or, pour  ${\mathcal
G}=\omega^{1}_{S}$ qui est réfléxif, on obtient
$$\xymatrix{\pi_{*}\pi^{*}(\omega^{1}_{S})\eq[r]\eq[d]&{\mathcal H}om(\pi_{*}{\mathcal O}_{X}, \omega^{1}_{S})\eq[r]&{\mathcal H}om(\pi_{*}\pi^{*}(\omega^{1}_{S}), \omega^{2}_{S})\eq[d]\\
\pi_{*}{\mathcal H}om(\pi^{*}(\omega^{1}_{S}),{\mathcal O}_{X})\eq[rr]&&\pi_{*}{\mathcal H}om(\pi^{*}(\omega^{1}_{S}), \omega^{2}_{X})} $$
et, donc, la réfléxivité de $\pi_{*}\pi^{*}(\omega^{1}_{S})$ et celle de $\pi^{*}(\omega^{1}_{S})$. Mais pour l'un comme pour l'autre, cela les contraint d'être localement libres puisque $X$ et $S$ sont des surfaces normales. Par conséquent, $\pi^{*}(\omega^{1}_{S})$ s'identifie nécessairement à $\Omega^{1}_{X}$ sachant que l'on a une image réciproque injective $\pi^{*}(\omega^{1}_{S})\rightarrow\Omega^{1}_{X}$ entre faisceaux réfléxifs (rappelons que $\omega^{1}_{S}={\mathcal L}^{1}_{S}$ faisceau des formes méromorphes se prolongeant holomorphiquement dans toute désingularisée de $S$ mais c'est aussi les formes holomorphes sur ${\Bbb C}^{2}$ invariantes sous l'action du groupe symétrique d'ordre $2$). Mais ceci est clairement impossible car, par exemple, la forme holomorphe $du$ n'est pas image réciproque d'une section de $\omega^{1}_{S}$ !
 \vspace{2mm}
 
\noindent
\begin{rem} Pour un morphisme de Cohen Macaulay $\pi:X\rightarrow S$, on a un isomorphisme canonique
$${\mathcal O}_{X}\simeq{\mathcal H}om(\omega^{n}_{X/S}, \omega^{n}_{X/S})$$
  En effet, on a 
 $${\rm I}\!{\rm R}{\mathcal H}om(\pi^{!}({\mathcal O}_{S}), \pi^{!}({\mathcal O}_{S}))\simeq {\rm I}\!{\rm R}{\mathcal H}om(\pi^{!}({\mathcal O}_{S}), {\rm I}\!{\rm R}{\mathcal H}om(\pi^{*}{\mathcal D}^{\bullet}_{S}, {\mathcal D}^{\bullet}_{X}))\simeq {\rm I}\!{\rm R}{\mathcal H}om(\pi^{!}({\mathcal O}_{S})\otimes^{{\rm I}\!{\rm L}}\pi^{*}{\mathcal D}^{\bullet}_{S}, {\mathcal D}^{\bullet}_{X}))$$
 Alors, $\pi$ étant de Cohen Macaulay, $$\pi^{!}({\mathcal O}_{S})[-n]=\omega^{n}_{X/S},\,\,\pi^{!}({\mathcal O}_{S})\otimes^{{\rm I}\!{\rm L}}\pi^{*}{\mathcal D}^{\bullet}_{S}\simeq \pi^{!}({\mathcal D}^{\bullet}_{S})={\mathcal D}^{\bullet}_{X}$$
D'où
$${\rm I}\!{\rm R}{\mathcal H}om(\omega^{n}_{X/S}, \omega^{n}_{X/S})\simeq {\rm I}\!{\rm R}{\mathcal H}om({\mathcal D}^{\bullet}_{X}, {\mathcal D}^{\bullet}_{X})={\mathcal O}_{X}$$
D'où, en particulier, 
$${\mathcal O}_{X}\simeq{\mathcal H}om(\omega^{n}_{X/S}, \omega^{n}_{X/S})$$
On en déduit que dans factorisation locale (comme dans la preuve ci-dessus) le faisceai $f_{*}{\mathcal O}_{X}$ est toujours réfléxif.\end{rem}
 \section{\color{blue}{Preuve du \theoremref{thm4}.}}
On va montrer que pour tout morphisme $\pi: X\rightarrow S$ propre et plat avec $\omega^{n}_{X/S}$ $S$-plat, le morphisme canonique de foncteurs
$${\rm I}\!{\rm R}\pi_{*}[n]\rightarrow{\rm I}\!{\rm L}\pi_{\#}$$
est un isomorphisme si et seulement si $\pi$ est de Cohen Macaulay.

\vspace{1mm}

\noindent
\par

$\bullet$ Supposons, $\pi$ de Cohen Macaulay. Alors, ${\mathcal H}^{j}(\pi^{!}({\mathcal O}_{S}))=0$ pour tout $j\not=-n$ et $\omega^{n}_{X/S}:={\mathcal H}^{-n}(\pi^{!}({\mathcal O}_{S}))$ est plat sur $S$. De plus, le morphisme canonique 
$$\pi^{!}{\mathcal A}^{\bullet}\otimes^{{\rm I}\!{\rm L}}{\rm I}\!{\rm L}\pi^{*}{\mathcal B}^{\bullet}\rightarrow\pi^{!}({\mathcal A}^{\bullet}\otimes^{{\rm I}\!{\rm L}}{\mathcal B}^{\bullet} )$$
est un isomorphisme. On a, donc, en particulier, 
$$\omega^{n}_{X/S} \otimes^{{\rm I}\!{\rm L}}\pi^{*}{\mathcal B}^{\bullet}\simeq\pi^{!}({\mathcal O}_{S})[-n]\otimes^{{\rm I}\!{\rm L}}\pi^{*}{\mathcal B}^{\bullet}\simeq \pi^{!}({\mathcal B}^{\bullet})[-n]$$
On a, alors, le diagramme commutatif
$$\xymatrix{{\rm I}\!{\rm R}\pi_{*}{\rm I}\!{\rm R}{\mathcal H}om({\mathcal F}^{\bullet}, \omega^{n}_{X/S}\otimes\pi^{*}{\mathcal N}^{\bullet})\eq[d]\eq[r]&{\rm I}\!{\rm R}{\mathcal H}om({\rm I}\!{\rm L}\pi_{\#}({\mathcal  F}^{\bullet}), {\mathcal N}^{\bullet})\ar[d]\\
{\rm I}\!{\rm R}\pi_{*}{\rm I}\!{\rm R}{\mathcal H}om({\mathcal F}^{\bullet}, \pi^{!}({\mathcal N}^{\bullet})[-n])\eq[r]&{\rm I}\!{\rm R}{\mathcal H}om({\rm I}\!{\rm R}\pi_{*}{\mathcal F}^{\bullet}[n], {\mathcal N}^{\bullet}))}$$
imposant à la seconde flèche verticale d'être un isomorphisme. En prenant ${\mathcal N}^{\bullet}:={\mathcal D}^{\bullet}_{S}$ et en dualisant cette flèche, on obtient l'égalité des foncteurs 
$${\rm I}\!{\rm R}\pi_{*}[n]={\rm I}\!{\rm L}\pi_{\#}$$
En effet, il suffit de remarquer que
$$\omega^{n}_{X/S}\otimes^{{\rm I}\!{\rm L}} {\rm I}\!{\rm L}\pi^{*}{\mathcal D}^{\bullet}_{S}\simeq\pi^{!}({\mathcal O}_{S})\otimes^{{\rm I}\!{\rm L}} {\rm I}\!{\rm L}\pi^{*}{\mathcal D}^{\bullet}_{S}[-n]\simeq\pi^{!}({\mathcal D}^{\bullet}_{S})[-n]\simeq {\mathcal D}^{\bullet}_{X}[-n]$$
et d'utiliser la définition 
$${\rm I}\!{\rm L}\pi_{\#}({\mathcal F}^{\bullet})={\rm I}\!{\rm D}_{S}({\rm I}\!{\rm R}\pi_{*}{\rm I}\!{\rm R}{\mathcal H}om({\mathcal F}^{\bullet}, {\mathcal D}^{\bullet}_{X}[-n]))\simeq{\rm I}\!{\rm D}_{S}({\rm I}\!{\rm R}{\mathcal H}om({\rm I}\!{\rm R}\pi_{*}{\mathcal F}^{\bullet}[n], {\mathcal D}^{\bullet}_{S}))= {\rm I}\!{\rm R}\pi_{*}{\mathcal F}^{\bullet}[n]$$
\indent $\bullet$ Réciproquement, supposons donnée la précédente relation entre ces foncteurs. Le diagramme précédent nous donne cette fois
$$\xymatrix{{\rm I}\!{\rm R}\pi_{*}{\rm I}\!{\rm R}{\mathcal H}om({\mathcal F}^{\bullet}, \omega^{n}_{X/S}\otimes\pi^{*}{\mathcal N}^{\bullet})\ar[d]\eq[r]&{\rm I}\!{\rm R}{\mathcal H}om({\rm I}\!{\rm L}\pi_{\#}({\mathcal  F}^{\bullet}), {\mathcal N}^{\bullet})\eq[dd]\\
{\rm I}\!{\rm R}\pi_{*}{\rm I}\!{\rm R}{\mathcal H}om({\mathcal F}^{\bullet}, \pi^{!}({\mathcal O}_{S})[-n]\otimes\pi^{*}{\mathcal N}^{\bullet})\ar[d]&\\
{\rm I}\!{\rm R}\pi_{*}{\rm I}\!{\rm R}{\mathcal H}om({\mathcal F}^{\bullet}, \pi^{!}({\mathcal N}^{\bullet})[-n])\eq[r]&{\rm I}\!{\rm R}{\mathcal H}om({\rm I}\!{\rm R}\pi_{*}{\mathcal F}^{\bullet}[n], {\mathcal N}^{\bullet}))}$$
De la bijectivité imposée à la composée des deux premières flèches verticales, on en déduit, alors, pour tout faisceau cohérent ${\mathcal G}$, un isomorphisme
$$\omega^{n}_{X/S}\otimes\pi^{*}{\mathcal G}\simeq\pi^{!}({\mathcal G}[-n]\simeq {\mathcal H}^{-n}(\pi^{!}{\mathcal G})$$ 
Mais, l'égalité des foncteurs $\pi^{\#}({\mathcal G}):=\omega^{n}_{X/S}\otimes\pi^{*}{\mathcal G}$ (naturellement exact à droite) et $\pi^{!}_{{\mathcal K}}({\mathcal G}):={\mathcal H}^{-n}({\mathcal G})$ (exact à gauche  cf \corollaryref{c2}) force leur exactitude. Mais, l'exactitude (à gauche) de  $\pi^{\#}$ est équivalent à la $S$-platitude de $\omega^{n}_{X/S}$ qui, étant de support $X$, entraîne la platitude de $\pi$. D'où
$$\pi^{!}({\mathcal O}_{S})[-n]\otimes\pi^{*}{\mathcal G}\simeq\pi^{!}({\mathcal G})[-n]$$ 
On, alors, nécessairement
$$\omega^{n}_{X/S}\otimes\pi^{*}{\mathcal G}\simeq\pi^{!}({\mathcal O}_{S})[-n]\otimes\pi^{*}{\mathcal G}$$
pour tout faisceau cohérent ${\mathcal G}$. D'où
$${\mathcal H}^{j}(\pi^{!}({\mathcal O}_{S}))=0,\,\forall\,j\not=-n,\,\,\omega^{n}_{X/S}:={\mathcal H}^{-n}(\pi^{!}({\mathcal O}_{S}))\,\, {\rm est \,plat\, sur\,} S$$
c'est-à-dire $\pi$ de Cohen Macaulay$\,\,\blacksquare$
\vspace{2mm}

\noindent
\begin{rem} 
{\bf(i)}  Pour  ${\mathcal F}^{\bullet}={\mathcal O}_{X}$, on a, pour tout faisceau cohérent ${\mathcal G}$ plat sur $S$, 
$$ {\rm I}\!{\rm L}\pi_{\#}({\mathcal O}_{X})\simeq{\rm I}\!{\rm R}\pi_{*}{\rm I}\!{\rm R}{\mathcal H}om( {\mathcal G}, \pi^{!}({\mathcal O}_{S}))\simeq{\rm I}\!{\rm R}{\mathcal H}om( {\rm I}\!{\rm R}\pi_{*}{\mathcal G}, {\mathcal O}_{S})$$
car   $\pi^{!}({\mathcal O}_{S})):={\rm I}\!{\rm R}{\mathcal H}om({\rm I}\!{\rm L}\pi^{*}{\mathcal D}^{\bullet}_{S}, {\mathcal D}^{\bullet}_{X})$. D'où, en particulier, 
$${\rm I}\!{\rm H}^{0}( {\rm I}\!{\rm L}\pi_{\#}({\mathcal O}_{X}[-n])={\rm I}\!{\rm H}^{-n}( {\rm I}\!{\rm L}\pi_{\#}({\mathcal O}_{X})\simeq {\rm I}\!{\rm H}om({\mathcal G}, \omega^{n}_{X/S})$$
Donnant des relations intéressantes pour ${\mathcal G}={\mathcal O}_{X}$ (i.e pour $\pi$ plat),
$${\rm I}\!{\rm H}^{-n}( {\rm I}\!{\rm L}\pi_{\#}({\mathcal O}_{X})\simeq \Gamma(S, \pi_{*}\omega^{n}_{X/S})$$
et pour  ${\mathcal G}=\omega^{n}_{X/S}$,
$${\rm I}\!{\rm H}^{-n}( {\rm I}\!{\rm L}\pi_{\#}({\mathcal O}_{X})\simeq\Gamma(S, \pi_{*}{\mathcal H}om(\omega^{n}_{X/S}, \omega^{n}_{X/S}))$$
{\bf(ii)} L'accouplement de Yoneda ou l'évaluation, on peut exhiber, et ce quelque soit la pondération, un morphisme canonique 
$${\rm I}\!{\rm L}\pi_{\#}({\mathcal A}^{\bullet})\otimes^{{\rm I}\!{\rm L}}{\mathcal B}^{\bullet}\rightarrow {\rm I}\!{\rm R}\pi_{*}({\mathcal A}^{\bullet} \otimes^{{\rm I}\!{\rm L}}\pi^{!}{\mathcal B}^{\bullet})$$
reposant sur l'exitence du morphisme de foncteurs
$${\rm I}\!{\rm L}\pi^{*}\rightarrow{\rm I}\!{\rm L}\pi^{\#}$$
\end{rem}
\vspace{2mm}

\noindent
En prenant la cohomologie de degré $0$, on a, en particulier, ${\rm I}\!{\rm R}^{n}\pi_{*}{\mathcal F}^{\bullet}\simeq \pi_{\#}$ et, donc, $\pi^{!}_{\mathcal K}\simeq\pi^{\#}$. Mais ceci est équivalent au fait que $\pi$ est de Cohen Macaulay d'après le \theoremref{th2}.\vspace{1mm}

\noindent
{\bf(ii)}$\Longrightarrow$ {\bf(i)}:\vspace{1mm}
Comme $\pi$ de Cohen Macaulay équivaut aux annulations ${\mathcal H}^{j}(\pi^{!}{\mathcal O}_{S})=0$ pour tout $j\not=-n$ et que, par conséquent $\pi^{!}{\mathcal O}_{S}\simeq \omega^{n}_{X/S}[n]$. Alors,
$${\rm I}\!{\rm R}{\mathcal H}om({\mathcal F}^{\bullet}, \omega^{n}_{X/S}[n]\otimes^{{\rm I}\!{\rm L}} {\rm I}\!{\rm L}\pi^{*}{\mathcal D}^{\bullet}_{S})\simeq{\rm I}\!{\rm R}{\mathcal H}om({\mathcal F}^{\bullet}, \pi^{!}{\mathcal O}_{S}\otimes^{{\rm I}\!{\rm L}} {\rm I}\!{\rm L}\pi^{*}{\mathcal D}^{\bullet}_{S}) $$
Comme $\pi$ est plat puisque de Cohen Macaulay, le morphisme canonique
$$\pi^{!}{\mathcal O}_{S}\otimes^{{\rm I}\!{\rm L}} {\rm I}\!{\rm L}\pi^{*}({\mathcal D}^{\bullet}_{S})\rightarrow\pi^{!}({\mathcal O}_{S}\otimes^{{\rm I}\!{\rm L}}{\mathcal D}^{\bullet}_{S})$$
est un isomorphisme. Ainsi, $\pi^{!}{\mathcal O}_{S}\otimes^{{\rm I}\!{\rm L}} {\rm I}\!{\rm L}\pi^{*}({\mathcal D}^{\bullet}_{S})\simeq\pi^{!}({\mathcal D}^{\bullet}_{S})={\mathcal D}^{\bullet}_{X}$. On en déduit donc que $${\rm I}\!{\rm R}{\mathcal H}om({\mathcal F}^{\bullet}[n], \omega^{n}_{X/S}[n]\otimes^{{\rm I}\!{\rm L}} {\rm I}\!{\rm L}\pi^{*}{\mathcal D}^{\bullet}_{S})\simeq{\rm I}\!{\rm R}{\mathcal H}om({\mathcal F}^{\bullet}[n],{\mathcal D}^{\bullet}_{X})$$
D'où, par définition et bidualité
$${\rm I}\!{\rm R}\pi_{*}{\mathcal F}^{\bullet}[n]\simeq{\rm I}\!{\rm L}\pi_{\#}({\mathcal F}^{\bullet})\,\blacksquare$$
\section{\color{blue}{Morphisme de dualité totale.}}
 Les résultats précédents nous permettent d'énoncer 

\Th\label{T4}{}  Soit $\pi:X\rightarrow S$ un morphisme propre, plat à fibres de dimension $n$ tel que $\pi^{!}_{{\mathcal K}}=\pi^{\#}$. Alors, le morphisme naturel 
$$\Delta^{j}_{{\mathcal F}, {\mathcal G}}: {\mathcal E}xt_{\pi}^{j}({\mathcal F}, \pi^{\#}({\mathcal G}))\rightarrow {\mathcal H}om({\rm I}\!{\rm R}^{n-j}\pi_{*}{\mathcal F}, {\mathcal G})$$
est un isomorphisme si et seulement si $\pi$ est de Cohen Macaulay.\rm
\begin{proof} En reprenant le diagramme de dualité de Flenner, on a 
$$\xymatrix{{\rm I}\!{\rm R}\pi_{*}{\rm I}\!{\rm R}{\mathcal H}om({\mathcal F}, \omega^{n}_{X/S}\otimes\pi^{*}{\mathcal G})[j]\ar[d]\eq[r]&{\rm I}\!{\rm R}{\mathcal H}om({\rm I}\!{\rm L}\pi_{\#}({\mathcal  F}), {\mathcal G})[j]\ar[d]\\
{\rm I}\!{\rm R}{\mathcal H}om({\rm I}\!{\rm R}^{n-j}\pi_{*}{\mathcal F}, {\mathcal G}))
&{\rm I}\!{\rm R}{\mathcal H}om({\rm I}\!{\rm R}\pi_{*}{\mathcal F}[n], {\mathcal G}))[j]\ar[l]}$$
dans lequel la flèche horizontale du bas est déduite du morphisme canonique  ${\rm I}\!{\rm R}^{n-j}\pi_{*}{\mathcal F}:={\mathcal H}^{0}({\rm I}\!{\rm R}\pi_{*}{\mathcal F}[n-j])\rightarrow {\rm I}\!{\rm R}\pi_{*}{\mathcal F}[n-j]$.
En prenant la cohomologie de degré $0$, on obtient le diagramme commutatif mettant en évidence la flèche désirée
 $$\xymatrix{{\mathcal E}xt_{\pi}^{j}({\mathcal F}, \pi^{\#}({\mathcal G}))\ar@[red][d]_{\Delta^{j}_{{\mathcal F}, {\mathcal G}}}\eq[r]&{\mathcal E}xt^{j}({\rm I}\!{\rm L}\pi_{\#}({\mathcal  F}), {\mathcal G})\ar[d]_{\psi^j_{{\mathcal F}, {\mathcal G}}}\\
{\mathcal H}om({\rm I}\!{\rm R}^{n-j}\pi_{*}{\mathcal F}, {\mathcal G}))&{\mathcal E}xt^{j}({\rm I}\!{\rm R}\pi_{*}{\mathcal F}[n], {\mathcal G})\ar[l]_{\phi^j_{{\mathcal F}, {\mathcal G}}}}$$
Les morphismes $\phi^j_{{\mathcal F}, {\mathcal G}}$ peuvent se voir commes des morphismes latéraux de la suite spectrale
$${\rm E}^{i,j}_{2}:={\mathcal E}xt^{j}({\mathcal H}^{-j}({\rm I}\!{\rm R}\pi_{*}{\mathcal F}[n]), {\mathcal G}))\Longrightarrow\, {\mathcal E}xt^{i+j}({\rm I}\!{\rm R}\pi_{*}{\mathcal F}[n], {\mathcal G})$$
Alors, si $\Delta^{j}_{{\mathcal F}, {\mathcal G}}$ est un isomorphisme pour tout $j$, cela entraine que $\psi^j_{{\mathcal F}, {\mathcal G}}$ est injectif. Mais cela équivaut à sa bijectivité. En effet,  on voit facilement que la cohomologie du triangle distingué
$${\rm I}\!{\rm R}\pi_{*}{\mathcal F}[n]\rightarrow{\rm I}\!{\rm L}\pi_{\#}({\mathcal  F})\rightarrow{\mathcal C}^{\bullet}_{u}[1]$$
donne, pour tout faisceau cohérent ${\mathcal G}$,  
$${\mathcal H}^{k}({\rm I}\!{\rm R}{\mathcal H}om({\mathcal C}^{\bullet}_{u}, {\mathcal G}))=0,\,\forall\,k$$
et, donc, $${\mathcal H}^{k}({\rm I}\!{\rm R}{\mathcal H}om({\mathcal C}^{\bullet}_{u}, {\mathcal G}^{\bullet}))=0$$
pour tout entier $k$ et tout complexe de faisceaux cohérents (ou à cohomologie cohérente). D'où, les foncteurs étant "way-out",  ${\rm I}\!{\rm R}{\mathcal H}om({\mathcal C}^{\bullet}_{u}, {\mathcal G}^{\bullet})=0$. En prenant ${\mathcal G}^{\bullet}={\mathcal D}^{\bullet}_{X}$ et en double dualisant, on a 
${\mathcal C}^{\bullet}_{u}=0$. Par conséquent $${\rm I}\!{\rm R}\pi_{*}{\mathcal F}[n]\simeq{\rm I}\!{\rm L}\pi_{\#}({\mathcal  F})$$
c'est-à-dire $\pi$ est de Cohen Macaulay.\vspace{1mm}

\noindent
On peut retrouver ce résultat, en raisonnant fibre par fibre. En effet, si  $\Delta^{j}_{{\mathcal F}, {\mathcal G}}$ est un isomorphisme, on a, en particulier, pour $S:=\{s\}$,
$${\rm E}xt^{j}({\mathcal F}, \omega^{n}_{X_s})\simeq \big({\rm H}^{n-j}(X_{s}, {\mathcal F})\big)^{'}$$
ce qui montre, en vertu des théorèmes de dualité absolus de Ramis-Ruget, que $\omega^{n}_{X_s}[n]$ est le complexe dualisant de $X_s$. Mais cela équivaut au fait que $X_s$ est de Cohen Macaulay.\vspace{1mm}

\noindent Réciproquement si le morphisme est de Cohen Macaulay, on considère, alors, une paramétrisation locale en un certain point de $X$ donnée par  $\xymatrix{X\ar@/_/[rr]_{\pi}\ar[r]^{f}&Y\ar[r]^{q}&S}$. Alors, 
$${\rm E}xt^{j}(f_{*}{\mathcal O}_{X}, q^{*}{\mathcal G} \otimes_{{\mathcal
O}_{Y}}\Omega^{n}_{Y/S})=0,\,\,\forall\,j>0$$
puisque le morphisme $\pi$ étant de Cohen Macaulay, $f$ est plat et, par suite, $f_{*}{\mathcal O}_{X}$ est localement libre.  De plus, on a
$${\rm E}xt^{j}({\mathcal F}, \bar\pi({\mathcal G}))\simeq {\rm E}xt^{j}(f_{*}{\mathcal F}, q^{*}{\mathcal G} \otimes_{{\mathcal
O}_{Y}}\Omega^{n}_{Y/S})$$
ce qui nous ramène à traiter la bijectivité de $\Delta^{j}_{{\mathcal F}, {\mathcal G}}$ pour une projection $q:Y:=S\times U\rightarrow S$ pour les faisceaux cohérents ${\mathcal F}$ à supports propres sur $S$ (ou $S$-propres). pour montrer que le morphisme  
$$\Delta^{j}_{{\mathcal F}, {\mathcal G}}: {\mathcal E}xt_{q}^{j}({\mathcal F}, q^{\#}({\mathcal G}))\rightarrow {\mathcal H}om({\rm I}\!{\rm R}^{n-j}q_{!}{\mathcal F}, {\mathcal G})$$
est un isomorphisme, il nous faut utiliser des outils propres à la géométrie analytique et de nature topologique que l'on a développer en détails dans [III]. \vspace{1mm}

\noindent
Pour  établir cette dualité, on se ramène au cas absolu en supposant $S$ lisse (qui est  le cas fondamental pour nos constructions) en vue d'appliquer la dualité de Serre et terminer par une formule de type Künneth.  Pour ce faire, on commence par le cas où  ${\mathcal G}$ est localement libre et donc, essentiellement, ${\mathcal G}:={\mathcal O}_{S}$. Le passage au cas absolu se fait en remarquant que l'isomorphisme canonique classique  $\displaystyle{\Omega^{n}_{Y/S}\simeq {\mathcal H}om(q^{*}(\Omega^{r}_{S}), \Omega^{n+r}_{Y})}$ nous permet d'écrire, pour tout ouvert de Stein $S'$, 
$${\rm E}xt^{j}(Y'; {\mathcal F}, \Omega^{n}_{Y/S})\simeq {\rm E}xt^{j}(Y'; {\mathcal F}\otimes q^{*}(\Omega^{r}_{S}), \Omega^{n+r}_{Y})$$
qui, par dualité de Serre, donne 
$$ {\rm Ext}^{j}(Y'; {\mathcal F}\otimes q^{*}(\Omega^{r}_{S}), \Omega^{n+r}_{Y})\simeq{\rm I}\!{\rm H}om_{\Bbb C}({\rm H}^{n-j+r}_{c}(Y',{\mathcal F}\otimes q^{*}(\Omega^{r}_{S})), {\Bbb C})$$
Une formule de décomposition de type Künneth topologique (cf [III])  nous donne
$${\rm H}^{n-j}_{c}(Y', {\mathcal F}){\displaystyle{\widehat{\otimes}_{\Bbb C}}} {\rm H}^{r}_{c}(S', 
 \Omega^{r}_{S})\simeq{\rm H}^{n+r-j}_{c}(Y',{\mathcal F}{\displaystyle{\otimes}_{{\mathcal O}_{Y}}}q^{*}(\Omega^{r}_{S}))$$
 sachant que $${\rm H}^{k}_{c}(S', 
 \Omega^{r}_{S})=0,\,\,\forall\,k\not=r$$
 Pour des raisons expliquées plus tard, on montre successivement:
$$\big({\rm H}^{n-j+r}_{c}(Y',{\mathcal F}\otimes q^{*}(\Omega^{r}_{S})\big)^{'}=\big({\rm H}^{n-j}_{c}(Y', {\mathcal F})\big)^{'}{\displaystyle{\widehat{\otimes}_{\Bbb C}}} \Gamma(S', {\mathcal O}_{S})\simeq {\rm I}\!{\rm H}om_{\Bbb C}({\rm H}^{n-j}_{c}(Y', {\mathcal F}), \Gamma(S', {\mathcal O}_{S}))$$

$${\rm I}\!{\rm H}om_{\Bbb C}({\rm H}^{n-j}_{c}(Y', {\mathcal F}), \Gamma(S', {\mathcal O}_{S}))\simeq{\rm I}\!{\rm H}om_{\Bbb C}(\Gamma(S', {\rm I}\!{\rm R}^{n-j}q_{!}{\mathcal F}), \Gamma(S', {\mathcal O}_{S}))$$

$${\rm I}\!{\rm H}om_{\Bbb C}(\Gamma(S', {\rm I}\!{\rm R}^{n-j}q_{!}{\mathcal F}), \Gamma(S', {\mathcal O}_{S}))\simeq {\rm I}\!{\rm H}om_{\Gamma(S, {\mathcal O}_{S})}(\Gamma(S', {\rm I}\!{\rm R}^{n-j}q_{!}{\mathcal F}), \Gamma(S', {\mathcal O}_{S}))$$
 et, enfin, grâce à la quasi cohérence (analytique)
 $$\Delta^{j}_{{\mathcal F}, {\mathcal O}_{S}}:{\rm E}xt^{j}_{q}({\mathcal F}, \Omega^{n}_{Y/S})\simeq {\mathcal H}om_{{\mathcal O}_{S}}({\rm I}\!{\rm R}^{n-j}q_{!}{\mathcal F}, {\mathcal O}_{S})$$
 et, par suite, la bijectivité de $\Delta^{j}_{{\mathcal F}, {\mathcal G}}$ pour tout faisceau cohérent localement libre ${\mathcal G}$. On en déduit facilement le cas général en raisonnant par récurrence sur la profondeur de ${\mathcal G}$ et, par conséquent, la bijectivité de
 $$\Delta^{j}_{{\mathcal F}, {\mathcal G}}: {\mathcal E}xt_{\pi}^{j}({\mathcal F}, \pi^{\#}({\mathcal G}))\rightarrow {\mathcal H}om({\rm I}\!{\rm R}^{n-j}\pi_{*}{\mathcal F}, {\mathcal G})$$
 dans ce cadre local. Pour globaliser, la nature local du problème sur la base permet de travailler au voisinage d'une fibre $X_{s_{0}}$. Comme le morphisme est propre, on peut recouvrir cette fibre par un nombre fini d'ouverts de Stein et relativement compacts, il suffit, alors, prendre les suites spectrales reliant le local au global pour conclure (ceux sont des méthodes classiques; pour la localisation, on peut reprendre les techniques de Grauert utilisées pour prouver la cohérence des images directes \cite{GR60}). 
\end{proof}
\section{\color{blue}{Une preuve élémentaire de l'isomorphisme de dualité relative.}}
On va montrer une variante un peu plus faible que le \theoremref{thm1}. 
\Th{}{}\label{thm1'} Soient $\pi:X\rightarrow S$ un morphisme propre à fibres de dimension constante $n$ d'espaces anlytiques réduits, des faisceaux cohérents ${\mathcal F}$ et $ {\mathcal G}$ sur $X$ et $S$ respectivement avec ${\mathcal G}$ sans torsion sur $S$. Alors, le morphisme naturel  
$$\Theta^{\pi}_{{\mathcal F}, {\mathcal G}}:{\rm I}\!{\rm H}om(X; {\mathcal F},\pi^{!}_{\mathcal K}{\mathcal G})
\rightarrow{\rm I}\!{\rm H}om(S; {\rm I}\!{\rm R}^{n}\pi_{*}{\mathcal
F}, {\mathcal G})$$
 est un isomorphisme satisfaisant toutes les propriétés fonctorielles annoncées dans le \theoremref{thm1}. \rm\vspace{2mm}

\noindent
Deux options s'offrent à nous, soit: \vspace{1mm}

\indent 
$\bullet$ établir le résultat dans le cas d'un morphisme lisse puis en déduire la cas général ou\vspace{1mm}

\indent 
$\bullet$ établir le résultat dans le cadre général pour ${\mathcal G}$ localement libre et passer au cas d'un faisceau sans torsion.\vspace{1mm}

\noindent 
Nous allons en donner les grandes lignes car la méthode est utilisée dans le cas non propre avec tous les détails. On suit le diagramme suivant

$$\xymatrix{\Theta^{\pi}_{{\mathcal F}, {\mathcal G}}\,{\rm bijective\,pour}\, {\mathcal G}\,{\rm localement\, libre}\ar@{=>}[d] \ar@{=>}[r] &\Theta^{\pi}_{{\mathcal F}, {\mathcal G}}\,{\rm injective\,pour}\,{\mathcal G}\,{\rm sans \,torsion}\ar@{=>}[d]  \\
\Theta^{\pi}_{{\mathcal F}, {\mathcal G}}\,{\rm bijective\,pour}\,{\mathcal G}\,{\rm sans \,torsion}&\Theta^{\pi}_{{\mathcal F}, {\mathcal G}}\,{\rm injective\,pour\,tout}\,{\mathcal G}\,\ar@{=>}[l] }$$
Pour alléger l'écriture et être cohérent avec la notation adoptée dans le cas non propre, on pose, pour tout morphisme $f$ d'espaces complexes réduits et à fibres de dimension $m$
$${\bar{f}}:=f^{!}_{\mathcal K}$$
Nous renvoyons le  lecteur à \cite{K7} pour les proriétés du foncteur $\bar\pi$ et, en particulier, la compatibilité de sa formation avec la composition des morphismes. Plus précisément, si$f:X\rightarrow Y$ et $g:Y\rightarrow S$ sont deux morphismes à fibres de dimension constante, on a $\overline{g\circ f}= \bar{f}\circ \bar g$.
\subsection {{\bf{Quelques propriétés d'invariance du morphisme $\Theta^{\pi}_{{\mathcal F}, {\mathcal G}}$}}}

\lem{}{}\label{Lem0} Soient $\pi:X\rightarrow S$ un morphisme universellement $n$-équidimensionnel d'espaces complexes réduits, $\phi:S\rightarrow S'$ une modification finie et $\pi'$ leur composée. Alors, \vspace{1mm}

\noindent
{\bf(i)} si $\Theta^{\pi'}_{{\mathcal F}, {\mathcal G}'}$ est un isomorphisme pour tout faisceau cohérent ${\mathcal G}'$ sur $S'$,  $\Theta^{\pi}_{{\mathcal F}, {\mathcal G}}$ l'est aussi pour tout faisceau cohérent ${\mathcal G}$ sur $S$,\vspace{1mm}

\noindent
{\bf(ii)} si $\Theta^{\pi}_{{\mathcal F}, {\mathcal G}}$ est un isomorphisme pour tout faisceau cohérent ${\mathcal G}$ sur $S$, $\Theta^{\pi'}_{{\mathcal F}, \phi_{*}{\mathcal G}}$ l'est aussi.\rm

\begin{proof}
\vspace{1mm}

\noindent  {\bf(i)} Supposons que $\Theta^{\pi'}_{{\mathcal F}, {\mathcal G}'}$ soit bijectif pour tout faisceau cohérent ${\mathcal G}'$ sur $S'$ et  considérons  la factorisation 
$$\xymatrix{X\ar[d]_{\pi}\ar[rd]^{\pi'}&\\
S\ar[r]_{\phi}&S'}$$
 On a, alors, le diagramme commutatif   
$$\xymatrix{\pi'_{*}{\mathcal H}om_{X}({\mathcal F}, \bar\pi'({\mathcal G}'))\eq[r]\eq[d]&{\mathcal H}om_{{S'}}({\rm I}\!{\rm R}^{n}\pi'_{*}{\mathcal F}, {\mathcal
G}')\eq[d]\\
\phi_{*}\pi_{*}{\mathcal H}om_{X}({\mathcal F}, \bar\pi(\bar\phi({\mathcal G'})))\eq[r]&{\mathcal  H}om_{{S'}}(\phi_{*}{\rm I}\!{\rm R}^{n}\pi_{*}{\mathcal F}, {\mathcal
G}')}$$
D'où, pour ${\mathcal G'}:=\phi_{*}({\mathcal G})$, l'isomorphisme
$$\phi_{*}\pi_{*}{\mathcal H}om_{X}({\mathcal F}, \bar\pi(\bar\phi(\phi_{*}({\mathcal G}))))\simeq{\mathcal  H}om_{{S'}}(\phi_{*}{\rm I}\!{\rm R}^{n}\pi_{*}{\mathcal F},  \phi_{*}({\mathcal G}))$$
On termine, alors, grâce au fait que:\vspace{1mm}

\indent
{\bf(a)} le morphisme canonique
$$\phi_{*}{\mathcal H}om_{{S}}({\rm I}\!{\rm R}^{n}\pi_{*}{\mathcal F}, {\mathcal
G}))\rightarrow{\mathcal H}om_{{S'}}(\phi_{*}{\rm I}\!{\rm R}^{n}\pi_{*}{\mathcal F}, \phi_{*}({\mathcal
G}))$$
est un isomorphisme puisque $\phi$ est une modification finie (cf {\bf[II], lemme 1.3, p.8}) et dont on déduit l'isomorphisme
$$\pi_{*}{\mathcal H}om_{X}({\mathcal F}, \bar\pi(\bar\phi(\phi_{*}({\mathcal G})))\simeq {\mathcal H}om_{{S}}({\rm I}\!{\rm R}^{n}\pi_{*}{\mathcal F}, {\mathcal G})$$
\indent {\bf(b)} par construction et définition ($\phi$ étant fini), on a 
$$\bar\phi\circ\phi_{*} ={\rm Id}$$
En effet, par définition 
$$\phi_{*}\bar\phi({\mathcal
G}')\simeq{\mathcal H}om(\phi_{*}{\mathcal O}_{S}, {\mathcal
G}')$$
qui, au vu de l'image réciproque ${\mathcal O}_{S'}\rightarrow\phi_{*}{\mathcal O}_{S}$, induit le morphisme de foncteurs $\phi_{*}\bar\phi\rightarrow {\rm Id}$ et l'identification $\bar\phi\circ\phi_{*} ={\rm Id}$ puisque 
$$\phi_{*}\bar\phi(\phi_{*}{\mathcal
G})\simeq{\mathcal H}om(\phi_{*}{\mathcal O}_{S}, \phi_{*}{\mathcal
G})\simeq\phi_{*}{\mathcal H}om({\mathcal O}_{S}, {\mathcal
G})\simeq\phi_{*}{\mathcal
G} $$
qui, par finitude de $\phi$, donne 
$$\bar\phi\phi_{*}{\mathcal G})\simeq{\mathcal G},\,\,\forall\,\,{\mathcal G}\in{\rm Coh}(S)$$
\vspace{1mm}

\noindent
D'où, finalement, la bijectivité attendue du morphisme 
$$\Theta^{\pi}_{{\mathcal F}, {\mathcal G}}:\pi_{*}{\mathcal H}om_{X}({\mathcal F}, \bar\pi({\mathcal G}))\rightarrow{\mathcal H}om_{{S}}({\rm I}\!{\rm R}^{n}\pi_{*}{\mathcal F}, {\mathcal G})$$
ou (globalement) 
$${\rm I}\!{\rm H}om_{X}({\mathcal F}, \bar\pi({\mathcal G}))\simeq{\rm I}\!{\rm  H}om_{{S}}({\rm I}\!{\rm R}^{n}\pi_{*}{\mathcal F}, {\mathcal
G})$$
\vspace{1mm}

\noindent
{\bf(ii)} D'après les relations ci-dessus, on construit facilement $\Theta^{\pi'}_{{\mathcal F}, \phi_{*}{\mathcal G}}$ au travers du diagramme
$$\xymatrix{\pi'_{*}{\mathcal H}om_{X}({\mathcal F}, \bar\pi'({\mathcal G'}))\eq[r]\eq[d]_{\Theta^{\pi'}_{{\mathcal F}, \phi_{*}{\mathcal G}}}&\phi_{*}\pi_{*}{\mathcal H}om_{X}({\mathcal F}, \bar\pi({\mathcal G})))\eq[d]^{\phi_{*}(\Theta^{\pi}_{{\mathcal F}, {\mathcal G}})}\\
{\mathcal  H}om_{{S'}}({\rm I}\!{\rm R}^{n}\pi'_{*}{\mathcal F}, \phi_{*}{\mathcal
G})\eq[r]&\phi_{*}{\mathcal  H}om_{{S'}}({\rm I}\!{\rm R}^{n}\pi_{*}{\mathcal F}, {\mathcal
G})}$$
Si l'on admet l'injectivité de $\Theta^{\pi}_{{\mathcal F}, {\mathcal G}}$ pour tout morphisme universellement équidimensionnel et tout couple de faisceaux cohérents, alors, on déduit de ce qui précède que $\Theta^{\pi'}_{{\mathcal F}, {\mathcal G'}}$ est un isomorphisme pour tout faisceau cohérent ${\mathcal G'}$ sans torsion sur $S'$. Il suffit, pour cela, d'installer ce faisceau dans une suite erxacte courte du type
$$\xymatrix{0\ar[r]&{\mathcal G'}\ar[r]&\phi_{*}{\mathcal G}\ar[r]&{\mathcal K}\ar[r]&0}$$
le faisceau ${\mathcal G}$ étant l'image réciproque de ${\mathcal G'}$ quotientée par la torsion.\vspace{1mm}

\noindent
On conclut, alors, grâce au diagramme commutatif
$$\xymatrix{{\rm I}\!{\rm H}om_{X}({\mathcal F}, \bar\pi'(({\mathcal G'}))\ar[d]_{\Theta^{\pi'}_{{\mathcal F}, {\mathcal G'}}}\ar@{^{(}->}[r]&
{\rm I}\!{\rm H}om_{X}({\mathcal F}, \bar\pi'(\phi_{*}{\mathcal G})))\ar[r]\eq[d]_{\Theta^{\pi'}_{{\mathcal F}, \phi_{*}{\mathcal G})}}&{\rm I}\!{\rm H}om_{X}({\mathcal F}, \bar\pi'({\mathcal K}))\ar@{^{(}->}[d]_{\Theta^{\pi}_{{\mathcal F}, {\mathcal K}}}\\
{\rm I}\!{\rm  H}om_{{S'}}({\rm I}\!{\rm R}^{n}\pi'_{*}{\mathcal F}, {\mathcal
G'}))\ar@{^{(}->}[r]&{\rm I}\!{\rm  H}om_{{S}}({\rm I}\!{\rm R}^{n}\pi'_{*}{\mathcal F}, \phi_{*}{\mathcal
G})\ar[r]&{\rm I}\!{\rm  H}om_{{S'}}({\rm I}\!{\rm R}^{n}\pi'_{*}{\mathcal F}, {\mathcal K})}$$
$\,\,\,\blacksquare$

\end{proof}
\lem{}{}\label{Lem1} Soit $\phi:S\rightarrow S'$ un morphisme fini et surjectif sur $S'$ lisse. Si $\Theta^{\pi'}_{{\mathcal F}, {\mathcal G}'}$ est un isomorphisme pour tout faisceau cohérent ${\mathcal G'}$ sans torsion sur $S'$, alors $\Theta^{\pi}_{{\mathcal F}, {\mathcal G}}$ l'est aussi pour tout faisceau cohérent ${\mathcal G}$ sans torsion sur $S$.\rm
\begin{proof}
On considère la factorisation 
$\xymatrix{X\ar[r]_{\pi}\ar@/^1pc/[rr]^{\pi'}&S\ar[r]_{\phi}&S'}$
 en supposant $\Theta^{\pi'}_{{\mathcal F}, {\mathcal G}'}$ bijective pour tout faisceau cohérent ${\mathcal G}'$ sans torsion sur $S'$. Alors, comme dans le lemme précédent et avec la même construction, on a le diagramme commutatif
$$\xymatrix{\pi'_{*}{\mathcal H}om_{X}({\mathcal F}, \bar\pi'({\mathcal G}'))\eq[r]\eq[d]&{\mathcal H}om_{{S'}}({\rm I}\!{\rm R}^{n}\pi'_{*}{\mathcal F}, {\mathcal
G'})\eq[d]\\
\phi_{*}\pi_{*}{\mathcal H}om_{X}({\mathcal F}, \bar\pi(\bar\phi({\mathcal G'})))\ar[d]_{\gamma}\eq[r]&{\mathcal H}om_{{S'}}(\phi_{*}{\rm I}\!{\rm R}^{n}\pi_{*}{\mathcal F}, {\mathcal
G}')\ar[d]_{\alpha}\ar@/_1pc/[ld]^{u}\\
\phi_{*}{\mathcal H}om_{{S}}({\rm I}\!{\rm R}^{n}\pi_{*}{\mathcal F}, \bar\phi({\mathcal
G}'))\ar@{^{(}->}[r]_{\beta}\ar@/_1pc/[ru]^{u'}&{\mathcal H}om_{{S'}}(\phi_{*}{\rm I}\!{\rm R}^{n}\pi_{*}{\mathcal F}, \phi_{*}\bar\phi({\mathcal
G}'))\ar@/_1pc/[u]_{\alpha'}}$$
dans lequel  $\beta$ est, ici, seulement injectif puisque, $\phi$ étant fini et surjectif, on a, pour tout faisceau cohérent ${\mathcal M}$ sur $S'$, une surjection naturelle
$$\phi^{*}\phi_{*}{\mathcal M}\rightarrow {\mathcal M}$$
de laquelle résulte le morphisme injectif
$${\mathcal H}om({\mathcal M}, {\mathcal G})\rightarrow {\mathcal H}om(\phi^{*}\phi_{*}{\mathcal M}, {\mathcal G})$$
auquel on applique le foncteur exact à droite $\phi_{*}$ et l'on termine par la formule d'adjonction pour avoir la flèche injective
$$\phi_{*}{\mathcal H}om({\mathcal M}, {\mathcal G})\rightarrow {\mathcal H}om(\phi_{*}{\mathcal M}, \phi_{*}{\mathcal G})$$
Comme on l'a dit plus haut, grâce à l'image réciproque ${\mathcal O}_{S'}\rightarrow\phi_{*}{\mathcal O}_{S}$ et la finitude de $\phi$, on a un morphisme
$$\phi_{*}\bar\phi({\mathcal
G}')\rightarrow {\mathcal
G}'$$
De plus, le morphisme trace $\phi_{*}{\mathcal O}_{S}\rightarrow{\mathcal O}_{S'}$ permet de construire 
un morphisme naturel et injectif
$${\mathcal G}'\rightarrow\phi_{*}\bar\phi({\mathcal G}')$$ 
D'où l'existence de $\alpha$ et $\alpha'$, avec injectivité du premier. \vspace{1mm}

\noindent
En fait, grâce au morphisme trace, on peut voir tout faisceau cohérent ${\mathcal
G}'$ sur $S'$ comme facteur direct du faisceau $\phi_{*}\bar\phi({\mathcal
G}')$ et, par conséquent, ${\mathcal H}om_{{S'}}(\phi_{*}{\rm I}\!{\rm R}^{n}\pi_{*}{\mathcal F}, {\mathcal
G}')$ comme facteur direct de ${\mathcal H}om_{{S'}}(\phi_{*}{\rm I}\!{\rm R}^{n}\pi_{*}{\mathcal F}, \phi_{*}\bar\phi({\mathcal G}'))$ et $\alpha'$ étant le morphisme de projection sur l'un des facteurs. \vspace{1mm}

\noindent La commutativité du diagramme nous montre que $u'$ est nécessairement surjectif et $\gamma$ injectif ($u$ aussi). Comme $\gamma$ est un isomorphisme générique et ce, au moins, sur la partie régulière du morphisme en vertu du cas lisse, $u'$ l'est aussi. Mais ${\mathcal G}'$ étant sans torsion, il en est de même de $\bar\phi({\mathcal G}')$ puisque c'est le cas pour $\phi_{*}\bar\phi({\mathcal G}'))$. Alors, comme les deux faisceaux $\phi_{*}{\mathcal H}om_{{S}}({\rm I}\!{\rm R}^{n}\pi_{*}{\mathcal F}, \bar\phi({\mathcal
G}'))$ et ${\mathcal H}om_{{S'}}(\phi_{*}{\rm I}\!{\rm R}^{n}\pi_{*}{\mathcal F}, {\mathcal
G}')$ sont sans torsion sur $S'$, $u'$ est nécessairement injectif donc bijectif. Cela entraine que $\gamma$ l'est aussi.  D'où l'isomorphisme
$$\pi_{*}{\mathcal H}om_{X}({\mathcal F}, \bar\pi(\bar\phi({\mathcal G}')))\simeq {\mathcal H}om_{{S}}({\rm I}\!{\rm R}^{n}\pi_{*}{\mathcal F}, \bar\phi({\mathcal
G}'))$$
et, donc, la bijectivité de $\Psi^{\pi}_{{\mathcal F}, \bar\phi({\mathcal G}')}$. 
$$\Theta^{\pi}_{{\mathcal F}, \bar\phi({\mathcal G}')}:{\rm I}\!{\rm H}om_{X}({\mathcal F}, \bar\pi(\bar\phi({\mathcal G}')))\simeq{\rm I}\!{\rm  H}om_{{S}}({\rm I}\!{\rm R}^{n}\pi_{*}{\mathcal F}, \bar\phi({\mathcal
G}'))$$
Pour tout faisceau cohérent sans torsion ${\mathcal G}'$ sur $S'$.\vspace{1mm}

\noindent
Il nous reste à montrer que c'est encore le cas pour tout faisceau cohérent sans torsion sur $S$. \vspace{1mm}

\noindent Soit  ${\mathcal G}$ un faisceau cohérent sans torsion sur $S$.  Alors,  le morphisme canonique et injectif ${\mathcal G}\rightarrow\bar\phi\phi_{*}({\mathcal
G})$, nous donne suite exacte naturelle
$$0\rightarrow {\mathcal G}\rightarrow\bar\phi\phi_{*}({\mathcal
G})\rightarrow {\mathcal K}\rightarrow 0$$
dont on déduit aisément la suite exacte
$$\xymatrix{0\ar[r]&{\mathcal G}\ar[r]&\bar\phi\phi_{*}({\mathcal
G})\ar[r]_{u}&\bar\phi\phi_{*}({\mathcal
K})}$$
donnant, à son tour,  la suite exacte courte
$$\xymatrix{0\ar[r]&{\mathcal G}\ar[r]&\bar\phi{\phi_{*}}({\mathcal
G})\ar[r]&{\rm Im}u\ar[r]&0}$$
puis le diagramme 
$$\xymatrix{0\ar[r]&\bar\pi({\mathcal G})\ar[r]&\bar\pi(\bar\phi\phi_{*}({\mathcal
G}))\ar[r]&\bar\pi({\rm Im}u)}$$
De laquelle résulte le diagramme commutatif
$$\xymatrix{{\rm I}\!{\rm H}om_{X}({\mathcal F}, \bar\pi({\mathcal G}))\ar[d]_{\Theta^{\pi}_{{\mathcal F}, {\mathcal G}}}\ar@{^{(}->}[r]&
{\rm I}\!{\rm H}om_{X}({\mathcal F}, \bar\pi(\bar\phi(\phi_{*}{\mathcal G})))\ar[r]\eq[d]_{\Theta^{\pi}_{{\mathcal F}, \bar\phi(\phi_{*}{\mathcal G})}}&{\rm I}\!{\rm H}om_{X}({\mathcal F}, \bar\pi(\bar\phi({\rm Im}u))\ar@{^{(}->}[d]_{\Theta^{\pi}_{{\mathcal F}, \bar\phi({\rm Im}u)}}\\
{\rm I}\!{\rm  H}om_{{S}}({\rm I}\!{\rm R}^{n}\pi_{*}{\mathcal F}, {\mathcal
G})\ar@{^{(}->}[r]&{\rm I}\!{\rm  H}om_{{S}}({\rm I}\!{\rm R}^{n}\pi_{*}{\mathcal F}, \bar\phi(\phi_{*}{\mathcal
G}))\ar[r]&{\rm I}\!{\rm  H}om_{{S}}({\rm I}\!{\rm R}^{n}\pi_{*}{\mathcal F}, \bar\phi({\rm Im}u))}$$
qui permet de conclure $\,\blacksquare$ \end{proof}
\vspace{1mm}

\noindent
\lem{}{}\label{Lem2} Si $\Theta^{\pi}_{{\mathcal F}, {\mathcal G}}$ est injectif pour tout faisceau cohérent ${\mathcal G}$ sans torsion alors il l'est pour ${\mathcal G}$ cohérent arbitraire.\rm
\begin{proof}
 On peut procéder soit  en  "filtrant" par les sous faisceaux de torsion (c-à-dire à un dévissage) soit  par récurrence sur la dimension de la base \footnote{On peut, dès lors, procéder par récurrence sur la dimension de $S$ sachant que le cas absolu est recouvert par la dualité de Serre \cite{S} ou celle de Andreotti-Kas \cite{AK} (le cas de la dimension $1$ est particulièrement simple puisqu' un faisceau est sans torsion si et seulement si il est localement libre).}. \vspace{1mm}
 
 \noindent $\bullet$ La méthode du dévissage consiste, en supposant ${\rm Supp}({\mathcal
G})=S$, à démarrer avec la suite exacte courte 
 $$ 0\rightarrow{\mathcal T}\rightarrow{\mathcal G}\rightarrow {\mathcal G}/{\mathcal T}\rightarrow 0$$
 avec ${\mathcal T}$ le sous faisceau de torsion de ${\mathcal G}$. On traite deux cas 
 selon que le morphisme soit lisse ou non.\vspace{1mm}
 
 \indent {{\bf{(i) le foncteur $\pi^{\#}$.}}} Rappelons que c'est le foncteur exact défini par ${\mathcal G}\rightarrow \pi^{*}({\mathcal G})\otimes \Omega^{N}_{X/S}$.
Comme il est exact, toute suite exacte courte donnera, par application de ce foncteur, une suite exacte courte et donc  le diagramme commutatif à lignes exactes
$$\xymatrix{{\rm I}\!{\rm H}om_{{X}}({\mathcal F},
 \pi^{\#}({\mathcal T}))\ar@{^{(}->}[r]^{\alpha}\ar[d]_{\Theta^{\pi}_{{\mathcal F},{\mathcal T}}}&{\rm I}\!{\rm
H}om_{{X}}({\mathcal F}, \pi^{\#}({\mathcal
G}))\ar[r]^{\beta}\ar[d]_{\Theta^{\pi}_{{\mathcal F},{\mathcal G}}}&{\rm I}\!{\rm
H}om_{{X}}({\mathcal F}, \pi^{\#}({\mathcal G}/{\mathcal T}))
\ar@{^{(}->}[d]_{\Theta^{\pi}_{{\mathcal F},{\mathcal G}/{\mathcal T}}}\\
{\rm I}\!{\rm  H}om_{{S}}({\rm I}\!{\rm R}^{N}\pi_{*}{\mathcal F}, {\mathcal
T})\ar@{^{(}->}[r]_{\alpha'}&{\rm I}\!{\rm H}om_{{S}}({\rm I}\!{\rm R}^{N}\pi_{*}{\mathcal
F}, {\mathcal G})\ar[r]_{\beta'}&{\rm I}\!{\rm  H}om_{{S}}({\rm I}\!{\rm
R}^{N}\pi_{*}{\mathcal F}, {\mathcal G}/{\mathcal T})}$$
qui montre que les morphismes $\Theta^{\pi}_{{\mathcal F},{\mathcal T}}$ et $\Theta^{\pi}_{{\mathcal F},{\mathcal G}}$ ont même noyau; ce qui  nous ramène  à l'étude de $\Theta^{\pi}_{{\mathcal F},{\mathcal T}}$ qui est associé à une projection au dessus d'une base de dimension strictement plus petite que celle de $S$.\vspace{1mm}

\noindent
 Si $S_{0}$ est le support de ${\mathcal T}$ et 
 $$\xymatrix{X_0\ar[r]^{\sigma_0}\ar[d]_{\pi_0}&X\ar[d]^{\pi}\\
 S_0\ar[r]_{i_0}&S}$$
 le diagramme cartésien déduit du changement de base $i_{0}:S_{0}\rightarrow S$, on a, en notant ${\mathcal F}_{0}:=\sigma^{*}_{0}{\mathcal F}$ et sachant que le faisceau des formes relatives commutent à tout changement de base,  les isomorphismes évidents
$${\rm I}\!{\rm H}om_{X}({\mathcal F},
 {\pi}^{\#}({\mathcal T}))\simeq {\rm I}\!{\rm H}om_{{X_{0}}}({\mathcal F}_{0}, 
 {\pi^{\#}_{0}}({\mathcal T}))$$ 
et 
$${\rm I}\!{\rm  H}om_{{S}}({\rm I}\!{\rm R}^{N}\pi_{*}{\mathcal F}, {\mathcal
T})\simeq {\rm I}\!{\rm  H}om_{{S_0}}({\rm I}\!{\rm R}^{N}{\pi_{0}}_{*}({\mathcal F}_{0}), {\mathcal T})$$
Si ${\mathcal T}$ est sans torsion, l'injectivité de $\Theta^{\pi_{0}}_{{\mathcal F}_{0},{\mathcal T}}$ entrainera celle de $\theta^{\pi}_{{\mathcal F},{\mathcal G}}$. Si ce n'est pas le cas, on réiterre le processus   par dévissage de la façon suivante:\vspace{1mm}

\noindent On initialise   en posant ${\rm dim}({\rm Supp}({\mathcal T}))=-1$ si  ${\rm
Supp}({\mathcal T})=\emptyset$ et note ${\mathcal T}_{0}={\mathcal T}:={\rm Tors}({\mathcal G})$, 
${\mathcal T}_{k}:={\rm Tors}({\mathcal T}_{k-1})$, $T_{k}:={\rm Supp}({\mathcal T}_{k})$ et ${\mathcal F}_{k}:={\mathcal F}\vert{T_{k}}$. On définit ainsi une suite `` décroissante'' de faisceaux (dans le sens où leurs supports forment une suite décroissante)  stationnaire à partir d'un certain rang $r$ avec  ${\mathcal T}_{k}=\emptyset$ pour tout $k>r$ et l'on raisonne  sur les suites exactes courtes
$$0\rightarrow{\mathcal T}_{k}\rightarrow {\mathcal T}_{k-1}\rightarrow
{{\mathcal T}_{k-1}}/{\mathcal T}_{k}\rightarrow 0$$ en utilisant les
diagrammes de cohomologie, tels que ci-dessus, munis des  morphismes
$\Theta^{\pi}_{{\mathcal F}_{k},{\mathcal T}_{k}}$,  $\Theta^{\pi}_{{\mathcal
F}_{k-1},{\mathcal T}_{k-1}}$, $\Theta^{\pi}_{{\mathcal F}_{k-1}, {{\mathcal
T}_{k-1}}/{\mathcal T}_{k}}$ avec  ${\mathcal T}_{r}$ sans torsion sur $S_{r}$ et, donc, $\Theta^{\pi}_{{\mathcal F}_{r},{\mathcal T}_{r}}$ injectif. On en déduit l'injectivité de  $\Theta^{\pi}_{{\mathcal F}_{r-1},{\mathcal T}_{r-1}}$ puis, en remontant ces données, celle de $\Theta^{\pi}_{{\mathcal F},{\mathcal G}}$.\vspace{1mm}

 \indent {{\bf{(ii) le foncteur $\bar\pi$.}}} Il exact à gauche seulement et la suite exacte courte nous donne la suite exacte
$$\xymatrix{0\ar[r]&{\bar\pi}({\mathcal T})\ar[r]&{\bar\pi}({\mathcal G})\ar[r]_{\alpha}&\bar\pi({\mathcal G}/{\mathcal T})}$$
et la suite exacte courte
$$\xymatrix{0\ar[r]&\bar\pi({\mathcal T})\ar[r]&\bar\pi({\mathcal G})\ar[r]&{\rm Im}\alpha\ar[r]&0}$$
de laquelle est déduit le diagramme commutatif
$$\xymatrix{&&&{\rm I}\!{\rm
H}om_{{X}}({\mathcal F}, {\bar\pi}({\mathcal G}/{\mathcal T}))\ar@{^{(}->}[ddd]^{\Theta^{\pi}_{{\mathcal F},{\mathcal G}/{\mathcal T}}}\\
{\rm I}\!{\rm H}om_{{X}}({\mathcal F},
 {\bar\pi}({\mathcal T}))\ar@{^{(}->}[r]^{\alpha}\ar[d]_{\Theta^{\pi}_{{\mathcal F},{\mathcal T}}}&{\rm I}\!{\rm
H}om_{{X}}({\mathcal F}, {\bar\pi}({\mathcal
G}))\ar[r]^{\beta}\ar[d]_{\Theta^{\pi}_{{\mathcal F},{\mathcal G}}}&{\mathcal A}\ar@{^{(}->}[ru]\ar@{^{(}->}[d]\ar[r]&\dots\\
{\rm I}\!{\rm  H}om_{{S}}({\rm I}\!{\rm R}^{N}\pi_{*}{\mathcal F}, {\mathcal
T})\ar@{^{(}->}[r]_{\alpha'}&{\rm I}\!{\rm H}om_{{S}}({\rm I}\!{\rm R}^{N}\pi_{*}{\mathcal
F}, {\mathcal G})\ar[r]_{\beta'}&{\mathcal B}\ar@{^{(}->}[rd]\ar[r]&\dots\\
&&&{\rm I}\!{\rm  H}om_{{S}}({\rm I}\!{\rm
R}^{N}\pi_{*}{\mathcal F}, {\mathcal G}/{\mathcal T})}$$
dans lequel on a désigné par ${\mathcal A}$ (resp. ${\mathcal B}$) le groupe  ${\rm I}\!{\rm H}om_{{X}}({\mathcal F}, {\rm Im}\alpha)$ (resp. son image par le morphisme injectif ${\Theta^{\pi}_{{\mathcal F},{\mathcal G}/{\mathcal T}}}$). \vspace{1mm}

\noindent Là encore, on voit que l'injectivité de $\Theta^{\pi}_{{\mathcal F},{\mathcal G}}$ est équivalente à celle de $\Theta^{\pi}_{{\mathcal F},{\mathcal T}}$ et que
$${\rm I}\!{\rm H}om_{X}({\mathcal F},
 {\bar\pi}({\mathcal T}))\simeq {\rm I}\!{\rm H}om_{{X_{0}}}({\mathcal F}_{0}, 
 {\bar\pi}_{0}({\mathcal T}))$$ 
$${\rm I}\!{\rm  H}om_{{S}}({\rm I}\!{\rm R}^{N}\pi_{*}{\mathcal F}, {\mathcal
T})\simeq {\rm I}\!{\rm  H}om_{{S_0}}({\rm I}\!{\rm R}^{N}{\pi_{0}}_{*}({\mathcal F}_{0}), {\mathcal T})$$
avec les notations de {\bf(i)}. On reprends alors le même raisonnement que dans le cas précédent pour achever la preuve$\blacksquare$
\end{proof}

On en déduit le
\cor{}{} \label{cor1} Si le \theoremref{thm1'} est vrai pour une base lisse, il l'est pour toute base réduite.\rm\vspace{1mm}

\noindent 
\subsection{\color{blue}{Le \theoremref{thm1'} déduit du cas des morphismes lisses.}}

\vspace{1mm}

\noindent  Cette approche consiste à utiliser la lissité générique du morphisme et repose fondamentalement sur la 
\Prop{}{}\label{P12} Soit $\pi:X\rightarrow S$ un morphisme lisse, propre à fibres de dimension $N$ sur  une variété analytique lisse de Stein $S$ de dimension finie $r$. Alors, pour tout couple de  faisceaux cohérents $({\mathcal F}, {\mathcal G}) $ de ${\rm Coh}(X)\times {\rm Coh}(S)$ avec ${\mathcal G}$ sans torsion, le morphisme 
$$\Theta^{\pi}_{{\mathcal F}, {\mathcal G}}:{\rm I}\!{\rm H}om(X; {\mathcal F},  \pi^{*}({\mathcal G})\otimes\Omega^{N}_{X/S}))
\rightarrow{\rm I}\!{\rm H}om(S; {\rm I}\!{\rm R}^{N}\pi_{*}{\mathcal
F}, {\mathcal G})$$
est bijectif.\rm
\vspace{1mm}

\noindent Rappelons que l'on note $\pi^{\#}$ le foncteur exact ${\mathcal G}\rightarrow \pi^{*}({\mathcal G})\otimes \Omega^{N}_{X/S}$. Cette proposition est une conséquence presqu'immédiate  de la 
\Prop{}{}\label{P12'} Soient $S$ une variété complexe lisse de dimension $r$ et $Z$ une variété analytique de compacte de dimension  $N$ et $\pi:S\times Z\rightarrow S$ la projection canonique. Alors $\Theta^{\pi}_{{\mathcal F}, {\mathcal G}}$ est un isomorphisme pour tout faisceau cohérent ${\mathcal G}$ sans torsion et ${\mathcal F}$ cohérent quelconque.\rm

\begin{proof}
Le morphisme  $\Theta^{\pi}_{{\mathcal F}, {\mathcal G}}$ apparait comme  la composée des morphismes naturels
$$\xymatrix{ {\rm I}\!{\rm H}om(X; {\mathcal F}, \pi^{\#}({\mathcal G}))\ar@[red][d]_{\Theta^{\pi}_{{\mathcal F}, {\mathcal G}}}\ar[r] &
{\rm I}\!{\rm H}om(S; {\rm I}\!{\rm R}^{N}\pi_{*}{\mathcal F}, {\rm I}\!{\rm R}^{N}\pi_{*}(\pi^{\#}({\mathcal G}))\eq[d]_{\alpha}\\
{\rm I}\!{\rm H}om(S; {\rm I}\!{\rm R}^{N}\pi_{*}{\mathcal F}, {\mathcal G})\eq[r]_{\beta}&{\rm I}\!{\rm H}om(S; {\rm I}\!{\rm R}^{N}\pi_{*}{\mathcal F}, {\mathcal G}\otimes{\rm I}\!{\rm R}^{N}\pi_{*}\Omega^{N}_{X/S})}$$
dans lequel $\alpha$ est induit par l'isomorphisme $${\rm I}\!{\rm R}^{N}\pi_{*}(\pi^{*}({\mathcal G})\otimes\Omega^{N}_{X/S})\simeq{\mathcal G}\otimes{\rm I}\!{\rm R}^{N}\pi_{*}\Omega^{N}_{X/S}$$
qui s'établit facilement à l'aide d'une formule de type Künneth classique:
$$\Gamma(S', {\rm I}\!{\rm R}^{N}\pi_{*}(\pi^{*}({\mathcal G})\otimes\Omega^{N}_{X/S}))={\rm
H}^{N}(S'\times Z, \pi^{*}({\mathcal G})\otimes{p_{S'}}^{\!\! *}\Omega^{N}_{Z})\simeq{\displaystyle{\bigoplus}_{i+j=N}}{\rm
H}^{i}(S', {\mathcal G})\widehat{\otimes_{\Bbb C}}{\rm
H}^{j}(Z, \Omega^{N}_{Z})$$
soit, $S'$ étant de Stein,
$$\Gamma(S', {\rm I}\!{\rm R}^{N}\pi_{*}(\pi^{*}({\mathcal G})\otimes\Omega^{N}_{X/S}))\simeq\Gamma(S', {\mathcal G})\widehat{\otimes_{\Bbb C}}{\rm
H}^{N}(Z, \Omega^{N}_{Z}) $$ 
{\bf(i) $\Theta^{\pi}_{{\mathcal F}, {\mathcal G}}$ est bijectif pour ${\mathcal G}$ localement libre.}\vspace{1mm}

\noindent
pour  ${\mathcal G}$  localement libre, la formule de projection nous donnent  les isomorphismes de ${\mathcal O}_{S}$-modules cohérents 
$$ \pi_{*}{\mathcal H}om({\mathcal F},\Omega^{N}_{X/S})\otimes {\mathcal G}\simeq  \pi_{*}{\mathcal H}om({\mathcal F}, \pi^{\#}({\mathcal G}))$$
$${\mathcal H}om({\rm I}\!{\rm R}^{N}\pi_{*}{\mathcal F}, {\mathcal O}_{S})\otimes {\mathcal G}\simeq {\mathcal H}om({\rm I}\!{\rm R}^{N}\pi_{*}{\mathcal F}, {\mathcal G})$$
qui nous montrent qu'il nous suffit de traiter le cas de ${\mathcal G}={\mathcal O}_{S}$ .\vspace{2mm}

\noindent
Notons  ${\mathcal A}:=\Gamma(S, {\mathcal O}_{S})$ et $X:=S\times Z$. Alors, pour tout ouvert de Stein $S'$ de $S$, on a  la factorisation 
$$\xymatrix{{\rm I}\!{\rm H}om(X'; {\mathcal F}, \Omega^{N}_{X/S})\ar[rr]^{\Phi^{\pi}_{{\mathcal F}, {\mathcal O}_{S}}}\ar[rd]_{\Theta^{\pi}_{{\mathcal F}, {\mathcal O}_{S}}}&&{\rm I}\!{\rm H}om_{\mathcal A}(\Gamma(S', {\rm I}\!{\rm R}^{N}\pi_{*}{\mathcal F}), \Gamma(S', {\mathcal O}_{S}))\\
&{\rm I}\!{\rm H}om(S'; {\rm I}\!{\rm R}^{N}\pi_{*}{\mathcal F}, {\mathcal O}_{S})\ar[ur]_{\chi^{\pi}_{{\mathcal F},{\mathcal O}_{S}}}&}$$
dans laquelle $\chi^{\pi}_{{\mathcal F},{\mathcal O}_{S}}$ est isomorphisme dû à (\cite{Fo}, Satz
2.1) disant que: \vspace{1mm}

\noindent {\it si $T$ est  un espace de Stein, ${\bf Coh}(T)$ la catégorie des faisceaux cohérents
et ${\bf SMod}(A)$ celle des modules de Stein sur l'algèbre de
Stein $A:=\Gamma(T, {\mathcal O}_{T})$, le foncteur section
globale $\Gamma(T, -): {\bf Coh}(T)\rightarrow {\bf SMod}({\mathcal
O}(T))$ est exact et pleinement fidèle. De plus, tout morphisme
$A$-linéaire de $A$-modules de Stein $f:M\rightarrow N$ est
automatiquement continu et $f(M)$ est un sous module fermé de
$N$}.\vspace{1mm}

\noindent
Ainsi, le morphisme canonique
$$\Gamma(T, {\mathcal H}om({\mathcal F}, {\mathcal G}))={\rm I}\!{\rm H}om(T; {\mathcal F}, {\mathcal G})\rightarrow {\rm I}\!{\rm H}om_{\Gamma(T, {\mathcal O}_{T})}(T; \Gamma(T, {\mathcal F}), \Gamma(T, {\mathcal G}))$$
est un isomorphisme.
\vspace{1mm}

\indent
Montrons que $\Phi^{\pi}_{{\mathcal F}, {\mathcal O}_{S}}$ est un isomorphisme de $\Gamma(S, {\mathcal O}_{S})$-modules topologiques.\vspace{1mm}

\noindent
 L'isomorphisme canonique classique  $\displaystyle{\Omega^{N}_{X/S}\simeq {\mathcal H}om(\pi^{*}(\Omega^{r}_{S}), \Omega^{N+r}_{X})}$ nous permet d'écrire
$${\rm I}\!{\rm H}om(X; {\mathcal F}, \Omega^{N}_{X/S})\simeq{\rm I}\!{\rm H}om(X; {\mathcal F}\otimes \pi^{*}(\Omega^{r}_{S}), \Omega^{N+r}_{X})$$
qui, par dualité de Serre, donne 
$${\rm I}\!{\rm H}om(X; {\mathcal F}\otimes \pi^{*}(\Omega^{r}_{S}), \Omega^{N+r}_{X})\simeq{\rm I}\!{\rm H}om_{\Bbb C}({\rm
H}^{N+r}_{c}(X,{\mathcal F}\otimes \pi^{*}(\Omega^{r}_{S})), {\Bbb C})$$
Une formule de décomposition de type Künneth  nous donne la décomposition
$${\rm H}^{N}(X, {\mathcal F}){\displaystyle{\widehat{\otimes}_{\Bbb C}}} {\rm H}^{r}_{c}(S, 
 (\Omega^{r}_{S})\simeq{\rm H}^{N+r}_{c}(X,{\mathcal F}{\displaystyle{\otimes}_{{\mathcal O}_{X}}}\pi^{*}(\Omega^{r}_{S}))$$
 et, par suite, 
$$\xymatrix{{\rm I}\!{\rm H}om(X; {\mathcal F}, \Omega^{N}_{X/S})\eq[ddd]_{\Phi^{\pi}_{{\mathcal F}, {\mathcal O}_{S}}}\eq[r]&{\rm I}\!{\rm H}om(X;{\mathcal F}\otimes \pi^{*}(\Omega^{r}_{S}), \Omega^{N+r}_{X})\eq[r]&{\rm H}om_{\Bbb C}({\rm
H}^{N+r}_{c}(X,{\mathcal F}\otimes \pi^{*}(\Omega^{r}_{S})), {\Bbb C})\eq[d]\\
&&{\rm H}om_{\Bbb C}({\rm
H}^{N}(X,{\mathcal F}), {\rm H}om_{\Bbb C}({\rm H}^{r}_{c}(S, 
 \Omega^{r}_{S}), {\Bbb C}))\eq[d]\\
&& {\rm H}om_{\Bbb C}({\rm
H}^{N}(X,{\mathcal F}), \Gamma(S, {\mathcal O}_{S}))\eq[d]\\
{\rm H}om_{\mathcal A}(\Gamma(S, {\rm I}\!{\rm R}^{N}{\mathcal F}), {\mathcal A})\eq[rr]&&{\rm H}om_{\Bbb C}(\Gamma(S, {\rm I}\!{\rm R}^{N}{\mathcal F}), \Gamma(S, {\mathcal O}_{S}))}$$
On en déduit, alors, que  $\Theta^{\pi}_{{\mathcal F}, {\mathcal O}_{S}}$ est un isomorphisme et que, par conséquent, $\Theta^{\pi}_{{\mathcal F}, {\mathcal G}}$ est un isomorphisme pour tout faisceau cohérent localement libre.\vspace{1mm}

\noindent
{\bf(ii) $\Theta^{\pi}_{{\mathcal F}, {\mathcal G}}$ est bijectif pour tout faisceau cohérent  ${\mathcal G}$ sans torsion sur $S$.}\vspace{1mm}

\noindent Ayant montré la bijectivité de  $\Theta^{\pi}_{{\mathcal F}, {\mathcal G}}$ pour ${\mathcal G}$ localement libre, il nous reste, pour terminer la preuve de notre assertion, à suivre  le schéma: 
$$\xymatrix{\Theta^{\pi}_{{\mathcal F}, {\mathcal G}}\,{\rm injective\,pour}\,{\mathcal G}\,{\rm sans \,torsion}\ar@{==>}[rd]\ar@{=>}[r]&\Theta^{\pi}_{{\mathcal F}, {\mathcal G}}\,{\rm injective\,pour\,tout}\,{\mathcal G}\ar@{=>}[d]\\
&\Theta^{\pi}_{{\mathcal F}, {\mathcal G}}\,{\rm bijective\,pour}\,{\mathcal G}\,{\rm sans \,torsion}}$$
\vspace{1mm}

\indent $\bullet$ {\bf{Injectivité pour ${\mathcal G}$ sans torsion.}}\vspace{1mm}
 
 \noindent
 Comme le problème est de nature local sur $S$, ${\mathcal G}$ peut être vu comme un sous faisceau d'un faisceau localement libre ${\mathcal L}$ de rang fini. L'exactitude du foncteur $\pi^{\#}$ nous donne, alors
 $$\xymatrix{{\rm I}\!{\rm H}om(X; {\mathcal F}, {\pi}^{\#}({\mathcal G}))\ar@{^{(}->}[r]\ar[d]_{\Theta^{\pi}_{{\mathcal F}, {\mathcal G}}}\ar[r]&{\rm I}\!{\rm H}om(X; {\mathcal F}, {\pi}^{\#}({\mathcal L}))\eq[d]_{\Theta^{\pi}_{{\mathcal F}, {\mathcal L}}}\\
{\rm I}\!{\rm H}om(S; {\rm I}\!{\rm R}^{N}\pi_{*}{\mathcal F}, {\mathcal G})\ar@{^{(}->}[r]&{\rm I}\!{\rm H}om(S; {\rm I}\!{\rm R}^{N}\pi_{*}{\mathcal F}, {\mathcal L})}$$
qui, au vu de la bijectivité de $\Theta^{\pi}_{{\mathcal F}, {\mathcal L}}$, donne l'injectivité désirée.  Injectivité que l'on  peut aussi déduire du fait que $\Theta^{\pi}_{{\mathcal F}, {\mathcal G}}$ est un isomorphisme (cf \propositionref{P12}) générique\footnote{il suffit de travailler sur ${\rm Reg}(\pi)\cap {\rm Reg}({\mathcal G})$, ${\rm Reg}({\mathcal G})$ étant le complémentaire du lieu singulier de ${\mathcal G}$.} entre faisceaux cohérents sans torsion. \vspace{1mm}

\indent
$\bullet$ {\bf{Injectivité pour ${\mathcal G}$ cohérent arbitraire.}}\vspace{1mm}
 
 \noindent
  Elle résulte du \lemmaref{Lem2}.\vspace{1mm}
  
\indent
$\bullet$ {\bf{bijectivité pour tout faisceau cohérent  ${\mathcal G}$ sans torsion sur $S$.}}\vspace{1mm}

\noindent
 On considère la suite exacte courte
$$0\rightarrow {\mathcal G}\rightarrow{\mathcal L}\rightarrow {\mathcal K}\rightarrow 0$$
 et  le  diagramme commutatif induit
 $$\xymatrix{0\ar[r]&{\rm I}\!{\rm H}om_{{X}}({\mathcal F},
 {\pi}^{\#}({\mathcal G}))\ar[r]^{\alpha}\ar[d]_{\Theta^{\pi}_{{\mathcal F},{\mathcal G}}}&{\rm I}\!{\rm
H}om_{{X}}({\mathcal F}, {\pi}^{\#}({\mathcal
L}))\ar[r]^{\beta}\eq[d]_{\Theta^{\pi}_{{\mathcal F}, {\mathcal L}}}&{\rm I}\!{\rm
H}om_{{X}}({\mathcal F}, {\pi}^{\#}({\mathcal K}))
\ar@{^{(}->}[d]_{\Theta^{\pi}_{{\mathcal F},{\mathcal K}}}\\
0\ar[r]&{\rm I}\!{\rm  H}om_{{S}}({\rm I}\!{\rm R}^{N}\pi_{*}{\mathcal F}, {\mathcal
G})\ar[r]_{\alpha'}&{\rm I}\!{\rm H}om_{{S}}({\rm I}\!{\rm R}^{N}\pi_{*}{\mathcal
F}, {\mathcal L})\ar[r]_{\beta'}&{\rm I}\!{\rm  H}om_{{S}}({\rm I}\!{\rm
R}^{N}\pi_{*}{\mathcal F}, {\mathcal K})}$$
 dans lequel la bijectivité de $\Theta^{\pi}_{{\mathcal F}, {\mathcal G}}$ découle clairement de la surjectivité de $\Theta^{\pi}_{{\mathcal F}, {\mathcal L}}$  et de l'injectivité de  $\Theta^{\pi}_{{\mathcal F}, {\mathcal K}}\,\blacksquare$ 
\end{proof}
\vspace{1mm}

\noindent
\subsubsection{\bf{La preuve de la \propositionref{P12}}} Elle devient évidente puisqu'il suffit de se ramener au cas d'une base lisse grâce aux résultats d'invariance, puis, le problème étant local sur la base, on se fixe un point $s_0$ et choisit un recouvrement fini par des compacts autour de la fibre $X_{s_{0}}$ et nous ramener à la situation de la \propositionref{P13}. Comme les propriétés fonctorielles ont déjà été prouvées dans le \theoremref{thm1},  il est inutile de les reprendre ici$\blacksquare$
\vspace{2mm}

\noindent
\subsubsection{\bf{La preuve du \theoremref{thm1'}}.} Nous allons procéder en deux étapes en montrant:\vspace{1mm}

\indent {\bf{(i) La bijectivité de $\Theta^{\pi}_{{\mathcal F}, {\mathcal G}}$ pour tout faisceau cohérent ${\mathcal G}$ sans torsion sur $S$ normal.}}\vspace{1mm}

\noindent 
On constate que  pour ${\mathcal G}$ de profondeur au moins deux (resp. sans torsion) sur $S$ normal, $\Theta^{\pi}_{{\mathcal F}, {\mathcal G}}$ est bijectif (resp. injectif). \vspace{1mm}

\noindent
 En effet, d'après \propositionref{P6'}, le faisceau ${\bar\pi}({\mathcal G}):=\pi^{!}_{{\mathcal K}}({\mathcal G})$ est aussi de profondeur au moins deux (resp. sans torsion)  sur $X$. Le morphisme $\pi$ étant équidimensionnel (et ouvert), le faisceau  cohérent $\pi_{*}{\mathcal H}om({\mathcal F}, {\bar\pi}({\mathcal G})) $ est, donc, de profondeur au moins deux (resp. sans torsion)  sur  $S$. Comme le faisceau ${\mathcal H}om( {\rm I}\!{\rm R}^{n}\pi_{*}{\mathcal
F}, {\mathcal G})$ est aussi de profondeur au moins deux (resp. sans torsion) sur $S$ et que $\Theta^{\pi}_{{\mathcal F}, {\mathcal G}}$ est un isomorphisme générique en dehors de la codimension deux, il est, alors, bijectif (resp. injectif\footnote{on peut aussi le voir en plongeant ${\mathcal G}$ dans un faisceau localement libre ou dans son bidual qui sont, tous deux, de profondeur au moins deux puisque $S$ est normal.}) sur $S$ tout entier.
\vspace{1mm}

\noindent Grâce au \lemmaref{Lem2}, on en déduit immédiatement son injectivité pour tout faisceau cohérent ${\mathcal G}$ et, par suite, sa bijectivité pour tout faisceau cohérent ${\mathcal G}$ sans torsion sur $S$. Pour s'en convaincre, il suffit de considérer la suite exacte courte (du plongement dans le bidual)
$$\xymatrix{0\ar[r]&{\mathcal G}\ar[r]&\widetilde{\mathcal G}:=(({\mathcal G})^{*})^{*}\ar[r]_{\alpha}&{\mathcal K}\ar[r]&0}$$
et le diagramme commutatif (cf \lemmaref{Lem2})
$$\xymatrix{&&&{\rm I}\!{\rm
H}om_{{X}}({\mathcal F}, {\bar\pi}({\mathcal K}))\ar@{^{(}->}[ddd]^{\Theta^{\pi}_{{\mathcal F},{\mathcal K}}}\\
{\rm I}\!{\rm H}om_{{X}}({\mathcal F},
 {\bar\pi}({\mathcal G}))\ar@{^{(}->}[r]\ar@{^{(}->}[d]_{\Theta^{\pi}_{{\mathcal F},{\mathcal G}}}&{\rm I}\!{\rm
H}om_{{X}}({\mathcal F}, {\bar\pi}(\widetilde{\mathcal
G}))\ar[r]\eq[d]_{\Theta^{\pi}_{{\mathcal F}, \widetilde{\mathcal G}}}&{\mathcal A}\ar@{^{(}->}[ru]\ar@{^{(}->}[d]\ar[r]&\dots\\
{\rm I}\!{\rm  H}om_{{S}}({\rm I}\!{\rm R}^{N}\pi_{*}{\mathcal F}, {\mathcal
G})\ar@{^{(}->}[r]_{\alpha'}&{\rm I}\!{\rm H}om_{{S}}({\rm I}\!{\rm R}^{N}\pi_{*}{\mathcal
F}, \widetilde{\mathcal G})\ar[r]_{\beta'}&{\mathcal B}\ar@{^{(}->}[rd]\ar[r]&\dots\\
&&&{\rm I}\!{\rm  H}om_{{S}}({\rm I}\!{\rm
R}^{N}\pi_{*}{\mathcal F}, {\mathcal K})}$$
dans lequel on a désigné par ${\mathcal A}$ (resp. ${\mathcal B}$) le groupe  ${\rm I}\!{\rm H}om_{{X}}({\mathcal F}, {\rm Im}\alpha)$ (resp. son image par le morphisme injectif $\Theta^{\pi}_{{\mathcal F},{\mathcal G}}$).

{\bf{(ii) La bijectivité de $\Theta^{\pi}_{{\mathcal F}, {\mathcal G}}$ pour tout faisceau cohérent ${\mathcal G}$ sans torsion sur $S$ réduit quelconque. }}\vspace{1mm}

\noindent
Soit $\nu:\hat{S}\rightarrow S$ la normalisation de $S$. Alors, d'après {\bf(ii)} du \lemmaref{Lem1}, $\Theta^{\pi}_{{\mathcal F}, \nu_{*}(\widehat{\mathcal G})}$ est bijective pour tout faisceau cohérent $\widehat{\mathcal G}$ sans torsion sur $\hat{S}$. On a, alors, pour tout faisceau cohérent ${\mathcal G}$ sans torsion sur $S$, la suite exacte courte
$$\xymatrix{0\ar[r]&{\mathcal G}\ar[r]&\nu_{*}(\widehat{\mathcal G})\ar[r]&{\mathcal K}\ar[r]&0}$$
où $\widehat{\mathcal G}:=(\nu^{*}{\mathcal G})/{\mathcal T}_{\nu}$ la préimage de ${\mathcal G}$ quotient par la torsion et ${\mathcal K}$ un faisceau de torsion. En reprenant le diagramme ci-dessus dans lequel on remplace simplement $\widetilde{\mathcal G}$ par $\nu_{*}(\widehat{\mathcal G})$ et utilisant la bijectivité (sa surjectivité) de $\Theta^{\pi}_{{\mathcal F}, \nu_{*}(\widehat{\mathcal G})}$ et l'injectivité de la troisième flèche verticale, la conclusion s'impose$\blacksquare$

\subsection{\bf{Le \theoremref{thm1'} déduit du cas où ${\mathcal G}$ est 
 localement libre.}}\vspace{1mm}

\noindent Ce point de vue consiste à se focaliser sur le faisceau ${\mathcal G}$ et non sur le morphisme. Pour ce faire, on commence par établir une variante local de la \propositionref{P12'}: 
\Prop{}{}\label{P14} Soient $q:S\times U\rightarrow S$ la projection canonique, ${\mathcal F}$ un faisceau cohérent à support propre sur $S$. Alors, le morphisme
$$\Theta^{q}_{{\mathcal F}, {\mathcal O}_{S}}:{\rm I}\!{\rm H}om(Y; {\mathcal F},\Omega^{n}_{Y/S})
\rightarrow{\rm I}\!{\rm H}om(S; {\rm I}\!{\rm R}^{n}q_{!}{\mathcal
F}, {\mathcal O}_{S})$$
est un isomorphisme.\rm
\begin{proof} L' isomorphisme naturel
$\Omega^{n}_{Y/S}\simeq{\mathcal H}om(q^{*}\Omega^{r}_{S}, \Omega^{n+r}_{Y})$ permet d'écrire
$${\rm I}\!{\rm H}om(Y; {\mathcal F},\Omega^{n}_{Y/S})\simeq{\rm I}\!{\rm H}om(Y; {\mathcal F}\otimes q^{*}\Omega^{r}_{S}, \Omega^{n+r}_{Y})$$
qui, grâce à la dualité de Serre, donne
$${\rm H}om(Y; {\mathcal F}\otimes q^{*}\Omega^{r}_{S}, \Omega^{n+r}_{Y})\simeq \big({\rm H}^{n+r}_{c}(Y, {\mathcal F}\otimes q^{*}\Omega^{r}_{S})\big)^{'}$$
La suite spectrale de Leray pour le morphisme $q$, donne 
$${\rm H}^{n+r}_{c}(Y, {\mathcal F}\otimes q^{*}\Omega^{r}_{S})\simeq {\rm H}^{r}_{c}(S, {\rm I}\!{\rm R}^{n}{q_{!}}({\mathcal F}\otimes q^{*}\Omega^{r}_{S}))$$
Or, par hypothèse sur le support de ${\mathcal F}$, ${\rm I}\!{\rm R}^{n}{q_{!}}({\mathcal F}\otimes q^{*}\Omega^{r}_{S})$ est cohérent et grâce à la formule de projection
$${\rm I}\!{\rm R}^{n}{q_{!}}({\mathcal F}\otimes q^{*}\Omega^{r}_{S})\simeq{\rm I}\!{\rm R}^{n}{q_{!}}({\mathcal F}) \otimes \Omega^{r}_{S} $$
En appliquant, encore une fois la dualité de Serre sur $S$, on obtient
finalement

$${\rm H}om(Y; {\mathcal F}\otimes q^{*}{\Omega^{r}_{S}}, \Omega^{n+r}_{Y})\simeq {\rm H}om(S; {\rm I}\!{\rm R}^{n}{q_{!}}({\mathcal F}) \otimes\Omega^{r}_{S}, \Omega^{r}_{S})$$
c'est-à-dire
$${\rm I}\!{\rm H}om(Y; {\mathcal F},\Omega^{n}_{Y/S})\simeq{\rm H}om(S; {\rm I}\!{\rm R}^{n}{q_{!}}({\mathcal F}), {\mathcal O}_{S} )$$
Remarquons que l'on pouvait aussi montrer l'égalité des foncteurs exacts à droite
${\mathcal T}_{1}({\mathcal G})={\mathcal G}\otimes{\rm I}\!{\rm R}^{N}\pi_{!}\Omega^{N}_{X/S}$ et ${\mathcal T}_{2}({\mathcal G})={\rm I}\!{\rm R}^{N}\pi_{!}(\pi^{*}({\mathcal G})\otimes\Omega^{N}_{X/S})$ en considérant des résolutions injectives de ${\mathcal G}\,\,\blacksquare$.
\end{proof}
\cor{}{}\label{coro2} Pour tout diagramme de factorisation locale d'un morphisme $\pi:X\rightarrow S$ propre  à fibres de dimension constante, 
$$\xymatrix{X\ar[rd]_{\pi}\ar[rr]^{f}&&S\times U:=Y\ar[ld]^{q}\\
&S&}$$
dans lequel $f$ est fini et surjectif et $q$ la projection canonique, le morphisme $$\Theta^{\pi}_{{\mathcal F}, {\mathcal O}_{S}}:{\rm I}\!{\rm H}om(X; {\mathcal F}, \omega^{n}_{\pi})
\rightarrow{\rm I}\!{\rm H}om(S; {\rm I}\!{\rm R}^{n}\pi_{*}{\mathcal
F}, {\mathcal O}_{S})$$
est un isomorphisme.\rm
\begin{proof}
Rappelons que, par construction, le faisceau $\omega^{n}_{\pi}$ (cf {\bf{Théorème 2}} de \cite{K3}) vérifie
$$f_{*}{\mathcal H}om({\mathcal F}, \omega^{n}_{\pi})\simeq{\mathcal H}om(f_{*}{\mathcal F},\Omega^{n}_{Y/S})$$
D'où, en appliquant, la proposition précédente à $f_{*}{\mathcal F}$,
$${\rm I}\!{\rm H}om(Y; f_{*}{\mathcal F},\Omega^{n}_{Y/S})\simeq{\rm H}om(S; {\rm I}\!{\rm R}^{n}{q_{!}}(f_{*}{\mathcal F}), {\mathcal O}_{S} )\simeq{\rm H}om(S; {\rm I}\!{\rm R}^{n}\pi_{*}{\mathcal F}, {\mathcal O}_{S} ) \,\,\blacksquare$$ 
\end{proof}
\cor{}{}\label{coro3} Pour tout morphisme $\pi:X\rightarrow S$ propre, surjectif et à fibres de dimension constante $n$, le morphisme 
$$\Theta^{\pi}_{{\mathcal F}, {\mathcal O}_{S}}:{\rm I}\!{\rm H}om(X; {\mathcal F}, \omega^{n}_{\pi})
\rightarrow{\rm I}\!{\rm H}om(S; {\rm I}\!{\rm R}^{n}\pi_{*}{\mathcal
F}, {\mathcal O}_{S})$$
est bijectif.\rm
\begin{proof}
   C'était un corollaire de (\cite{K3}, {\bf{théorème 1.2}}).  Comme cela a déjà été fait dans la construction du morphisme ${\rm I}\!{\rm R}^{n}\pi_{!}\omega^{n}_{\pi}\rightarrow {\mathcal O}_{S}$ du théorème cité  et sera repris dans le cas non propre, il nous suffit de procéder au découpage du morphisme pour nous ramener au cas local et appliquer ce qui précède$\,\blacksquare$
\end{proof}
\cor{}{}\label{coro3'}Pour tout morphisme $\pi:X\rightarrow S$ propre, surjectif et à fibres de dimension constante $n$, le morphisme 
$$\Theta^{\pi}_{{\mathcal F}, {\mathcal G}}:{\rm I}\!{\rm H}om(X; {\mathcal F}, {\bar\pi}({\mathcal G}))
\rightarrow{\rm I}\!{\rm H}om(S; {\rm I}\!{\rm R}^{n}\pi_{*}{\mathcal
F}, {\mathcal G})$$
est bijectif pour tout faisceau localement libre ${\mathcal G}$.\rm
\vspace{1mm}

\noindent
\subsection{\color{blue}{La preuve du \theoremref{thm1'}}} 

Elle devient claire au vu de ce qui précède. 
  En effet, ${\mathcal G}$ étant sans torsion, on peut l'installer dans une suite exacte courte du type 
$$0\rightarrow {\mathcal G}\rightarrow{\mathcal L}\rightarrow {\mathcal K}\rightarrow 0$$
avec ${\mathcal L}$ localement libre de rang fini.  Par abus de notation, on a le diagramme commutatif induit (cf \lemmaref{Lem2})
 $$\xymatrix{0\ar[r]&{\rm I}\!{\rm H}om_{{X}}({\mathcal F},
 {\bar\pi}({\mathcal G}))\ar[r]^{\alpha}\ar[d]_{\Theta^{\pi}_{{\mathcal F},{\mathcal G}}}&{\rm I}\!{\rm
H}om_{{X}}({\mathcal F}, {\bar\pi}({\mathcal
L}))\ar[r]^{\beta}\eq[d]_{\Theta^{\pi}_{{\mathcal F}, {\mathcal L}}}&{\rm I}\!{\rm
H}om_{{X}}({\mathcal F}, {\bar\pi}({\mathcal K}))
\ar@{^{(}->}[d]_{\Theta^{\pi}_{{\mathcal F},{\mathcal K}}}\\
0\ar[r]&{\rm I}\!{\rm  H}om_{{S}}({\rm I}\!{\rm R}^{N}\pi_{*}{\mathcal F}, {\mathcal
G})\ar[r]_{\alpha'}&{\rm I}\!{\rm H}om_{{S}}({\rm I}\!{\rm R}^{N}\pi_{*}{\mathcal
F}, {\mathcal L})\ar[r]_{\beta'}&{\rm I}\!{\rm  H}om_{{S}}({\rm I}\!{\rm
R}^{N}\pi_{*}{\mathcal F}, {\mathcal K})}$$
 dans lequel la bijectivité de $\Theta^{\pi}_{{\mathcal F}, {\mathcal G}}$ découle clairement de la surjectivité de $\Theta^{\pi}_{{\mathcal F}, {\mathcal L}}$ donnée par le \corollaryref{coro3'}  et de l'injectivité de  $\Theta^{\pi}_{{\mathcal F}, {\mathcal K}}$ donnée par le \lemmaref{lem2}.\vspace{1mm}
 
 \noindent
 On peut remarquer qu'en vertu du \lemmaref{Lem0} ou du \lemmaref{Lem1}, on peut supposer $S$ normal ou même lisse. Dans ce cas, ${\mathcal G}$ est localement libre en dehors de la codimension deux au moins dans $S$. Par conséquent, si ${\mathcal G}$ est de profondeur au moins deux (resp. de torsion), on a la bijectivité (resp. injectivité) de $\Theta^{\pi}_{{\mathcal F}, {\mathcal L}}$ pour les mêmes raisons que celles évoquées dans la preuve utilisant la lissité générique du morphisme. On termine, alors, en invoquant le \lemmaref{Lem2} garantissant l'injectivité pour tout faisceau cohérent ${\mathcal G}$ que l'on exploite pour une suite exacte courte du type ci-dessus dans laquelle on peut remplacer ${\mathcal L}$ par le bidual qui est de profondeur au moins deux$\,\blacksquare$

\end{document}